\documentclass[reqno,a4paper,final,10pt]{amsart}
\usepackage[british]{babel}
\usepackage{amssymb,amstext,amsthm,eucal,amscd}
\usepackage[width=15cm, height=22cm, centering]{geometry}
\usepackage{array}
\usepackage[svgnames,table]{xcolor}
\usepackage[T1]{fontenc}
\usepackage[utf8]{inputenc}
\usepackage{euler}
\usepackage{mathrsfs}
\usepackage{verbatim}
\usepackage[final]{microtype}
\usepackage[nosort]{cite}
\usepackage{adjustbox}
\usepackage{longtable}
\usepackage{pdflscape}
\usepackage{booktabs} 
\usepackage{makecell}
\usepackage{caption}
\usepackage{appendix}
\usepackage{chngcntr}
\usepackage{apptools}
\usepackage{enumitem}
\usepackage{romanbar}
\usepackage{rotating}
\usepackage{tikz,pgfplots}
\usepgfplotslibrary{fillbetween}
\usetikzlibrary{calc,arrows,arrows.meta,cd,shapes.geometric,positioning,through,decorations.markings,decorations.text}
\tikzset{>=latex}
\usepackage[unicode,hidelinks,colorlinks=true,linkcolor=blue,citecolor=blue,pagebackref=true]{hyperref}
\renewcommand*{\backref}[1]{}
\renewcommand*{\backrefalt}[4]{%
  \ifcase #1%
  \or [Page~#2.]%
  \else [Pages~#2.]%
  \fi%
}
\newcommand{\be}{\boldsymbol{e}}
\newcommand{\bsigma}{\boldsymbol{\sigma}}
\newcommand{\J}{\boldsymbol{J}}
\newcommand{\B}{\boldsymbol{B}}
\renewcommand{\P}{\boldsymbol{P}}

\newcommand{\bx}{\boldsymbol{x}}
\newcommand{\by}{\boldsymbol{y}}
\newcommand{\bz}{\boldsymbol{z}}
\newcommand{\gl}{\mathfrak{gl}}
\newcommand{\g}{\mathfrak{g}}
\newcommand{\h}{\mathfrak{h}}
\newcommand{\s}{\mathfrak{s}}

\newcommand{\sS}{{\mathsf{S}}}

\newcommand{\GL}{{\mathsf{GL}}}
\newcommand{\T}{{\mathsf{T}}}
\newcommand{\SL}{{\mathsf{SL}}}
\newcommand{\Diag}{{\mathsf{Diag}}}
\newcommand{\so}{\mathfrak{so}}
\newcommand{\su}{\mathfrak{su}}
\newcommand{\fsl}{\mathfrak{sl}}

\newcommand{\fr}{\mathfrak{r}}
\newcommand{\ft}{\mathfrak{t}}
\newcommand{\fb}{\mathfrak{b}}
\newcommand{\fk}{\mathfrak{k}}

\newcommand{\fdiag}{\mathfrak{diag}}
\newcommand{\CC}{\mathbb{C}}
\newcommand{\PP}{\mathbb{P}}
\newcommand{\RR}{\mathbb{R}}
\newcommand{\ord}{\operatorname{ord}}
\newcommand{\wt}{\operatorname{wt}}
\newcommand{\ad}{\operatorname{ad}}
\newcommand{\id}{\operatorname{id}}
\newcommand{\im}{\operatorname{im}}
\newcommand{\tr}{\operatorname{tr}}
\newcommand{\Tr}{\operatorname{Tr}}
\newcommand{\Ad}{\operatorname{Ad}}
\newcommand{\Hom}{\mathrm{Hom}}
\newcommand{\End}{\mathrm{End}}
\newcommand{\Aut}{\mathrm{Aut}}
\newcommand{\Aff}{\mathsf{Aff}}
\newcommand{\vecspan}{\operatorname{span}}

\renewcommand{\tilde}{\widetilde}
\newcommand{\aff}{\mathfrak{aff}}
\renewcommand{\d}{\partial}
\newcommand{\ba}{\boldsymbol{a}}
\newcommand{\q}{\mathfrak{q}}
\DeclareMathOperator{\starwedge}{\tilde\curlywedge}
\newcommand\ext[1]{\scalebox{.9}[1]{$\bigwedge^{\!#1}$}}
\theoremstyle{plain}
\newtheorem{lemma}{Lemma}[section]
\newtheorem{proposition}[lemma]{Proposition}
\newtheorem{theorem}[lemma]{Theorem}
\newtheorem{corollary}[lemma]{Corollary}
\theoremstyle{definition}
\newtheorem{definition}[lemma]{Definition}
\newtheorem{remark}[lemma]{Remark}
\newtheorem{example}[lemma]{Example}
\newtheorem{thmx}{Theorem}

\newtheorem*{thm}{Theorem}
\allowdisplaybreaks[0]
\numberwithin{equation}{section}
\AtAppendix{\counterwithin{lemma}{section}}
\newsavebox{\pmatrixbox}
\newenvironment{colorpmatrix}
  {\begin{lrbox}{\pmatrixbox}
   \mathsurround=0pt
   $\displaystyle
   \begin{pmatrix}}
  {\end{pmatrix}$%
   \end{lrbox}%
   \usebox{\pmatrixbox}%
   \kern-\wd\pmatrixbox
   \makebox[0pt][l]{$\left(\vphantom{\usebox{\pmatrixbox}}\right.$}%
   \kern\wd\pmatrixbox}

\begin{document}

\title[Coisotropy-one homogeneous spacetimes I]{Coisotropy-one homogeneous spacetimes I:\\ (3,2)-kinematical Lie algebras with $d\geq 3$}

\author[Figueroa-O'Farrill]{José Miguel Figueroa-O'Farrill}
\address[JMF]{Maxwell Institute and School of Mathematics, The University
  of Edinburgh, James Clerk Maxwell Building, Peter Guthrie Tait Road,
  Edinburgh EH9 3FD, Scotland, United Kingdom}
\email[JMF]{\href{mailto:j.m.figueroa@ed.ac.uk}{j.m.figueroa@ed.ac.uk}, ORCID: \href{https://orcid.org/0000-0002-9308-9360}{0000-0002-9308-9360}}

\author[Lorenzo-Naveiro]{Juan Manuel Lorenzo-Naveiro}
\address[JML]{David and Judi Proctor Department of Mathematics, University of Oklahoma, Physical Sciences Center, 601 Elm Avenue, Norman 73019, Oklahoma, United States of America}
\email[JML]{\href{mailto:jnaveiro@ou.edu}{jnaveiro@ou.edu}, ORCID: \href{https://orcid.org/0000-0002-8683-3067}{0000-0002-8683-3067}}

\begin{abstract}
  This is the first of several papers devoted to the classification of
  coisotropy-one homogeneous spacetimes.  It is the first step in
  extending the classification of spatially isotropic homogeneous
  spacetimes of kinematical Lie groups to the much vaster context of
  homogeneous spacetimes that are not spatially isotropic; that is,
  where the rotations no longer act transitively on the space of
  directions in each tangent space.  In this first paper we classify
  (up to isomorphisms preserving the rotational subalgebra) a
  class of the relevant transitive Lie algebras, the so-called
  (3,2)-kinematical Lie algebras, for the coisotropy-one case in
  spatial dimension at least $4$.  As a by-product of this
  classification, we classify two-dimensional real representations of
  Bianchi Lie algebras up to a suitable notion of equivalence.
\end{abstract}



\maketitle
\tableofcontents

\section{Introduction}
\label{sec:introduction}

More than half a century ago, Bacry and Lévy-Leblond
\cite{Bacry:1968zf} started what could be termed the kinematical
Erlanger programme, in effect adding ``time'' to Klein's famous
programme of studying geometries via their symmetry groups. In their
original paper they classified kinematical Lie algebras for a class of
four-dimensional spacetimes. The corresponding kinematical Lie groups
act on the spacetimes as a composition of rotations, boosts and
spatial and temporal translations. Their classification was subject to
the assumptions of spatial isotropy and the existence of time-reversal
and parity automorphisms. These latter ``by no means compelling''
\cite{Bacry:1968zf} assumptions were lifted in work of Bacry and Nuyts
\cite{Bacry:1986pm}, who still worked in four spacetime dimensions and
still subject to spatial isotropy. More recently, this classification
was extended to arbitrary dimension via deformation theory:
\cite{Figueroa-OFarrill:2017ycu} (recovering the classification of
Bacry and Nuyts), \cite{Figueroa-OFarrill:2017tcy} (in higher
spacetime dimension) and \cite{Andrzejewski:2018gmz} (in three
spacetime dimensions), but always under the assumption of spatial
isotropy. Based on these results, the simply-connected spatially
isotropic homogeneous spacetimes (i.e., the kinematical Klein
geometries) were classified in \cite{Figueroa-OFarrill:2018ilb},
resulting in a small class of geometries, including lorentzian,
galilean and carrollian, which were further studied in
\cite{Figueroa-OFarrill:2019sex}.

A natural continuation of this programme is to lift the assumption of
spatial isotropy and the natural first step is to consider the case of
coisotropy-one: where the rotational subalgebra acts on the tangent
sphere (i.e., the space of directions) with codimension-one
orbits. This is a much richer class of spacetime geometries and their
classification, which is substantially harder, is the goal of a series
of papers commencing with this one.

We recall that a $(d+2)$-dimensional spatially isotropic homogeneous
spacetime admits a transitive action of a spatially isotropic
kinematical Lie group $\mathsf{K}$. Its Lie algebra $\fk$ is a real
Lie algebra of dimension $(d+2)(d+3)/2$ with a Lie subalgebra
$\fr \cong \so(d+1)$ under which the adjoint representation of $\fk$,
decomposes into $\ad \oplus 2 U \oplus \RR$, where $\ad$,
$U \cong \RR^{d+1}$ and $\RR$ are the adjoint, vector and scalar
representations of $\fr$, respectively.  Such Lie algebras were the
subject of the classical papers already mentioned
\cite{Bacry:1968zf,Bacry:1986pm} (for $d=2$) as well as the more
recent papers \cite{Figueroa-OFarrill:2017ycu,
  Figueroa-OFarrill:2017tcy, Andrzejewski:2018gmz} which deal with
general $d$. The Klein pair corresponding to a spatially isotropic
spacetime is $(\fk,\h)$, where $\h \cong \fr \oplus U$ as an
$\fr$-module. Among such Klein pairs, the simplest are those
corresponding to Minkowski, Galilei and Carroll spacetimes. We collect
them in Table~\ref{tab:spacetimes}. The notation in the table is such
that $\fk$ is spanned by $J_{ab},B_a,P_a,H$ and $\h$ is spanned by
$J_{ab},B_a$, with $J_{ab}$ a basis for $\fr$. Geometrically we may
interpret $J_{ab}$ as generating rotations, $B_a$ as generating
boosts, $P_a$ generating spatial translations and $H$ generating time
translations. The table omits the indices and displays the Lie
brackets of $\fk$ in an abbreviated -- yet unambiguous -- notation
explained in \cite[§2.2]{Figueroa-OFarrill:2018ilb}. We omit the Lie
brackets common to all kinematical Lie algebras; that is, those that
involve the rotational subalgebra $\fr$, which in this abbreviated
notation are $[\J,\J] = \J$, $[\J, \B] = \B$, $[\J,\P] = \P$ and
$[\J,H]=0$.

\begin{table}[h!]
  \centering
  \caption{Minkowski, Galilei and Carroll spacetimes.}
  \label{tab:spacetimes}
  \rowcolors{2}{blue!10}{white}
    \begin{tabular}{l|*{3}{>{$}l<{$}}}\toprule
     \multicolumn{1}{c|}{Spacetime} & \multicolumn{3}{c}{Nonzero Lie brackets}\\ \midrule
      Minkowski & [H,\B] = -\P & [\B,\B] = \J & [\B,\P] = H \\
      Galilei & [H,\B] = -\P & &\\
      Carroll & & & [\B,\P] = H \\\bottomrule
    \end{tabular}
  \end{table}

  The minimal breaking of spatial isotropy (where the rotational
  algebra $\so(d+1)$ breaks to $\so(d)$) in these spacetimes will
  result in a coisotropy-one homogeneous spacetime whose transitive
  Lie algebra belongs to one of two classes of Lie algebras we now
  describe.  To this end it is convenient to recall a definition due
  to Morand \cite{Morand:2023emw}.

\begin{definition}
  An \textbf{$(s,v)$-kinematical Lie algebra} (KLA) with
  $n$-dimensional spatial isotropy is a real Lie algebra $\g$ with a
  Lie subalgebra $\fr \cong \so(n)$ and relative to which $\g$
  decomposes as an $\fr$-module as $\g = \fr \oplus v \RR^n \oplus s
  \RR$, where $\RR^n$ and $\RR$ denote, respectively, the standard
  (vector) and scalar irreducible representations of $\so(n)$.
\end{definition}

This definition encompasses many classes of Lie algebras relevant to
the description of spacetimes.  For example, $(1,1)$-KLAs are
aristotelian Lie algebras, whereas $(1,2)$-KLAs are the original
spatially-isotropic kinematical Lie algebras.  The $(2,1)$-KLAs are
the Lifshitz Lie algebras, classified in
\cite{Figueroa-OFarrill:2022kcd}, whereas $(2,2)$-KLAs are the
boost-extended Lifshitz Lie algebras, which are the subject of ongoing
work which will be presented in \cite{BoostedLifshitz}.  The following
structural result will be shown in that paper.

\begin{thm}[\hspace{-.1pt}\cite{BoostedLifshitz}]
  Every $(2,2)$-KLA $\g$ is isomorphic to an extension-by-derivation
  of a $(1,2)$-KLA $\fk$; that is, it fits in a (split) exact sequence
  of the form
  \begin{equation}
    \begin{tikzcd}
      0 \arrow[r] & \fk \arrow[r] & \g \arrow[r] & \RR D \arrow[r] & 0.
    \end{tikzcd}
  \end{equation}
  The $(2,2)$-KLAs obtained by extending $\fk$ by $D$ and by $D'=
  \theta \circ (\lambda D - \mu \ad_H) \circ \theta^{-1}$, where
  $\theta \in \Aut(\fk)$ is an automorphism and $\lambda,\mu \in \RR$
  with $\lambda \neq 0$, are isomorphic.
\end{thm}

This theorem is analogous to a theorem proved in
\cite{Figueroa-OFarrill:2022kcd}, which shows that every $(2,1)$-KLA
is an extension-by-grading-element of a $(1,1)$-KLA.

Let us now break spatial isotropy in any of the spacetimes in
Table~\ref{tab:spacetimes}.  We assume that there is a distinguished
nonzero spatial vector $u \in \fk/\h$, which breaks the rotational
symmetry from $\fr$ to the stabiliser of $u$, which is a subalgebra
$\s = \so(u^\perp)$ isomorphic to $\so(d)$. Under $\s$, the vector
representation decomposes as $U = u^\perp \oplus \RR u$, where
$u^\perp \cong \RR^d$ is the vector representation of $\s$.  We define
$\g$ to be the normaliser of $u \cdot \P$ in $\fk$ and consider the
Klein pair $(\g, \g \cap \h)$ as a model of a coisotropy-one
homogeneous spacetime.  For Minkowski and Carroll spacetimes, $\g \cap
\h$ does not contain the boost $u \cdot \B$ along $u$, since $[u \cdot
\B, u \cdot \P] = |u|^2 H$; whereas for Galilei spacetime $\g \cap \h$
does include $u \cdot \B$.  This gives rise to two classes of Lie
algebras $\g$.

If $\g$ does not include the boost $u \cdot \B$, then as a
representation of $\s = \so(u^\perp)$, $\g$ decomposes as $\g = \s
\oplus 2 u^\perp \oplus 2 \RR$, where the two scalars are the
translation generators along time and along $u$.  Therefore $\g$ is a
$(2,2)$-KLA.

If $\g$ does include the boost $u \cdot \B$, then it has an additional
$\so(u^\perp)$ scalar, namely $u \cdot \B$ itself.  In that case,
$\g = \s \oplus 2 u^\perp \oplus 3 \RR$ as a representation of $\s$
and hence it is a $(3,2)$-KLA, a class of Lie algebras whose
classification is the main result of this paper (for $d\geq 3$) and a
following paper (for $d=2$).  A further paper will contain the
classification of the corresponding simply-connected coisotropy-one
homogeneous spacetimes; that is, the corresponding spacetime Klein
geometries.

Let us now summarise the paper and explain its organisation.  In
Section~\ref{sec:struc-32-klas} we prove several structure theorems
for $(3,2)$-KLAs.  We consider $d \geq 3$.  For $d=1$ there is no
distinction between the scalar and vector representations, so any
five-dimensional Lie algebra is a $(3,2)$-KLA.  Such Lie algebras were
classified up to isomorphism in \cite{MR155871} and we shall have
nothing more to say about them.  For $d=2$ there is no distinction
between the scalar and adjoint representations and this, together with
the fact that the vector representation admits a larger endomorphism
ring than when $d\geq 3$, makes the classification problem
substantially more complicated.  In order to keep the papers to a
reasonable size, we will present the $d=2$ classification in a
separate paper.

In this first paper, our treatment breaks naturally into two cases and
this will be reflected throughout: the generic case ($d>3$)
and the special case $d=3$.  The main results in
Section~\ref{sec:struc-32-klas} are contained in
Theorems~\ref{th:generic-d-32-structure-theorem} (for $d>3$) and
\ref{thm:3d-32StructureTheorem-simplified} (for $d=3$).

\begin{thmx}\label{th:generic-d-32-structure-theorem}
    Let $\g=\so(V)\oplus (V\otimes W)\oplus \fb$ be a $(3,2)$-KLA with $\dim V=d> 3$.
    Then the Lie bracket of $\g$ is determined by the following algebraic data:
    \begin{enumerate}[label=(\roman*), ref=\roman*]
        \setlength{\itemsep}{0pt}
        \item\label{thItem:Bianchi} a (Bianchi) Lie algebra structure on $\fb$,
        \item\label{thItem:2DimRep} a two-dimensional Lie algebra
            representation $\rho \colon \fb\to \gl(W)$,
        \item\label{thItem:PhiForm} and a $\fb$-equivariant $\fb$-valued
            $2$-form $\varphi\colon \ext{2} W \to \fb$.
    \end{enumerate}
    For all $X\in \so(V)$, $v,v_i\in V$, $w,w_i\in W$ and $B\in \fb$,
    the Lie algebra structure is given by the brackets
    \begin{equation}\label{eq:generic-d-lie-bracket}
        \begin{split}
            [X,v\otimes w]={}&Xv\otimes w, \\
            [X,B]={}&0, \\
            [B,v \otimes w]={}&v\otimes \rho(B)w, \\
            [v_1 \otimes w_1, v_2\otimes w_2]={}&\alpha(w_1,w_2) v_1
            \curlywedge v_2 + \left<v_1,v_2\right> \varphi(w_1 \wedge w_2),
        \end{split}
    \end{equation}
    where $\alpha \colon W \times W \to \RR$ is the unique
    $\fb$-invariant symmetric bilinear form given by
    \begin{equation}\label{eq:generic-jacobi-relation}
        \rho(\varphi(w_1\wedge w_2))w_3=\alpha(w_2,w_3)w_1 - \alpha(w_1,w_3)w_2, \quad w_i \in W.
    \end{equation}
\end{thmx}

\begin{thmx}\label{thm:3d-32StructureTheorem-simplified}
    Let $\g=\so(V)\oplus (V\otimes W)\oplus \fb$ be a $(3,2)$-KLA with $\dim V=d= 3$.
    Then the Lie bracket of $\g$ is determined, after perhaps performing a suitable linear isomorphism that leaves $\so(V)$ unchanged, by the following algebraic data:
    \begin{enumerate}[label=(\roman*), ref=\roman*]
        \setlength{\itemsep}{0pt}
        \item a Bianchi Lie algebra structure on $\fb$,
        \item a two-dimensional Lie algebra representation $\rho \colon \fb\to \gl(W)$,
        \item a $\fb$-valued $2$-form $\varphi\colon \ext{2} W \to \fb$ that is $\fb$-equivariant,
        \item and a $\fb$-equivariant $W$-valued symmetric bilinear form
            $\sigma \colon W \times W \to W$.
    \end{enumerate}
    Letting $\alpha \in (\sS^2W^*)^\fb$ be defined by the relation
    \begin{equation}\label{eq:3d-main-alpha-formula}
        \begin{split}
            \alpha(w_2,w_3) w_1 - \alpha(w_1,w_3)w_2 ={}& \rho(\varphi(w_1\wedge w_2)) w_3  \\
                                                        &+ \sigma(\sigma(w_2,w_3),w_1) - \sigma(\sigma(w_1,w_3),w_2),
        \end{split}
    \end{equation}
    we have that $\alpha$, $\sigma$ and $\varphi$
    must satisfy the following equations:
    \begin{equation}\label{eq:3d-algebraic-data-relations}
        \begin{split}
            0 &= \alpha(\sigma(w_1,w_3),w_2) - \alpha(\sigma(w_2,w_3),w_1), \\
            0 &= \varphi(\sigma(w_1,w_2) \wedge w_3 + \sigma(w_2,w_3) \wedge w_1 +  \sigma(w_3,w_1) \wedge w_2).
        \end{split}
    \end{equation}
    for all $w_i \in W$.  The Lie algebra structure is given by the brackets
    \begin{equation}\label{eq:3d-lie-bracket}
        \begin{split}
            [X,v\otimes w] ={}& Xv\otimes w, \\
            [X,B] ={}& 0, \\
            [B,v \otimes w]  ={}& v \otimes \rho(B)w, \\
            [v_1\otimes w_1, v_2\otimes w_2] ={}& \alpha(w_1,w_2) v_1 \curlywedge v_2 + \left<v_1, v_2\right> \varphi(w_1 \wedge w_2) \\
                                                &+ (v_1 \times v_2) \otimes \sigma(w_1, w_2),
        \end{split}
    \end{equation}
    for all $X\in \so(V)$, $v,v_i\in V$, $w,w_i\in W$ and $B\in \fb$.
\end{thmx}

Let us denote by $\Phi=(\fb,\rho,\alpha,\varphi,\sigma)$ (with the
tacit understanding that not all objects appear for all $\dim V$) the
data defining a $(3,2)$-KLA $\g = \g(\Phi)$. In
Section~\ref{sec:relative-isomorphisms} we investigate how two
$5$-tuples $\Phi$ and $\Phi'$ resulting in isomorphic $(3,2)$-KLAs
$\g(\Phi)$ and $\g(\Phi')$ are related. We fix the Lie subalgebra
$\so(V)$ of $\g$ once and for all and consider isomorphisms
$\g(\Phi) \to \g(\Phi')$ which restrict to the identity on $\so(V)$.
These isomorphisms will be known as relative isomorphisms between
$\g(\Phi)$ and $\g(\Phi')$.

For $\dim V > 3$ the main results are contained in
Proposition~\ref{prop:generic-d-relative-isomorphisms}, which states 
that any relative isomorphism $f \colon \g(\Phi) \to \g(\Phi')$ satisfies
\begin{equation}
  f(X)=X, \qquad f(B) = \psi(B) \qquad\text{and}\qquad f(v \otimes w)
  = v \otimes h(w),
\end{equation}
for $X \in \so(V)$, $B \in \fb$ and $v\otimes w \in V \otimes W$,
where $\psi \colon \fb \to \fb'$ is a Lie algebra isomorphism and $h
\colon W \to W'$ is a vector space isomorphism relating $\Phi$ and
$\Phi'$ as follows
\begin{equation}
  \rho' \circ \psi = \Ad(h) \circ \rho, \qquad \alpha = h^*\alpha',
  \qquad\text{and}\qquad \psi \circ \varphi = h^* \varphi',
\end{equation}
where $(\Ad(h) \circ \rho)(B) := h \circ \rho(B) \circ h^{-1}$ for all
$B \in \fb$.

The results for $\dim V = 3$ are contained in
Proposition~\ref{prop:3d-relative-isomorphisms}, which states that the
isomorphism $f \colon \g(\Phi) \to \g(\Phi')$ is given by
\begin{equation}
  f(X)=X, \qquad f(B) = \psi(B) \qquad\text{and}\qquad f(v \otimes w)
  = \mu(w) \varepsilon(v) + v \otimes h(w),
\end{equation}
for $X \in \so(V)$, $B \in \fb$ and $v\otimes w \in V \otimes W$,
where $\varepsilon \colon V \to \so(V)$ is given by $\varepsilon(v) v'  = v
\times v'$.  Here again $\psi \colon \fb \to \fb'$ is a Lie algebra
isomorphism and $h  \colon W \to W'$ is a vector space isomorphism,
but we have a new ingredient: a $\fb$-invariant $\mu \in W^*$.  Now
$\Phi$ and $\Phi'$ are related as follows:
\begin{equation}
  \rho' \circ \psi = \Ad(h) \circ \rho, \quad \psi \circ \varphi =
  h^* \varphi', \quad \alpha - \mu \circ \sigma + \mu \otimes \mu =
  h^*\alpha' \quad\text{and}\quad h \circ \sigma - \mu \otimes h - h
  \otimes \mu = h^* \sigma'.
\end{equation}

A common ingredient in all $(3,2)$-KLAs for $d\geq 3$ is a
two-dimensional real representation $\rho \colon \fb \to \gl(W)$ of a
Bianchi Lie algebra $\fb$. Once we choose a basis for $W$, the image
of $\rho$ is a Lie subalgebra of $\gl(2,\RR)$ and their classification
(up to conjugacy) is reviewed in Section~\ref{sec:gl2-subalgebras}.
This is summarised in Table~\ref{tab:gl2-subalgebras}.  In that
section we also list those subalgebras of $\gl(2,\RR)$ which leave
invariant a linear form $\mu$, a symmetric bilinear form $\alpha$, as
well as those which leave invariant a nonzero symmetric bilinear
$\sigma \colon \RR^2 \times \RR^2 \to \RR^2$.

In Section~\ref{sec:2d-reps-bianchi} we classify
real two-dimensional representations $\rho \colon \fb \to \gl(W)$
of Bianchi Lie algebras. The notion of equivalence which is useful for
the classification of $(3,2)$-KLAs is slightly weaker than the usual one:
two representations $\rho \colon \fb \to \gl(W)$
and $\rho' \colon \fb' \to \gl(W')$ are weakly equivalent if
there exists an isomorphism $h \colon W \to W'$ \emph{and} an
isomorphism $\psi \colon \fb \to \fb'$ such that
$\rho' \circ \psi = \Ad(h) \circ \rho$, where
$\Ad(h) \colon \gl(W) \to \gl(W')$ is defined by
$\Ad(h) f = h \circ f \circ h^{-1}$.  In practice, we fix a Bianchi
Lie algebra $\fb$ in each isomorphism class (as recorded in
Table~\ref{tab:Bianchi}) and we take $\psi \in \Aut(\fb)$ in the
definition of weak equivalence of representations.  For every isomorphism class of
Bianchi Lie algebras, we determine the possible kernels and images of
a representation $\rho$ and determine the possible $\fb$-invariant
$\alpha \in \sS^2 W^*$, $\mu \in W^*$ and $\fb$-equivariant $\varphi
\colon \ext{2} W \to \fb$ and $\sigma \colon \sS^2 W \to W$. These are
summarised in a number of tables in that section.

In Section~\ref{sec:classification-3-2} we bring all ingredients
together and classify the $(3,2)$-KLAs for $d \geq 3$.  The main results of this section are contained in Theorems~\ref{thm:generic-d-32-KLAs-classification} (for $d>3$) and \ref{thm:3d-32-KLAs-classification} (for $d=3$).

\begin{thmx}\label{thm:generic-d-32-KLAs-classification}
  Let $\g=\so(V)\oplus(V\otimes W)\oplus \fb$ be a $(3,2)$-KLA for which $\dim V = d > 3$ and consider the algebraic data $(\fb,\rho,\varphi,\alpha)$ associated with $\g$.
  Then, exactly one of the following situations occurs:
  \begin{enumerate}[label=(\Roman*), ref=\Roman*]
    \setlength{\itemsep}{0pt}
  \item\label{th-item:generic-d-split-abelian-extension}
    The subspace $V\otimes W$ is abelian (that is, $\varphi = 0$).
    In this case, $\g$ is relatively isomorphic to the split abelian extension of $\h=\so(V)\oplus \fb$ by $V\otimes W$ induced by exactly one of the representations given in Section~\ref{sec:2d-reps-bianchi}.
  \item\label{th-item:generic-d-phi}
    The subspace $[V\otimes W,V\otimes W]$ is a nonzero subspace of $\fb$ (that is, $\varphi \neq 0$ but $\alpha = 0$).
    In this case, $\g$ is relatively isomorphic to exactly one of the $(3,2)$-KLAs described in Table~\ref{tab:generic-d-class-II}.
  \item\label{th-item:generic-d-phi-alpha}
    The subspace $[V\otimes W,V\otimes W]$ has nontrivial projection both onto $\so(V)$ and $\fb$ (that is, $\alpha\neq 0$).
    In this case, $\g$ is relatively isomorphic to exactly one of the $(3,2)$-KLAs described in Table~\ref{tab:generic-d-class-III}.
  \end{enumerate}
\end{thmx}

\begin{thmx}\label{thm:3d-32-KLAs-classification}
  Let $\g=\so(V)\oplus(V\otimes W)\oplus \fb$ be a $(3,2)$-KLA for which $\dim V = d = 3$ and consider the algebraic data $(\fb,\rho,\varphi,\alpha,\sigma)$ associated with $\g$.
  Then, up to relative isomorphism, exactly one of the following situations occurs:
  \begin{enumerate}[label=(\Roman*), ref=\Roman*]
    \setlength{\itemsep}{0pt}
  \item\label{th-item:3d-class-I}
    The subspace $V\otimes W$ is abelian (that is, $\varphi$ and $\sigma$ vanish).
    In this case, $\g$ is relatively isomorphic to the split abelian extension of $\h=\so(V)\oplus \fb$ by $V\otimes W$ obtained from exactly one of the representations described in Section~\ref{sec:2d-reps-bianchi}.
  \item\label{th-item:3d-class-II}
    The subspace $[V\otimes W,V\otimes W]$ is a nonzero subspace of $\fb$ (that is, $\sigma$ and $\alpha$ are zero whereas $\varphi \neq 0$).
    In this case, $\g$ is relatively isomorphic to exactly one of the $(3,2)$-KLAs appearing in Table~\ref{tab:generic-d-class-II}.
  \item\label{th-item:3d-class-III}
    The subspace $[V\otimes W,V\otimes W]$ is contained in $\h=\so(V)\oplus \fb$ and has nontrivial projection both onto $\so(V)$ and $\fb$ (that is, $\sigma = 0$ whereas $\alpha \neq 0$).
    In this case, $\g$ is relatively isomorphic to exactly one of the $(3,2)$-KLAs appearing in Table~\ref{tab:generic-d-class-III}.
  \item\label{th-item:3d-class-IV} There is a choice of subspace
    $V\otimes W \subset \g$ complementary to $\so(V) \oplus \fb$ which
    is a nonabelian ideal of $\g$ (equivalently, after a relative
    isomorphism, $\varphi = 0$ whereas $\sigma \neq 0$) and $\g$ is
    not relatively isomorphic to a split abelian extension.  In this
    case, $\g$ is relatively isomorphic to exactly one of the
    $(3,2)$-KLAs appearing in Table~\ref{tab:3d-32-KLAs}.
  \end{enumerate}
\end{thmx}

The paper concludes with two appendices. Appendix~\ref{app:Bianchi}
is devoted to Bianchi Lie algebras and their automorphisms, mostly to
set notation and collect some results which are used in
Section~\ref{sec:2d-reps-bianchi} when we classify two-dimensional
real representations of the Bianchi Lie algebras up to our notion of
equivalence, which uses the Lie algebra automorphisms in an essential
way. In Appendix~\ref{app:binary-cubics} we study a relevant action of
$\GL(2,\RR) \ltimes (\RR^2)^*$ on $\Hom(\sS^{2}\RR^2,\RR^2)$ and
classify its orbits, a problem which we show, for the cases of
interest, to be equivalent to the classification of binary cubic forms
of weight $-1$ under the action of $\GL(2,\RR)$. We then summarise
their classification using classical invariant theory.

\section{Structure of (3,2)-kinematical Lie algebras}
\label{sec:struc-32-klas}

The aim of this section is to prove
Theorems~\ref{th:generic-d-32-structure-theorem}
and~\ref{thm:3d-32StructureTheorem-simplified} stated in the
Introduction. To this effect, we consider for a $d$-dimensional
euclidean vector space $(V,\langle\cdot,\cdot\rangle)$ the
$\so(V)$-representation
\begin{equation}
  \g = \so(V)\oplus V\oplus V \oplus \RR \oplus \RR \oplus \RR,
\end{equation}
where $V$ denotes the standard representation and $\RR$ denotes the
trivial representation.
Our aim is thus to determine all possible Lie brackets on $\g$
that turn it into a $(3,2)$-KLA.
We remark that for $d>3$ (the \emph{generic}
case), the $\so(V)$-modules $\so(V)$, $V$ and $\RR$ are pairwise
non-isomorphic, but if $d=3$ then $\so(V) \cong V$.  We also let
$W=\RR^2$ and $\fb=\RR \oplus \RR\oplus \RR$, so that we may rewrite
our $(3,2)$-KLA as
\begin{equation}
  \g = \so(V) \oplus (V\otimes W) \oplus \fb,
\end{equation}
where under $\ad_{\so(V)}$ we have $V\otimes W$ as the tensor product
of the standard representation and the two-dimensional trivial
representation of $\so(V)$, and $\fb$ is the three-dimensional trivial
representation of $\so(V)$. Note that for $d>3$, the isotypical
components of $\g$ as an $\so(V)$-representation are precisely
$\so(V)$ (if $d=4$, $\so(V)$ splits into two non-isomorphic
irreducible representations), $V\otimes W$ and $\fb$.

We will first treat the generic case and then separately the case
$d=3$, which is somewhat more complicated, just as with
the classification of spatially isotropic kinematical Lie algebras.
We will not consider further the case $d=1$, since any real
five-dimensional Lie algebra fits the definition or $d=2$, which is
the subject of a second paper.

The following notation will prove useful throughout this section:
given a direct sum $X=Y\oplus Z$ of vector spaces, a vector $x\in X$
and a subspace $W\subseteq X$, we denote by $x_{Y}$ and $W_{Y}$ the
projection of $x$ and $W$ into the subspace $Y$ relative to the
previous direct sum decomposition (and we similarly define $x_{Z}$ and
$W_{Z}$).

\subsection{Generic case: $d>3$}
\label{sec:generic-case:-d3}

We now proceed to the proof of
Theorem~\ref{th:generic-d-32-structure-theorem}. Consider an arbitrary
$(3,2)$-KLA $\g$ for which $d>3$. In order to understand the Lie
bracket of $\g$, we need to exploit the Jacobi identity, which tells
us that the bracket $\ext{2} \g \to \g$ is actually a $\g$-module
homomorphism. In particular, it is also an $\so(V)$-equivariant map.
Because the projection maps
\begin{equation}
  \g\to\so(V), \qquad \g\to V\otimes W, \qquad \g \to \fb,
\end{equation}
are all $\so(V)$-equivariant, so are their compositions
\begin{equation}
  \ext{2}\g\to\so(V), \qquad \ext{2}\g\to V\otimes W, \qquad \ext{2}\g \to \fb,
\end{equation}
with the bracket.

Let us consider the restriction of the bracket to $\ext{2} \fb$, which
is a three-dimensional trivial $\so(V)$-module.
Because neither $\so(V)$ nor $V\otimes W$ have trivial subrepresentations
of positive dimension, the
induced maps $\ext{2} \fb\to \so(V)$ and $\ext{2} \fb \to V\otimes W$
are identically zero. This means that the image of the Lie bracket
restricted to $\fb$ lies again in $\fb$, so we have obtained the
following result.

\begin{proposition}
  \label{prop:bianchi}
  If $d>3$, the subspace $\fb$ is a Lie subalgebra of $\g$.
\end{proposition}

In particular, the three-dimensional Lie algebra $\fb$ is a Bianchi Lie
algebra, as these were classified in the first instance by Bianchi
\cite{Bianchi,MR1900159}. Since we will make extensive use of this
classification later on, we recall this classification in
Appendix~\ref{app:Bianchi}, where we choose a representative from each
isomorphism class of Bianchi Lie algebra in Table~\ref{tab:Bianchi}
and also record their automorphism groups.

We now discuss the Lie bracket of elements of $\fb$ with elements of
$V\otimes W$. The Jacobi identity, together with the equality
$[\so(V),\fb]=0$, implies that for every $B\in \fb$ the maps
$\lambda=\operatorname{proj}_{\so(V)}\circ \ad(B)\colon V\otimes W\to
\so(V)$ and
$\xi=\operatorname{proj}_{\fb}\circ \ad (B)\colon V\otimes W\to \fb$
are $\so(V)$-equivariant. If $\lambda \neq 0$, then the restriction of
$\lambda$ to any irreducible invariant subspace of $V\otimes W$ not
contained in $\ker \lambda$ induces an isomorphism $V\cong \so(V)$,
which is not possible for $d>3$, so $\lambda = 0$. Similarly, if
$\xi\neq 0$, its restriction to an irreducible invariant subspace not
contained in $\ker \xi$ would give an injective map $V\to \fb$, which
is not possible for $d>3$, so $\xi = 0$ as well. All in all, we deduce
that $[\fb,V\otimes W]\subseteq V \otimes W$, and we now describe
explicitly this bracket.

For any $B \in \fb$, the map
$\ad(B) \colon V \otimes W \to V \otimes W$ is
$(\so(V)\oplus \fb)$-equivariant.  Since $V$ is an irreducible
$\so(V)$-module, equivariance under $\so(V)$ says\footnote{We may use
  that
  $\ad_B \in \End_{\so(V)}(V \otimes W) = \End_{\so(V)}(V) \otimes
  \End(W)$, but since $V$ is an irreducible representation (whose
  complexification remains irreducible as a complex representation) of
  $\so(V)$, Schur's Lemma says that $\End_{\so(V)}(V) = \RR \id_V$,
  from where the result follows.} that there exists a linear map
$\rho\colon \fb\to\gl(W)$ such that $[B,v\otimes w]=v\otimes\rho(B)w$
for every $v\in V$ and $w\in W$.
Equivariance under $\fb$ then implies that
$\rho$ is a Lie algebra homomorphism.  In summary, we have the
following:

\begin{proposition}
  The Lie bracket $\fb\times (V\otimes W)\to \g$ is given by the equation
  \begin{equation}
    [B,v\otimes w]=v\otimes \rho(B)w,
  \end{equation}
  where $\rho\colon \fb\to\gl(W)$ is a two-dimensional representation of $\fb$.
\end{proposition}
	
The only bracket left for us to understand is that of elements in
$V\otimes W$ with themselves, which defines an $(\so(V) \oplus
\fb)$-equivariant linear map $\ext{2}(V \otimes W) \to \g$.  We may
decompose both spaces into $(\so(V)\oplus \fb)$-modules and hence the
bracket defines a linear map:
\begin{equation}
  \label{eq:vw-bracket}
  \left( \ext{2} V \otimes \sS^2 W\right) \oplus \left( \sS_0^2 V
    \otimes \ext{2} W  \right)  \oplus \left( \RR_{\Tr} \otimes
    \ext{2} W \right) \to \so(V) \oplus \left(V \otimes W\right) \oplus \fb,
\end{equation}
where $\sS^2 V = \sS_0^2 V \oplus \RR_{\Tr}$ is the decomposition into
traceless and trace symmetric $2$-tensors; that is, the dual
$\RR^*_{\Tr}$ is the line spanned by the euclidean inner product in
$\sS^2 V^*$.  We will make use of the $\so(V)$-equivariant isomorphism
$\iota \colon \ext{2} V \to \so(V)$ which sends
$v_1 \wedge v_2 \mapsto v_1 \curlywedge v_2$, where
$v_1 \curlywedge v_2 \in \so(V)$ is the skew-symmetric endomorphism of
$V$ defined by
\begin{equation}
  (v_1 \curlywedge v_2) v = \left<v,v_2\right> v_1- \left<v,v_1\right>
  v_2.
\end{equation}
The $\so(V)$-modules $\ext{2} V$, $V$, $\sS_0^2 V$ and $\RR_{\Tr}$
are mutually non-isomorphic, and hence the only nontrivial components
of the linear map in \eqref{eq:vw-bracket} are
\begin{equation}\label{eq:nontrivial-components-bracket-vw}
  \ext{2} V \otimes \sS^2 W \to \so(V) \qquad\text{and}\qquad
  \RR_{\Tr} \otimes \ext{2} W \to \fb.
\end{equation}
For $d\neq 4$ these modules are also irreducible, but for $d=4$,
$\ext{2} V$ decomposes into two non-isomorphic irreducible
three-dimensional submodules, neither of which are isomorphic to
either $V$, $\sS_0^2 V$ or $\RR_{\Tr}$.  This means that there is a
second $\so(V)$-equivariant isomorphism
$\iota \circ \star \colon \ext{2} V \to \so(V)$ obtained by
precomposing the isomorphism $\iota$ defined above with the Hodge star
operator\footnote{This requires orienting $V$ and, as in $d=3$, we
  assume that a choice of orientation has been made.  The argument is
  manifestly independent of this choice.}
$\star \colon \ext{2} V \to \ext{2} V$.  This sends
$v_1 \wedge v_2 \mapsto v_1 \starwedge v_2$ which, up to an
inconsequential sign, is given by
\begin{equation}
  (v_1 \starwedge v_2) v = \star (v_1 \wedge v_2 \wedge v).
\end{equation}

Taking into account the above isomorphism(s) $\ext{2} V \cong \so(V)$,
the first map in~\eqref{eq:nontrivial-components-bracket-vw} is
determined by some $\alpha \in \sS^2 W^*$ and, for
$d=4$, an additional $\tilde\alpha \in \sS^2W^*$, whereas the second
map is determined by some $\varphi \colon \ext{2} W \to \fb$. In
terms of these maps, the Lie bracket becomes
\begin{equation}
  [v_1 \otimes w_1, v_2 \otimes w_2] = \alpha(w_1, w_2) v_1
  \curlywedge v_2 + \tilde\alpha(w_1,w_2) v_1 \starwedge v_2 + \left<v_1, v_2\right> \varphi(w_1 \wedge w_2),
\end{equation}
where the middle term on the right-hand side is only present for
$d=4$.  Equivariance of the bracket under $\fb$, which is part of the
Jacobi identity, says that $\alpha$ and $\tilde\alpha$ are
$\fb$-invariant and $\varphi$ is $\fb$-equivariant.

A consequence of $\fb$-equivariance of $\varphi$ is that the image of $\rho\circ \varphi$ lies in $\fsl(W)$.
To see this, note that for every $B\in \fb$ and $w_1,\,w_2\in W$, we have
\begin{equation}
  \Tr(\rho(B))\varphi(w_1\wedge w_2)=\varphi(B\cdot (w_1\wedge w_2))=[B,\varphi(w_1\wedge w_2)],
\end{equation}
so applying $\rho$ to both sides we obtain
\begin{equation}
  \Tr(\rho(B))\rho(\varphi(w_1\wedge w_2))=[\rho(B),\rho(\varphi(w_1\wedge w_2))],
\end{equation}
and the right hand side has trace zero because of the identity $\Tr(XY)=\Tr(YX)$, so by letting $B=\varphi(w_1\wedge w_2)$ we deduce that
\begin{equation}
  \Tr(\rho(\varphi(w_1\wedge w_2)))^2 = 0, \quad w_1,w_2\in W,
\end{equation}
which implies that $\im (\rho \circ \varphi)\subseteq \fsl(W)$.
	
Summarising, by imposing the skew-symmetry and $(\so(V)\oplus
\fb)$-equivariance conditions, we have obtained an explicit expression
of the Lie bracket of $\g$ in terms of the Lie algebra structure of
$\fb$, a $2$-dimensional representation $\rho\colon \fb\to \gl(W)$, two
$\fb$-invariant symmetric bilinear forms $\alpha,\, \tilde\alpha \in
\sS^{2}W^{*}$ and a $\fb$-equivariant $\fb$-valued $2$-form
$\varphi\colon \ext{2} W \to \fb$.  We write $\alpha,\,\tilde\alpha
\in (\sS^2W^*)^\fb$ and $\varphi \in \Hom_\fb(\ext{2}W,\fb)$.

In addition, we can determine further relations between these objects by
imposing the full Jacobi identity.  To this end, we define the
\emph{Jacobiator} $J\colon \g\times \g\times \g\to \g$ given by the
equation
\begin{equation}{\label{eq:JacobiatorDefinition}}
  J(x,y,z)=[x,[y,z]]+[y,[z,x]]+[z,[x,y]].
\end{equation}
This is actually a $\g$-valued $3$-form, and for $\g$ to be a Lie
algebra it needs to be identically zero.  Equivariance under $\so(V)
\oplus \fb$ says that $J(x,y,z)=0$ whenever one of the arguments lies
in $\so(V)$ or $\fb$.  Hence we only need to compute $J(x,y,z)$
in the case that $x,\,y,\,z \in V \otimes W$.

We calculate the nested bracket for all $v_i \in V$ and $w_i \in W$:
\begin{multline}
  [[v_1 \otimes w_1, v_2 \otimes w_2], v_3 \otimes w_3]= \alpha(w_1,w_2) \left( \left<v_2,v_3\right> v_1\otimes w_3-
    \left<v_1,v_3\right> v_2\otimes w_3 \right)\\
  + \tilde\alpha(w_1,w_2) \star(v_1 \wedge v_2 \wedge v_3) \otimes
  w_3 + \left<v_1,v_2\right> v_3 \otimes \rho(\varphi(w_1 \wedge w_2)) w_3.
\end{multline}
The Jacobi identity says that
\begin{equation}
  [[v_1 \otimes w_1, v_2 \otimes w_2], v_3 \otimes w_3] +  \text{cyclic} = 0.
\end{equation}
This has to hold for all $v_1,v_2,v_3 \in V$, hence in particular if
we take them to be mutually orthogonal.  In that case, only the
$\tilde\alpha$ terms contribute (only in the case that $d=4$) and the
Jacobi identity becomes
\begin{equation}
  \tilde\alpha(w_1,w_2) w_3 + \text{cyclic} = 0.
\end{equation}
Let $(\be_1,\be_2)$ be a basis for $W$.  Then taking $w_1 = w_2 = w_3
= \be_i$, we find that $\tilde\alpha(\be_1,\be_1) =
\tilde\alpha(\be_2,\be_2) = 0$.  In addition, if we let $w_1=w_2 = \be_1$ and $w_3
= \be_2$ we see that $\tilde\alpha(\be_1,\be_2) = 0$ and thus
$\tilde\alpha = 0$, so the case $d=4$ reduces to the generic case
$d>3$.

From now on we will take
\begin{equation}
  [v_1 \otimes w_1, v_2 \otimes w_2] = \alpha(w_1, w_2) v_1
  \curlywedge v_2 + \left<v_1, v_2\right> \varphi(w_1 \wedge w_2), 
\end{equation}
for $d>3$.  The Jacobi identity is now
\begin{multline}
  [[v_1 \otimes w_1, v_2 \otimes w_2], v_3 \otimes w_3]=
  \alpha(w_1,w_2) \left( \left<v_2,v_3\right> v_1\otimes w_3-
    \left<v_1,v_3\right> v_2\otimes w_3 \right)\\
  + \left<v_1,v_2\right> v_3 \otimes \rho(\varphi(w_1 \wedge w_2)) w_3,
\end{multline}
which, after rearranging, becomes
\begin{equation}
  \left<v_1,v_2\right> v_3 \otimes \left(  \rho(\varphi(w_1 \wedge
    w_2)) w_3 + \alpha(w_1,w_3) w_2 - \alpha(w_2,w_3) w_1\right) +
  \text{cyclic} = 0.
\end{equation}
In particular this has to hold when $v_i$ are such that
$\left<v_1,v_2\right> \neq 0$ and $v_3 \neq 0$, but $\left<v_3,v_1\right> =
\left<v_3,v_2\right> = 0$.  In that case, only the first term
survives, giving the following relation:
\begin{equation}
  \label{eq:phiAlphaRelation}
  \rho(\varphi(w_1 \wedge w_2)) w_3 = \alpha(w_2,w_3) w_1 - \alpha(w_1,w_3) w_2
\end{equation}
for all $w_i \in W$.  Conversely if equation~\eqref{eq:phiAlphaRelation} is
satisfied, then the Jacobi identity is satisfied.

In summary, we have proved that the Jacobiator is identically zero
precisely when $\varphi \in \Hom_{\fb}(\ext{2}W,\fb)$
and $\alpha\in (\sS^2W^*)^\fb$ satisfy
equation~\eqref{eq:phiAlphaRelation}.

It turns out that the relation~\eqref{eq:phiAlphaRelation} completely determines $\alpha$ from $\rho$ and $\varphi$.
Indeed, because $W$ is two-dimensional we have an isomorphism
$\phi\colon \sS^2 W^* \to \Hom(\ext{2} W, \fsl(W))$ of
$\gl(W)$-modules, sending $\alpha \mapsto \phi_\alpha$, where
\begin{equation}\label{eq:iso-s2W*-homwedge2slW}
  \phi_\alpha(w_1\wedge w_2)w_3=\alpha(w_2,w_3)w_1-\alpha(w_1,w_3)w_2.
\end{equation}
In particular, it is also an isomorphism of $\fb$-modules, so it
restricts to a bijective correspondence between $(\sS^2 W^*)^\fb$ and
$\Hom_\fb(\ext{2} W,\fsl(W))$.  As a consequence, if $\varphi\colon
\ext{2}W \to \fb$ is $\fb$-equivariant, the composition $\rho\circ
\varphi$ is in $\Hom_\fb(\ext{2}W,\fsl(W))$,
so~\eqref{eq:iso-s2W*-homwedge2slW} yields a unique
$\alpha\in (\sS^2 W^*)^\fb$ such that~\eqref{eq:phiAlphaRelation} is
satisfied.

\begin{remark}
    Let us look at the expression of the isomorphism~\eqref{eq:iso-s2W*-homwedge2slW} in coordinates.
    Suppose that $\{\be_1,\be_2\}$ is a basis of $W$, and let $\alpha\in \sS^2 W^*$ have coordinates $\alpha_{ij}=\alpha(\be_i,\be_j)$.
    Then the matrix of the endomorphism $\phi_\alpha(\be_1\wedge\be_2)\colon W \to W$ is readily given by
    \begin{equation}
        [\phi_\alpha(\be_1\wedge\be_2)]=\begin{pmatrix} \alpha_{12} & \alpha_{22} \\ -\alpha_{11} & -\alpha_{12} \end{pmatrix}.
    \end{equation}
    As a consequence, given a linear map $\phi\colon \ext{2} W \to \fsl(W)$, if $[\phi^i{}_j]$ is the matrix of $\phi(\be_1\wedge\be_2)$ with respect to our basis, then $\phi=\phi_\alpha$ for the bilinear form $\alpha$ with coefficients
    \begin{equation}
        (\alpha_{ij})=\begin{pmatrix}
            -\phi^2{}_1 & \phi^1{}_1 \\
            \phi^1{}_1 & \phi^1{}_2    
        \end{pmatrix}.
    \end{equation}
    This means that, given $\rho\colon \fb \to \gl(W)$ and $\varphi\in \Hom_\fb(\ext{2}W,\fb)$, we have $\rho \circ \varphi \colon \ext{2} W \to \fsl(W)$ and the unique symmetric form $\alpha\in(\sS^2 W^*)^\fb$ satisfying~\eqref{eq:phiAlphaRelation} is given in coordinates by
    \begin{equation}
        (\alpha_{ij})=\begin{pmatrix}
            -\rho^2{}_1 & \rho^1{}_1 \\
            \rho^1{}_1 & \rho^1{}_2    
        \end{pmatrix},
    \end{equation}
    where $[\rho^i{}_j]$ is the matrix of $\rho(\varphi(\be_1\wedge\be_2))$ with respect to the basis $\{\be_1,\be_2\}$.
\end{remark}

Gathering our calculations, we have shown that the algebra structure on
$\g$ is determined by a Bianchi algebra $\fb$, a Lie algebra
representation $\rho \colon \fb\to \gl(W)$ and a $\fb$-equivariant map
$\varphi\in \Hom_{\fb}(\ext{2} W, \fb)$.
Conversely, if we start with a triple $(\fb,\rho,\varphi)$ satisfying the
aforementioned conditions, it is clear that the bilinear operation on $\g$
defined by~\eqref{eq:generic-d-lie-bracket} induces a $(3,2)$-KLA structure on $\g$.
This finishes the proof of Theorem~\ref{th:generic-d-32-structure-theorem}.

As a consequence of the expression~\eqref{eq:generic-d-lie-bracket}, we also
obtain that $(3,2)$-KLAs with $d>3$ are naturally $\mathbb{Z}_2$-graded:

\begin{corollary}
  Any $(3,2)$-KLA $\g=\so(V)\oplus(V\otimes W)\oplus \fb$ with $\dim V = d > 3$ has a $\mathbb{Z}_2$-gradation given by
  \begin{equation}
    \g_{\bar{0}}=\so(V)\oplus\fb, \qquad \g_{\bar{1}}=V\otimes W.
  \end{equation}
\end{corollary}

\subsection{Algebras with $d=3$}
In this section we prove Theorem~\ref{thm:3d-32StructureTheorem-simplified}.
The main difference when $d=3$ is that now $\ext{2}V \cong V$ as
$\so(V)$-modules. 

Let $\g=\so(V)\oplus(V\otimes W)\oplus \fb$ be a $(3,2)$-KLA with
three-dimensional spatial isotropy.
As a representation of $\so(V)$, we see that $\g$ contains two
isotypical components: the component $\so(V)\oplus (V\otimes W)$
corresponding to the standard representation and the component $\fb$
corresponding to the trivial one.

Note that Proposition~\ref{prop:bianchi} still holds in this case.
Indeed, since $\fb$ is a trivial $\so(V)$-module, so
is $\ext{2}\fb$ and hence, by $\so(V)$-equivariance, the only
possible nonzero bracket is $\ext{2}\fb \to \fb$, making $\fb$ into a
Lie subalgebra of $\g$.  Since $\so(V) \cong V$ as $\so(V)$-modules,
the bracket between $B \in \fb$ and $v\otimes w \in V \otimes W$ can
now have a component along $\so(V)$:
\begin{equation}
  [B, v \otimes w] = v \otimes \rho(B)w + \lambda_B(w) \varepsilon(v),
\end{equation}
for some $\rho \colon \fb \to \gl(W)$ and
$\lambda \colon \fb \to W^*$, sending $B \mapsto \lambda_B$, and where
$\varepsilon(v)\in\so(V)$ is defined by
$\varepsilon(v)u = v \times u$. Note that
$\varepsilon(u \times v) = - u \curlywedge v$.

The $(\fb,\fb,V\otimes W)$ Jacobi identity is tantamount to imposing
$\fb$-equivariance of this bracket.  Let $B,C \in \fb$ and $v\otimes w
\in V \otimes W$.  Then the Jacobi identity reads
\begin{equation}
  [B,[C, v \otimes w]] - [C,[B,v\otimes w]] = [[B,C],v \otimes w]
\end{equation}
and this has two components:
\begin{itemize}
\item a component along $V \otimes W$:
  \begin{equation}
    v \otimes \left( \rho(B)\rho(C) - \rho(C)\rho(B)
      -\rho([B,C])\right)w = 0,
  \end{equation}
  which says that $\rho \colon \fb \to \gl(W)$ is a representation;
\item and a component along $\so(V)$:
  \begin{equation}
    \label{eq:cocycle}
    \left( \lambda_B(\rho(C)w) - \lambda_C(\rho(B)w) -
      \lambda_{[B,C]}(w)\right) \varepsilon(v) = 0.
  \end{equation}
\end{itemize}

\begin{proposition}
  \label{prop:H1}
  Equation~\eqref{eq:cocycle} says that $\lambda \in C^1(\fb;W^*)
  = \Hom(\fb,W^*)$ is a Chevalley--Eilenberg cocycle and hence defines
  a class $[\lambda] \in H^1(\fb;W^*)$.
\end{proposition}

\begin{proof}
  Let $C^p(\fb;W^*) = \Hom(\ext{p}\fb, W^*)$ and let $\d \colon
  C^p(\fb;W^*) \to C^{p+1}(\fb;W^*)$ denote the Chevalley--Eilenberg
  differential.  We will only need the expressions for $p=0,1$.  If
  $\theta \in W^*$, then $\d \theta \in C^1(\fb;W^*)$ is given by
  \begin{equation}
    (\d\theta)(B) = \rho^*(B) \theta := - \theta \circ \rho(B)
  \end{equation}
  and if $\lambda \in C^1(\fb,W^*)$ then $\d\lambda \in C^2(\fb;W^*)$
  is given by
  \begin{equation}
    \begin{split}
      \d\lambda(B,C) &= \rho^*(B) \lambda_C - \rho^*(C) \lambda_B - \lambda_{[B,C]}\\
      &= - \lambda_C \circ \rho(B) + \lambda_B \circ \rho(C) - \lambda_{[B,C]}.
    \end{split}
  \end{equation}
  Applying this to $w \in W$ and setting it equal to zero, we find
  precisely equation~\eqref{eq:cocycle}.  Therefore that equation
  becomes $\d\lambda=0$ so that it defines a class $[\lambda] \in
  H^1(\fb;W^*)$.
\end{proof}

It remains to consider the bracket of two elements in $V \otimes W$,
which gives a linear map
\begin{equation}
  \ext{2}(V \otimes W) = \left( \ext{2} V \otimes \sS^2 W \right)
  \oplus \left( \sS_0^2 V \otimes \ext{2} W \right) \oplus \left(
    \RR_{\Tr} \otimes \ext{2} W \right) \to \so(V) \oplus \left( V
    \otimes W\right) \oplus \fb.
\end{equation}
Equivariance under $\so(V)$ now says that the only possible nonzero
components of this linear map are
\begin{equation}
  \label{eq:3d-brackets}
      \ext{2} V \otimes \sS^2 W \to \so(V) \oplus \left( V \otimes W \right) \quad
      \text{and}\quad
      \RR_{\Tr} \otimes \ext{2} W \to \fb.
\end{equation}
The first linear map itself breaks up into two components.
The first component is
\begin{equation}
  \ext{2} V \otimes \sS^2 W  \to\so(V),
\end{equation}
which given the isomorphism $\ext{2} V \cong \so(V)$ in the previous
section, is characterised by a symmetric bilinear form $\alpha \colon W
\times W \to \RR$.  The second component is new:
\begin{equation}
  \ext{2} V \otimes \sS^2 W  \to V \otimes W.
\end{equation}
The isomorphism $\ext{2} V \to V$ is essentially the vector cross
product $v_1 \wedge v_2 \mapsto v_1 \times v_2$ and hence this
component is characterised by a $W$-valued symmetric bilinear form
$\sigma\colon W \times W \to W$.  The second linear map in
equation~\eqref{eq:3d-brackets} is as in the generic case and is
characterised by a $\fb$-valued $2$-form
$\varphi\colon \ext{2} W \to \fb$.

In summary, the bracket is given by
\begin{equation}
  [v_1 \otimes w_1, v_2 \otimes w_2] =\alpha(w_1,w_2) v_1 \curlywedge
  v_2 + \left<v_1, v_2\right> \varphi(w_1 \wedge w_2) + (v_1 \times
  v_2) \otimes \sigma(w_1,w_2),
\end{equation}
for all $v_i \in V$ and $w_i \in W$.

The Jacobi identity with $B \in \fb$ says that the brackets are
$\fb$-equivariant.  Let $B \in \fb$ and consider
\begin{multline}
  [B,[v_1 \otimes w_1,v_2\otimes w_2]] = \left<v_1,v_2\right>
  [B,\varphi(w_1\wedge w_2)]\\ + (v_1 \times v_2) \otimes
  \rho(B)\sigma(w_1,w_2) + \lambda_B(\sigma(w_1,w_2)) \varepsilon(v_1
  \times v_2).
\end{multline}
The Jacobi identity says that
\begin{equation}
  [B,[v_1 \otimes w_1,v_2\otimes w_2]] = [[B,v_1 \otimes w_1], v_2
  \otimes w_2] + [v_1 \otimes w_1, [B, v_2 \otimes w_2]],
\end{equation}
and this equation has three components according to $\g = \so(V) \oplus \fb
\oplus (V \otimes W)$.  The $\so(V)$-component yields
\begin{equation}
  \left( \alpha(\rho(B)w_1,w_2) + \alpha(w_1,\rho(B)w_2)   \right) v_1
  \curlywedge v_2 - \lambda_B(\sigma(w_1,w_2)) \varepsilon(v_1 \times v_2) = 0.
\end{equation}
Using that $\varepsilon(v_1\times v_2) = - v_1 \curlywedge v_2$, we
arrive at
\begin{equation}
  \alpha(\rho(B)w_1,w_2) + \alpha(w_1,\rho(B)w_2) + \lambda_B(\sigma(w_1,w_2)) = 0,
\end{equation}
which measures the lack of $\fb$-equivariance of $\alpha$:
\begin{equation}
  \label{eq:B-alpha-lambda}
  (B\cdot \alpha)(w_1, w_2) = \lambda_B(\sigma(w_1,w_2)).
\end{equation}

The $\fb$-component simply says that $\varphi \colon \ext{2}W \to
\fb$ is $\fb$-equivariant:
\begin{equation}
  \varphi(\rho(B) w_1 \wedge w_2) + \varphi(w_1 \wedge \rho(B)w_2) =
  [B,\varphi(w_1 \wedge w_2)].
\end{equation}
Finally, the $V\otimes W$-component says that
\begin{equation}
  \label{eq:B-sigma-lambdad-expanded}
  \sigma(\rho(B)w_1, w_2) + \lambda_B(w_1)w_2 + \sigma(w_1,
  \rho(B)w_2) + \lambda_B(w_2)w_1 - \rho(B)\sigma(w_1,w_2)= 0,
\end{equation}
which says that $\lambda$ also measures the failure of
$\sigma\colon \sS^2 W \to W$ being $\fb$-equivariant:
\begin{equation}
  \label{eq:B-sigma-lambda}
  (B \cdot \sigma)(w_1,w_2) = \lambda_B(w_1)w_2 + \lambda_B(w_2)w_1.
\end{equation}
This has an important consequence:

\begin{proposition}
  The cohomology class $[\lambda] \in H^1(\fb;W^*)$ is zero; that is, $\lambda = \d \mu$ for some $\mu \in W^*$.
\end{proposition}

\begin{proof}
  The bilinear form $\sigma \colon W \times W \to W$ is equivalent, by
  currying, to a linear map $\hat \sigma \colon W \to \End(W)$ sending
  $w \mapsto \hat\sigma_w$, where $\hat\sigma_w (w') =
  \sigma(w,w')$.  This allows us to abstract $w_2$ from
  equation~\eqref{eq:B-sigma-lambdad-expanded} in order to obtain the
  following equation in $\End(W)$:
  \begin{equation}
    \rho(B) \circ \hat\sigma_{w_1} - \hat\sigma_{w_1} \circ
    \rho(B) - \hat\sigma_{\rho(B) w_1} = \lambda_B(w_1) \id_W +
    w_1 \otimes \lambda_B.
  \end{equation}
  We now take the trace $\Tr \colon \End(W) \to \RR$ and use that 
  $\Tr(XY) = \Tr (YX)$ to arrive at
  \begin{equation}
    -\Tr \hat\sigma_{\rho(B)w_1} = 3 \lambda_B(w_1).
  \end{equation}
  Let us define $\mu \in W^*$ by $\mu := \frac13 \Tr \circ
  \hat\sigma$.  In terms of $\mu$, and abstracting $w_1$
  from the previous equation, we arrive at
  \begin{equation}
    \lambda_B = - \mu \circ \rho(B) = \rho^*(B) \mu,
  \end{equation}
  where $\rho^*(B)$ is the dual representation $\rho^*\colon \fb \to
  \gl(W^*)$.  In other words, $\lambda = \d\mu$.
\end{proof}

This means that we may actually set $\lambda = 0$ without changing the
isomorphism class of $\g$.

\begin{proposition}\label{prop:automorphisms-from-W*}
  The Lie algebra $\g$ is isomorphic to one with $\lambda = 0$.
\end{proposition}

\begin{proof}
  Let $\lambda = \d\mu$. We define $j_\mu \in \GL(\g)$ to be the invertible
  endomorphism characterised by
  \begin{equation}
    j_\mu(v \otimes w + X + B) =v \otimes w - \mu(w)\varepsilon(v) + X + B,
  \end{equation}
  for all $X\in\so(V),\,v\in V,\,w\in W$ and $B\in\fb$.
  The inverse of $j_\mu$ satisfies
  \begin{equation}
    j^{-1}_\mu (v \otimes w + X + B) = v \otimes w + \mu(w)\varepsilon(v) + X + B,
  \end{equation}
  for all $v \in V$, $w \in W$, $X \in \so(V)$ and $B \in \fb$.  Then
  the transformed Lie bracket $[x,y]' := j_\mu^{-1}[j_\mu(x),
  j_\mu(y)]$ is given by the (unchanged) brackets involving elements in
  $\so(V)$ and, in addition,
  \begin{equation}
    \label{eq:3d32brackets}
    \begin{split}
      [B,v \otimes w]'={}&v\otimes \rho(B)w, \\
      [v_1\otimes w_1, v_2\otimes w_2]'={}&\alpha'(w_1,w_2) v_1 \curlywedge v_2 + \left<v_1, v_2\right> \varphi(w_1 \wedge w_2) + (v_1 \times v_2) \otimes \sigma'(w_1, w_2),
    \end{split}
  \end{equation}
  for all $v,v_i\in V$, $w,w_i\in W$ and $B\in \fb$, where
  \begin{equation}
    \begin{split}
      \alpha' &= \alpha - \mu \circ \sigma +  \mu \otimes \mu\\
      \sigma' &= \sigma - \mu \otimes \id - \id \otimes \mu.
    \end{split}
  \end{equation}
  This proves the assertion.\qedhere
\end{proof}

It remains to consider the Jacobi identity with three vectors in $V
\otimes W$, but now setting $\lambda = 0$.  The nested bracket has
again three components:
\begin{equation}
  \begin{split}
    [[v_1\otimes w_1, v_2\otimes w_2],v_3 \otimes w_3]_{\so(V)} ={}& \alpha(\sigma(w_1,w_2),w_3) (v_1 \times v_2) \curlywedge v_3 \\
    [[v_1\otimes w_1, v_2\otimes w_2],v_3 \otimes w_3]_{\fb} ={}& \left<v_1 \times v_2, v_3\right> \varphi(\sigma(w_1,w_2) \wedge w_3),\\
    [[v_1\otimes w_1, v_2\otimes w_2],v_3 \otimes w_3]_{V \otimes W} ={}& \alpha(w_1,w_2) \left( \left<v_2,v_3\right> v_1 - \left<v_1,v_3\right> v_2\right)\otimes w_3\\
    & +  \left<v_1, v_2\right> v_3 \otimes \rho(\varphi(w_1 \wedge w_2)) w_3\\
    & + ((v_1 \times v_2) \times v_3) \otimes \sigma(\sigma(w_1,w_2),w_3).
  \end{split}
\end{equation}
Each component must vanish separately in the Jacobi identity.
The $\so(V)$ terms give the equation
\begin{equation}
  \label{eq:3d-soV-jacobi}
  \alpha(\sigma(w_1,w_2),w_3) (v_1 \times v_2) \curlywedge v_3 +  \text{cyclic} = 0,
\end{equation}
whereas the $\fb$ terms give
\begin{equation}
  \label{eq:3d-b-jacobi}
  \left<v_1 \times v_2, v_3\right> \varphi(\sigma(w_1,w_2) \wedge w_3) + \text{cyclic} = 0
\end{equation}
and the $V \otimes W$ terms give
\begin{multline}
  \label{eq:3d-VW-jacobi}
  \alpha(w_1,w_2) \left( \left<v_2,v_3\right> v_1 - \left<v_1,v_3\right> v_2\right)\otimes w_3  +  \left<v_1, v_2\right> v_3 \otimes  \rho(\varphi(w_1 \wedge w_2)) w_3\\
  + ((v_1  \times v_2) \times v_3) \otimes \sigma(\sigma(w_1,w_2),w_3) + \text{cyclic}= 0.
\end{multline}
Using the fact that $u \curlywedge v = - \varepsilon (u \times v)$ and
the well-known cross product identity
\begin{equation}
  \label{eq:cross-prod-ident}
  (v_1 \times v_2) \times v_3 = \left<v_1, v_3\right> v_2 - \left<v_2, v_3\right> v_1,
\end{equation}
we may rearrange the terms in equation~\eqref{eq:3d-soV-jacobi} to give
\begin{equation}
  \label{eq:3d-soV-jacobi-2}
  \left<v_1,v_2\right> \left( \alpha(\sigma(w_3,w_1),w_2) - \alpha(\sigma(w_3,w_2),w_1) \right)
  \varepsilon(v_3) + \text{cyclic}= 0.
\end{equation}
This has to hold for every choice of $v_1,\,v_2,\,v_3$.
Consequently, if we take $v_1=v_2$ to be nonzero and $v_3$ to be
a nonzero vector that is orthogonal to $v_1$ and $v_2$, this component
of the Jacobi identity becomes
\begin{equation}
  \label{eq:3d-jacobi-1}
  \alpha(\sigma(w_3,w_2),w_1) = \alpha(\sigma(w_3,w_1),w_2)
\end{equation}
for all $w_i \in W$.  It is clear that if this equation holds, then so
does equation~\eqref{eq:3d-soV-jacobi-2}.

The expression $\left<v_1 \times v_2, v_3 \right>$ in
equation~\eqref{eq:3d-b-jacobi} is cyclically invariant, hence that
equation becomes simply
\begin{equation}
  \label{eq:3d-jacobi-2}
  \varphi(\sigma(w_1,w_2) \wedge w_3) +  \varphi(\sigma(w_2,w_3) \wedge w_1) +   \varphi(\sigma(w_3,w_1) \wedge w_2)  = 0.
\end{equation}

In equation~\eqref{eq:3d-VW-jacobi}, the expression $(v_1 \times v_2)
\times v_3$ can be expanded using the cross product identity~\eqref{eq:cross-prod-ident}.
Inserting this into equation~\eqref{eq:3d-VW-jacobi} and rearranging,
we may rewrite that equation as
\begin{multline}
  \left<v_1,v_2\right> v_3 \otimes \left( \rho(\varphi(w_1\wedge w_2)) w_3 -
  \alpha(w_2,w_3) w_1 + \alpha(w_1,w_3) w_2\right. \\ \left.  +
  \sigma(\sigma(w_2,w_3),w_1) - \sigma(\sigma(w_1,w_3),w_2)\right) +
\text{cyclic}= 0.
\end{multline}
Choosing $v_i \in V$ so that $\left<v_1,v_2\right>$ and $v_3$ are nonzero but $\left<v_1,v_3\right> =
\left<v_2,v_3\right> = 0$, we obtain the equation
\begin{equation}
  \label{eq:3d-jacobi-3}
  \alpha(w_2,w_3)w_1-\alpha(w_1,w_3)w_2=\rho(\varphi(w_1\wedge w_2))w_3+\sigma(\sigma(w_2,w_3),w_1)-\sigma(\sigma(w_1,w_3),w_2) 
\end{equation}
for all $w_i \in W$.  Conversely, if this holds, then so does
equation~\eqref{eq:3d-VW-jacobi}.  Similarly to the case $d>3$, the
bilinear form $\alpha$ can be recovered from $\rho$, $\varphi$ and
$\sigma$. Indeed, suppose that we are given maps $\varphi\colon
\ext{2} W \to \fb$ and $\sigma\colon \sS^2 W \to W$ such that
$\varphi$ is $\fb$-equivariant, and define a map $\phi\colon \ext{2} W
\to \gl(W)$ via
\begin{equation}
  \phi(w_1\wedge w_2)w_3=\rho(\varphi(w_1\wedge w_2))w_3+\sigma(\sigma(w_2,w_3),w_1)-\sigma(\sigma(w_1,w_3),w_2), \quad w_i\in W.
\end{equation}
This map actually takes values in $\fsl(W)$.  Indeed, the first
summand yields a traceless endomorphism because of $\fb$-equivariance
of $\varphi$.  As for the other two summands, we may recognise them as
the commutator $[\hat{\sigma}_{w_1},\hat{\sigma}_{w_2}]$, where
$\hat{\sigma}_w$ is the endomorphism taking $w'\mapsto \sigma(w,w')$
for all $w\in W$.  In particular, the trace
$\Tr([\hat{\sigma}_{w_1},\hat{\sigma}_{w_2}])=0$, so we conclude that
$\im \phi \subseteq \fsl(W)$. Owing to the
isomorphism~\eqref{eq:iso-s2W*-homwedge2slW}, we can guarantee the
existence of a unique $\alpha\in \sS^2 W^*$ that is related to $\rho$,
$\varphi$ and $\sigma$ via~\eqref{eq:3d-jacobi-3}.  In particular, if
$\sigma$ is also $\fb$-equivariant, then $\phi$ is $\fb$-invariant
(due to equivariance of $\sigma$ and $\varphi$), and $\alpha\in (\sS^2
W^*)^\fb$.

\begin{remark}
    Just as in the generic case,
    we can compute the coordinates of $\alpha$ from the rest of ingredients.
    Let $\be_1,\be_2$ be a basis for $W$ and put $w_1 = \be_1$ and $w_2 = \be_2$
    in equation~\eqref{eq:3d-jacobi-3}. We define the matrix
    $[\rho^j{}_i]$ by
    \begin{equation}
        \rho(\varphi(\be_1\wedge \be_2)) \be_i = \rho^j{}_i \be_j. 
    \end{equation}
    In addition, we consider $\hat\sigma_i := \sigma(\be_i,-)$ and define coefficients $[\hat\sigma_1,
    \hat\sigma_2]^j{}_i$ by means of the equation
    \begin{equation}
        \sigma(\be_1,\sigma(\be_2, \be_i)) - \sigma(\be_2,
        \sigma(\be_1,\be_i)) =  [\hat\sigma_1, \hat\sigma_2] \be_i = \be_j [\hat\sigma_1, \hat\sigma_2]^j{}_i.
    \end{equation}
    Then, arguing in the same fashion as in the case $d>3$, we obtain
    \begin{equation}
        \begin{split}
            \alpha_{11} &= - \rho^2{}_1 - [\hat\sigma_1, \hat\sigma_2]^2{}_1,\\
            \alpha_{12} &= \rho^1{}_1 + [\hat\sigma_1, \hat\sigma_2]^1{}_1= - \rho^2{}_2 - [\hat\sigma_1, \hat\sigma_2]^2{}_2,\\
            \alpha_{22} &= \rho^1{}_2 +[\hat\sigma_1, \hat\sigma_2]^1{}_2.
        \end{split}
    \end{equation}
\end{remark}

In conclusion, we have seen that
(up to isomorphism that acts trivially on $\so(V)$)
the Lie bracket of a $(3,2)$-KLA $\g$ with three-dimensional
spatial isotropy is determined by a Bianchi Lie algebra $\fb$,
a two-dimensional real representation $W$ of $\fb$,
a $\fb$-equivariant map $\varphi\in\Hom_\fb(\ext{2}W,\fb)$
and a $\fb$-equivariant symmetric bilinear map
$\sigma\colon W\times W\to W$.
Moreover, taking into consideration the unique symmetric
bilinear map $\alpha\in\sS^2W^*$ characterised by~\eqref{eq:3d-jacobi-3},
then $\varphi$, $\sigma$ and $\alpha$ must satisfy the
relations~\eqref{eq:3d-jacobi-1} and~\eqref{eq:3d-jacobi-2}.
Conversely, given a quadruple $(\fb,W,\varphi,\sigma)$ that
satisfies all of these properties, the bilinear operation on $\g$
given by~\eqref{eq:3d-lie-bracket} yields a $(3,2)$-KLA structure
on $\g$.
All in all, we have proved Theorem~\ref{thm:3d-32StructureTheorem-simplified}.

\section{Relative isomorphisms between $(3,2)$-kinematical Lie algebras}
\label{sec:relative-isomorphisms}

In this section, we describe a suitable notion of isomorphism between
two $(3,2)$-kinematical Lie algebras.  Recall that our final objective
is the classification of coisotropy-one homogeneous spacetimes up to
isomorphisms and coverings.  Consequently, for our purposes it is
enough to classify these algebras up to isomorphisms that preserve
rotations.

A first observation is the following:
let $\g$ and $\g'$ be two $(3,2)$-KLAs, which may be written as
\begin{equation}\label{eq:two-KLA-decompositions}
    \g=\so(d)\oplus\RR^d\oplus\RR^d\oplus\RR\oplus\RR\oplus\RR,\quad
    \g'=\so(d')\oplus\RR^{d'}\oplus\RR^{d'}\oplus\RR\oplus\RR\oplus\RR,
\end{equation}
where we have decomposed $\g$ (respectively, $\g'$) as a direct sum of $\so(d)$-modules (respectively, $\so(d')$-modules).
Note that, in writing the decompositions in~\eqref{eq:two-KLA-decompositions}, we are implicitly
choosing two embeddings $i\colon\so(d)\to\g$ and $j\colon\so(d')\to\g'$ given by
$i(X)=X$ for all $X\in\so(d)$ and $j(X)=X$ for all $X\in\so(d')$.
Since
\begin{equation}
	\dim \g = \frac{d(d-1)}{2}+2d+3
\end{equation}
and similarly with $\g'$, it is clear that for $\g$ and $\g'$ to be isomorphic we need $d=d'$.
Moreover, we will say that a map $f\colon\g\to\g'$ is a \textit{relative isomorphism} if it is a Lie
algebra isomorphism and $j=f\circ i$.

Clearly, a relative isomorphism $f\colon\g\to\g'$ sends the rotational subalgebra of $\g$
to that of $\g'$.
However, the notion of relative isomorphism depends on the choice of embeddings
$i\colon\so(d)\to\g$ and $j\colon\so(d)\to\g'$.
In addition, a Lie algebra isomorphism $f\colon\g\to\g'$ that preserves $\so(d)$
may restrict to a nontrivial automorphism of $\so(d)$.
Nevertheless, the next proposition shows that the existence of an isomorphism
preserving the rotational subalgebras is equivalent to the existence of a relative
isomorphism between two $(3,2)$-KLAs.

\begin{proposition}\label{prop:normaliser}
    Consider two $(3,2)$-kinematical algebras $\g$, $\g'$ with $d$-dimensional spatial isotropy, where $d\geq3$,
    and choose decompositions for $\g$ and $\g'$ as in~\eqref{eq:two-KLA-decompositions}.
    If there exists a Lie algebra isomorphism $f\colon \g \to \g'$
    that sends the rotational subalgebra of $\g$
    to the rotational subalgebra of $\g'$,
    then there exists a relative isomorphism from $\g$ to $\g'$.
\end{proposition}

\begin{proof}
    For the sake of clarity, we write
    \begin{equation}
        \g=\so(V)\oplus V \oplus V \oplus \fb, \quad
        \g'=\so(V)\oplus V \oplus V \oplus \fb',
    \end{equation}
    where $V=\RR^d$
    and $\fb$, $\fb'$ are three-dimensional vector spaces.
    Under the adjoint representation, the action of $\so(V)$ on each copy of $V$ is given by the standard representation,
    whereas $\fb$ becomes the three-dimensional trivial representation of $\so(V)$.
    The situation is analogous for $\g'$.

    Consider a Lie algebra isomorphism $f\colon \g \to \g'$ that satisfies $f(\so(V))=\so(V)$.
    We define $\eta \colon \so(V)\to \so(V)$ as the restriction of $f$ to $\so(V)$,
    which is a Lie algebra automorphism.
    We claim that $\eta=\Ad(g)$ for some $g\in \mathsf{O}(V)$.
    To see this, we consider two possible situations.

    Firstly, suppose that $\dim V = d \neq 8$.
    In this case, it is well-known that the automorphism group $\Aut(\so(V))$ coincides with $\Ad(\mathsf{O}(V))$,
    so our claim follows immediately.

    Now, assume that $\dim V = d = 8$.
    Here we see that the automorphism group is the semidirect product $\mathfrak{S}_3 \ltimes \Ad(\mathsf{SO}(V))$, where $\mathfrak{S}_3$ is the symmetric group on three elements.
    More precisely, $\mathfrak{S}_3$ is the outer automorphism group of $\so(V)$.
    In addition, the map $\eta$ induces a new $\so(V)$-module structure on $V$
    that turns it into an irreducible eight-dimensional representation $V^\eta$.
    Consequently, $V^\eta$ is equivalent to either the standard representation $V$ or to one of the two half-spin representations $\Delta^+,\,\Delta^-$.

    Let us label the vector part of $\g$ as $V_1 \oplus V_2$,
    where $V_1 = V_2 = V$ are standard representations of $\so(V)$.
    We define maps $\kappa \colon V \to \so(V)$, $\lambda_i \colon V \to V_i$ (with $i=1,2$)
    and $\zeta\colon V \to \fb'$ by the equation
    \begin{equation}
        f(v)=\kappa(v)+\lambda_1(v)+\lambda_2(v)+\zeta(v), \quad v\in V_1.
    \end{equation}
    Because $f$ is a Lie algebra isomorphism we must have $[f(X),f(v)]=f(Xv)$
    for every $X\in \so(V)$ and $v\in V_1$, and this equation amounts to
    \begin{align*}
        \kappa(Xv)={}&[\eta(X),\kappa(v)], & \lambda_1(Xv)={}&\eta(X)\lambda_1(v), & \lambda_2(Xv)={}&\eta(X)\lambda_2(v), & \zeta(Xv)={}&0.
    \end{align*}
    The last equation above implies that $\zeta$ is the zero map.
    The first one implies that $\kappa$ is an $\so(V)$-equivariant map between the standard representation $V$
    and the representation $\so(V)^\eta$ whose underlying vector space is $\so(V)$
    and where the action of $\so(V)$ on itself is given by the map $\ad \circ \eta$.
    Because both $V$ and $\so(V)^\eta$ are irreducible representations of different dimensions, we immediately obtain $\kappa = 0$.
    Furthermore, since $f$ is an isomorphism, we have that either $\lambda_1$ or $\lambda_2$ is nonzero,
    and both turn out to be equivariant maps from $V$ to $V^\eta$, due to the second and third equations above.
    As $V$ and $V^\eta$ are irreducible, we conclude that $V\cong V^\eta$.
    This is possible if and only if $\eta$ is an inner automorphism
    (that is, of the form $\Ad(g)$ for some $g\in\mathsf{SO}(V)$)
    or if the class of $\eta$ in $\mathrm{Out}(\so(V))=\Aut(\so(V))/\Ad(\mathsf{SO}(V))=\mathfrak{S}_3$ is the transposition that interchanges the half-spin representations (this is equivalent to $\eta=\Ad(g)$ for some $g\in \mathsf{O}(V)$, $\det g = -1$). 
    All in all, we have shown that $\eta$ is the adjoint action of an element of $\mathsf{O}(V)$, as desired.

    We have deduced that the restriction of $f$ to $\so(V)$
    is of the form $\Ad(g)$ for some $g\in \mathsf{O}(V)$.
    We define a linear map $s_g\colon \g \to \g$ by declaring that
    \begin{equation}\label{eq:embedding-o(d)}
        s_g(X)={}gXg^{-1}, \quad
        s_g(v)={}gv, \quad
        s_g(B)={}B,
    \end{equation}
    for each $X\in \so(V)$, $v\in V_i$ and $B\in \fb$.
    One checks using Theorems~\ref{th:generic-d-32-structure-theorem}
    and~\ref{thm:3d-32StructureTheorem-simplified} that $s_g$ is a Lie algebra automorphism.\footnote{If $d=3$ and $\sigma \neq 0$, then we need to further impose that 
    $g\in\mathsf{SO}(V)$ so that $s_g$ is a homomorphism.
    This is not restrictive since $\Aut(\so(3))=\Ad(\mathsf{SO}(3))$.}
    Thus, the map $f\circ s_g^{-1}$ is a relative isomorphism from $\g$ to $\g'$,
    as desired. \qedhere
\end{proof}

Motivated by this result, we say that two $(3,2)$-KLAs $\g$ and $\g'$ (with $d\geq 3$) are \textit{relatively isomorphic}
if there exists a Lie algebra isomorphism $f\colon\g\to\g'$ that sends the rotational
subalgebra of $\g$ to the rotational subalgebra of $\g'$.
This notion is equivalent to the existence of a relative isomorphism between $\g$
and $\g'$.
In addition, if $\fr$ denotes the rotational subalgebra of $\g$,
we define the \textit{relative automorphism group} $\Aut(\g,\fr)$
of $\g$ as the group
of all Lie algebra automorphisms of $\g$ that act as the identity on $\fr$.

Observe that Proposition~\ref{prop:normaliser} can also be used
to describe the normaliser of $\fr\cong\so(V)$ in the full automorphism group $\Aut(\g)$.
Indeed, suppose that $\g=\so(V)\oplus (V\otimes W)\oplus \fb$ is a $(3,2)$-KLA where $\dim V = d\geq 3$.
Denote by $\mathsf{N}$ the normaliser of $\so(V)$ in $\Aut(\g)$.
Note that $\Aut(\g,\fr)$ is a normal subgroup of $\mathsf{N}$.

If $d$ is even, then we can define a Lie group embedding $s\colon \mathsf{O}(V)\to\Aut(\g)$
by letting $s_g=\Ad(g)\oplus(g\otimes\id_W)\oplus\id_{\fb}$ for every $g\in\mathsf{O}(V)$.
Proposition~\ref{prop:normaliser} then shows that $\mathsf{N}=\mathsf{O}(V)\Aut(\g,\fr)$,
where the intersection $\mathsf{O}(V)\cap\Aut(\g,\fr)=\{s_{\id_V}, s_{-\id_V}\}$ is finite.
This implies that the normaliser of $\so(V)$ in $\Aut(\g)$ is the almost semidirect product
of $\mathsf{O}(V)$ and $\Aut(\g,\fr)$.

Now, assume $d$ is odd.
The restriction of $s$ to $\mathsf{SO}(V)$ is also a Lie group embedding
(including when $d=3$, where $s\colon\mathsf{O}(V)\to\Aut(\g)$ may not be a homomorphism).
In this case, the proof of Proposition~\ref{prop:normaliser} yields that
$\mathsf{N}=\mathsf{SO}(V)\Aut(\g,\fr)$, where the intersection
$\mathsf{SO}(V)\cap\Aut(\g,\fr)=\{s_{\id_V}\}$ is trivial.
Thus, we obtain $\mathsf{N}=\mathsf{SO}(V)\ltimes\Aut(\g,\fr)$.

All in all, we have proved the following result:

\begin{corollary}
    Let $\g$ be a $(3,2)$-KLA with $d$-dimensional spatial isotropy, where $d\geq 3$,
    and let $\mathsf{N}=\mathsf{N}_{\Aut(\g)}(\so(V))$ be the normaliser
    of the rotational subalgebra $\fr\cong\so(V)$.
    Then the following assertions are true:
    \begin{enumerate}[label=\emph{(\roman*)}]
        \item
            If $d$ is odd, then $\mathsf{N}$ splits as the semidirect product
            $\mathsf{SO}(V)\ltimes \Aut(\g,\fr)$.
        \item
            If $d$ is even, then $\mathsf{N}$ splits as the
            almost semidirect product $\mathsf{O}(V)\Aut(\g,\fr)$,
            where the intersection $\mathsf{O}(V)\cap \Aut(\g,\fr)=\{
            s_{\id_V}, s_{-\id_V} \}$, where the Lie group embedding $s : \mathsf{O}(V) \to
            \Aut(\g,\fr)$ is given by $s_g = \Ad(g)\oplus(g\otimes\id_W)\oplus\id_{\fb}$.
    \end{enumerate}
\end{corollary}

For the remainder of this section, we concern ourselves with determining when two
$(3,2)$-KLAs $\g$ and $\g'$ are relatively isomorphic in terms of their corresponding
algebraic objects described in Theorems~\ref{th:generic-d-32-structure-theorem}
and~\ref{thm:3d-32StructureTheorem-simplified}.
Moreover, we will determine the set of all relative isomorphisms between
two (relatively isomorphic) $(3,2)$-KLAs.

\subsection{Generic case: $d>3$}
\label{sec:isomorphisms-d-generic}
Consider two $(3,2)$-KLA $\g$ and $\g'$, which we write as
\begin{equation}
	\g=\so(V)\oplus (V\otimes W) \oplus \fb, \quad
  \g'=\so(V)\oplus (V\otimes W') \oplus \fb',
\end{equation}
where $\dim V = d > 3$.
Then, according to Theorem~\ref{th:generic-d-32-structure-theorem} the
Lie brackets of $\g$ and $\g'$ are determined by triples
$(\fb,\rho,\varphi)$ and $(\fb',\rho',\varphi')$, where $\fb$ is a
Bianchi algebra, $\rho\colon \fb \to \gl(W)$ is a two-dimensional
representation and $\varphi\colon \ext{2} W \to \fb$ is a
$\fb$-equivariant map (and likewise for $\fb'$, $\rho'$ and $\varphi'$).

Suppose a relative isomorphism $f\colon \g \to \g'$ is given.  As
$f\vert_{\so(V)}$ is the identity map, $f$ is also an isomorphism of
$\so(V)$-modules.  In particular, $f$ sends the isotypical components
of $\g$ to those of $\g'$.  For $d>3$, the isotypical components of
$\g$ are $\so(V)$ (when $d=4$, $\so(V)$ splits into two), $V\otimes W$
and $\fb$, and identically for $\g'$.  As a consequence, we have
$f(V\otimes W)=V\otimes W'$ and $f(\fb)=\fb'$, and the equivariance of
$f\vert_{V\otimes W}$ yields that there exists a map $h\colon W\to W'$
such that $f(v\otimes w)=v\otimes h(w)$ for all $v\otimes w\in
V\otimes W$.  We have thus shown that $f$ decomposes as
\begin{equation}
	f=\id_{\so(V)}\oplus(\id_{V}\otimes h)\oplus \psi,
\end{equation}
for a certain linear isomorphism $h\colon W\to W'$ and a Lie algebra
isomorphism $\psi\colon \fb\to \fb'$.

Let $B\in \fb$, $v\otimes w\in V\otimes W$.
The equations
\begin{equation}
  \begin{split}
    f([B,v\otimes w]) &= f(v\otimes \rho(B)w)=v\otimes h(\rho(B)w), \\
    [f(B),f(v\otimes w)] &= [\psi(B),v\otimes h(w)]=v\otimes \rho'(\psi(B))h(w)
  \end{split}
\end{equation}
imply that the maps $h \circ \rho(B)$ and $\rho'(\psi(B))\circ h$ are equal.
For a linear isomorphism $g\colon W\to W'$, let us define a map $\Ad(g)\colon \gl(W)\to\gl(W')$ via $\Ad(g)T=g\circ T\circ g^{-1}$.
Then we can rewrite the previous condition as $\rho'\circ \psi = \Ad(h)\circ \rho$.
Furthermore, choose $v_{1},v_{2}\in V$ and $w_{1},w_{2}\in W$.
We have
\begin{equation}
  \begin{split}
    f([v_1 \otimes w_1, v_2 \otimes w_2]) &= f(\alpha(w_1 , w_2)v_1 \curlywedge v_2 + \langle v_1, v_2 \rangle\varphi(w_1 \wedge w_2)) \\
    &= \alpha(w_1, w_2)v_1 \curlywedge v_2 + \langle v_1, v_2 \rangle \psi(\varphi(w_1 \wedge w_2)),
  \end{split}
\end{equation}
whereas
\begin{equation}
  \begin{split}
    [f(v_1 \otimes w_1),f(v_2 \otimes w_2)] &= [v_1 \otimes h(w_1), v_2 \otimes h(w_2)] \\
    &= \alpha' (h(w_1),h(w_2))v_1 \curlywedge v_2 + \langle v_1 , v_2 \rangle \varphi'(h(w_1)\wedge h(w_2)) \\
    &= (h^* \alpha')(w_1 , w_2)v_1 \curlywedge v_2 + \langle v_1 , v_2 \rangle (h^*\varphi') (w_1 \wedge w_2).
  \end{split}
\end{equation}
Therefore, we must have $\alpha = h^* \alpha'$ and $\psi \circ \varphi
= h^*\varphi'$.  It is not hard to see that the latter equation,
together with $\rho' \circ \psi = \Ad(h) \circ \rho$, implies the former.

Conversely, suppose $h\colon W\to W'$ is a linear isomorphism and
$\psi\colon \fb \to \fb'$ is a Lie algebra isomorphism such that
$\rho' \circ \psi = \Ad(h) \circ \rho$ and $\psi \circ \varphi =
h^*\varphi'$.  Then from the equalities above, it is clear that the
map $f\colon \g \to \g'$ defined by $f=\id_{\so(V)}\oplus (\id_V
\otimes h)\oplus \psi$ is a relative isomorphism.

We summarise this discussion in the following result:

\begin{proposition}\label{prop:generic-d-relative-isomorphisms}
  Let $\g=\so(V)\oplus (V\otimes W)\oplus \fb$ and $\g' = \so(V)\oplus
  (V\otimes W')\oplus \fb'$ be two $(3,2)$-KLAs for which $\dim V =
  d>3$, and let $(\fb,\rho, \varphi)$ and $(\fb',\rho', \varphi')$
  be the algebraic data defining the Lie brackets of $\g$ and $\g'$.
  Then $\g$ and $\g'$ are relatively isomorphic if and only
  if there exist
  \begin{enumerate}[label=\emph{(\roman*)}]
  \item a vector space isomorphism $h\colon W \to W'$,
  \item and a Lie algebra isomorphism $\psi\colon \fb \to \fb'$,
  \end{enumerate}
  satisfying the following equations:
  \begin{equation}
    \rho'\circ \psi = \Ad(h)\circ \rho \qquad\text{and}\qquad \psi \circ \varphi = h^*\varphi'.
  \end{equation}
  The map
  \begin{equation}
    f=\id_{\so(V)} \oplus (\id_V\otimes h)\oplus \psi
  \end{equation}
  is a relative isomorphism from $\g$ to $\g'$.
  Furthermore, every relative isomorphism arises in this manner.
\end{proposition}

\subsection{Case $d=3$}
\label{sec:isomorphisms-d-3}

We now study the possible isomorphisms between $(3,2)$-KLAs for which $d=3$.
We consider two $(3,2)$-KLAs of the form
\begin{equation}
	\g = \so(V)\oplus (V \otimes W)\oplus \fb, \quad \g' = \so(V)\oplus(V\otimes W')\oplus \fb',
\end{equation}
where $\dim V = d = 3$.  Due to
Theorem~\ref{thm:3d-32StructureTheorem-simplified},
we can assume that (up to a relative isomorphism) the Lie brackets
of $\g$ and $\g'$ are determined by the $5$-tuples
$(\fb,\rho,\phi,\sigma,\alpha)$ and $(\fb',\rho',\phi',\sigma',\alpha')$,
where $\fb$ is a Bianchi algebra, $\rho\colon\fb\to\gl(W)$ is a representation,
$\phi\in\Hom(\ext{2}W,\fb)$, $\alpha\in\sS^2W^*$, 
and $\sigma\in\Hom(\sS^2W,W)$ are such that the relations~\eqref{eq:3d-algebraic-data-relations}
are satisfied (and similarly for $\fb',\rho',\phi',\sigma',\alpha'$).
We aim to compute all relative isomorphisms from $\g$ to $\g'$.

Suppose that $f \colon \g \to \g'$ is a relative isomorphism.  Then
$f$ is also an $\so(V)$-module isomorphism, so it must preserve the
isotypical decompositions of $\g$ and $\g'$.  For $d=3$, the isotypical
components of $\g$ are $\so(V)\oplus(V\otimes W)$ and $\fb$, as
$\so(V)\cong V$ (and the same holds for $\g '$).  This implies that
$f(\fb)=\fb'$, while $f(V\otimes W)$ is a subspace of
$\so(V)\oplus (V\otimes W')$ complementary to $\so(V)$.  The
restriction $\psi\colon \fb \to \fb'$ of $f$ to $\fb$ is a Lie algebra
isomorphism, while the maps
$v\otimes w \mapsto f(v\otimes w)_{\so(V)}$ and
$v\otimes w \mapsto f(v\otimes w)_{V\otimes W}$ are both
$\so(V)$-equivariant.  Since $\Hom_{\so(V)}(V,V)\cong \RR$ and
$\Hom_{\so(V)}(V,\so(V))\cong \RR$ is generated by the map
$\varepsilon\colon V \to \so(V)$, $\varepsilon(v)w=v\times w$, we see
that there exists a covector $\mu\in W^*$ and a linear isomorphism
$h\colon W\to W'$ such that
\begin{equation}
	f(v\otimes w)=\mu(w)\varepsilon(v)+v\otimes h(w), \quad v\in V,w\in W.
\end{equation}

Consider $B\in \fb$, $v\in V$ and $w\in W$.
Because $f$ is an isomorphism, the elements
\begin{equation}
	\begin{aligned}
		f([B,v\otimes w])={}&f(v\otimes \rho(B)w)=\mu(\rho(B)w)\varepsilon(v)+v\otimes h(\rho(B)w)
	\end{aligned}
\end{equation}
and
\begin{equation}
  \begin{aligned}\relax
    [f(B),f(v\otimes w)]={}&[\psi(B),\mu(w)\varepsilon(v)+v\otimes h(w)]=v\otimes \rho'(\psi(B))h(w)
  \end{aligned}
\end{equation}
are equal, and because $B,\,v,\,w$ are arbitrary, this is possible if and only if
\begin{equation}
  h(\rho(B)w)=\rho'(\psi(B))h(w),\quad\mu(\rho(B)w)=0,
\end{equation}
for all choices of $B,\,v,\,w$.
The first condition is equivalent to the equation
\begin{equation}\label{eq:3d-isomorphism-bvw1}
	\rho'\circ \psi = \Ad(h)\circ \rho,
\end{equation}
while the second condition can be read as
\begin{equation}\label{eq:3d-isomorphism-bvw2}
	\partial \mu = 0,
\end{equation}
where we recall that $\partial \colon C^\bullet(\fb;W^*)\to C^{\bullet+1}(\fb;W^*)$ is the Chevalley-Eilenberg differential.
In other words, $\mu$ is invariant under the dual representation of $\rho$.

Now, let $v_{1},v_{2}\in V$ and $w_{1},w_{2}\in W$.
We have
\begin{equation}
  \begin{aligned}
    f([v_1 \otimes w_1, v_2 \otimes w_2])={}&f(\alpha(w_1,w_2)v_1\curlywedge v_2+(v_1\times v_2)\otimes\sigma(w_1,w_2)+\langle v_1 , v_2 \rangle \varphi(w_1 \wedge w_2)) \\
    ={}& (\alpha(w_1,w_2)-\mu(\sigma(w_1, w_2)))v_1\curlywedge v_2+(v_1\times v_2)\otimes h(\sigma(w_1,w_2))\\
    \quad &+ \langle v_1 , v_2 \rangle \psi(\varphi(w_1 \wedge w_2))
  \end{aligned}
\end{equation}
while using the equality $[\varepsilon(v_1),\varepsilon(v_2)]=\varepsilon(v_1\times v_2)=-v_1 \curlywedge v_2$ we see that
\begin{equation}
  \begin{aligned}\relax
    [f(v_1 \otimes w_1),f(v_2 \otimes w_2)]={}&[\mu(w_1)\varepsilon(v_1)+v_1 \otimes h(w_1),\mu(w_2)\varepsilon(v_2)+v_2 \otimes h(w_2)]\\
    ={}&((h^* \alpha')(w_1,w_2)-\mu(w_1)\mu(w_2))v_1 \curlywedge v_2\\
    &+\langle v_1 , v_2 \rangle (h^*\varphi')(w_1 \wedge w_2)\\
    &+(v_1\times v_2)\otimes (\mu(w_1)h(w_2)+\mu(w_2)h(w_1)+(h^* \sigma')(w_1,w_2)).\\
  \end{aligned}
\end{equation}
These two elements are equal for all choices of $v_i,w_i$ if and only if
\begin{equation}\label{eq:3d-isomorphism-vw1}
	h\circ \sigma = h^*\sigma' + \mu \otimes h + h \otimes \mu, \quad \psi\circ \varphi = h^*\varphi'.
\end{equation}
and
\begin{equation}\label{eq:3d-isomorphism-vw2}
  \alpha=h^* \alpha'+\mu \circ \sigma - \mu \otimes \mu.
\end{equation}

Conversely, suppose that we have a linear isomorphism $h\colon W \to W'$, a Lie algebra isomorphism $\psi\colon \fb \to \fb'$ and a one-form $\mu\in W^*$ that satisfy the relations~\eqref{eq:3d-isomorphism-bvw1}, \eqref{eq:3d-isomorphism-bvw2} and~\eqref{eq:3d-isomorphism-vw1}.
Note that we do not require that~\eqref{eq:3d-isomorphism-vw2} is satisfied.
Define a map $f\colon \g \to \g'$ defined by $f\vert_{\so(V)\oplus\fb}=\id_{\so(V)}\oplus \psi$ and
\begin{equation}
	f(v\otimes w)=\mu(w)\varepsilon(v)+v\otimes h(w), \quad v\in V,w\in W,
\end{equation}

We claim that $f$ is a relative isomorphism (which in particular implies~\eqref{eq:3d-isomorphism-vw2}).
To see this, consider the Lie bracket $[\cdot,\cdot]_{*}$ on $\g'$ that turns $f\colon \g \to (\g',[\cdot,\cdot]_{*})$ into a Lie algebra isomorphism.
This Lie bracket satisfies
\begin{equation}
  \begin{split}
    [X,v\otimes w']_{*}={}&f([X,v\otimes h^{-1}(w')-\mu(h^{-1}(w'))\varepsilon(v)])=Xv\otimes w', \\
    \qquad [X,B']_*={}&f([X,\psi^{-1}(B')])=0,
  \end{split}
\end{equation}
for all $X\in\so(V)$, $v\in V$, $w'\in W'$ and $B'\in \fb'$, so $(\g',[\cdot,\cdot]_*)$ is a $(3,2)$-KLA.
From the definition of $f$ and the equations~\eqref{eq:3d-isomorphism-bvw1}, \eqref{eq:3d-isomorphism-bvw2} and \eqref{eq:3d-isomorphism-vw1} we see that the algebraic data associated with $[\cdot ,\cdot]_*$ is $\fb_*=\fb'$, $\rho_*=\rho'$, $\varphi_*=\varphi'$ and $\sigma_*=\sigma'$, and because $\alpha_*$ is obtained from the rest of ingredients, we have necessarily $\alpha_*=\alpha'$.
This gives that $\g$ and $\g'$ (with the original bracket) are isomorphic, and~\eqref{eq:3d-isomorphism-vw2} holds.

We conclude the following:
\begin{proposition}\label{prop:3d-relative-isomorphisms}
  Let $\g=\so(V)\oplus(V\otimes W)\oplus \fb$ and $\g'=\so(V)\oplus
  (V\otimes W')\oplus \fb'$ be two $(3,2)$-KLAs for which $\dim V = d
  = 3$ and let $(\fb,\rho,\varphi,\sigma,\lambda)$ and
  $(\fb',\rho',\varphi',\sigma',\lambda')$ be the algebraic data defining the
  Lie brackets of $\g$ and $\g'$ respectively.
  Assume furthermore that both $\lambda$ and $\lambda'$ vanish.
  Then $\g$ and $\g'$
  are isomorphic relative to $\so(V)$ if and only if there exist
  \begin{enumerate}[label=\emph{(\roman*)}]
  \item a linear isomorphism $h\colon W \to W'$,
  \item a Lie algebra isomorphism $\psi \colon \fb \to \fb'$,
  \item and a covector $\mu\in W^*$,
  \end{enumerate}
  satisfying the following equations:
  \begin{equation}\label{eq:d-3-relative-isomorphism-equations}
    \begin{aligned}
      \rho' \circ \psi ={}& \Ad(h)\circ \rho, & \psi \circ \varphi={}& h^*\varphi', \\
      h\circ \sigma ={}& h^* \sigma' + \mu \otimes h + h \otimes \mu, & \partial \mu ={}& 0.
    \end{aligned}
  \end{equation}
  The map $f\colon \g \to \g'$ defined by
  \begin{equation}
    \begin{split}
      f(X)={}&X, \\
      f(v\otimes w)={}&\mu(w)\varepsilon(v)+v\otimes h(w), \\
      f(B)={}&\psi(B)
    \end{split}
  \end{equation}
  for all $X\in \so(V)$, $v\otimes w\in V\otimes W$ and $B\in \fb$, is a relative isomorphism from $\g$ to $\g'$.
  Furthermore, every relative isomorphism arises in this manner.
\end{proposition}

\section{Lie subalgebras of $\gl(2,\RR)$}\label{sec:gl2-subalgebras}

As we already saw in Section~\ref{sec:struc-32-klas}, part of the
information that describes the general structure of a $(3,2)$-KLA is
given by a two-dimensional representation of a Bianchi Lie algebra.
Because of this, it is necessary to develop a systematic approach to
determining all possible representations of this type by understanding
their possible kernels and images (that is, the ideals of all Bianchi
algebras and the subalgebras of $\gl(2,\RR)$). Thus, the aim of this
section is to provide a classification result of all subalgebras of
$\gl(2,\RR)$ up to conjugacy. Although this classification has been
done several times in the
literature~\cite{MR442045,MR2004893,chapovskyi}, in this paper we
repeat the calculation to fix notation and for completeness.

The key observation is that $\gl(2,\RR)=\fsl(2,\RR)\oplus \RR I$ is
actually a direct sum of a simple Lie algebra and the one-dimensional
abelian algebra. In general, one can compute the Lie subalgebras of a
direct sum of algebras by means of the so-called Goursat lemma. This
result was introduced originally in~\cite{goursat} in order to compute
the subgroups of a product of two groups, but the result and its proof
are essentially the same in the category of Lie algebras (see, e.g.,
\cite[Section~3]{Figueroa-OFarrill:2015bsi}). In what follows, we
briefly describe the lemma.

Let $\g=\g_1\oplus \g_2$ be a direct sum of Lie algebras and consider the projections $\pi_i\colon \g\to \g_i$, which are surjective Lie algebra homomorphisms.
To each Lie subalgebra $\h\subseteq \g$, we associate the following pieces of data:
\begin{itemize}
  \item $\h_i=\pi_i(\h)$, which is a Lie subalgebra of $\g_i$ for each $i\in \{1,2\}$,
  \item $\q=\h\cap \g_2$, which is an ideal of $\h_2$,
  \item and the (well-defined) surjective Lie algebra homomorphism $f\colon \h_1\to \h_2 / \q$ given by letting
  \begin{equation}
    f(X)=Y+\q, \quad \text{where } (X,Y)\in\h.
  \end{equation}
\end{itemize}
\begin{lemma}[Goursat]\label{lemma:goursat}
  Assume the notation from above.
  The map $\h\mapsto (\h_1,\h_2,\q,f)$ establishes a bijective correspondence between the set of all Lie subalgebras of $\g$ and the set of all $4$-tuples $(\h_1,\h_2,\q,f)$ such that each $\h_i$ is a Lie subalgebra of $\g_i$, $\q$ is an ideal of $\h_2$ and $f\colon \h_1 \to \h_2/\q$ is a surjective homomorphism.

  The inverse map is given as follows:
  for a $4$-tuple $(\h_1,\h_2,\q,f)$ satisfying the aforementioned properties, the corresponding subalgebra is defined as
  \begin{equation}
    \h=\{(X,Y)\in\g\mid X\in \h_1, f(X)=Y+\q\}.
  \end{equation}
\end{lemma}

Since we are interested in determining subalgebras of $\g$ up to
conjugacy, we need to take into account how the action of an inner
automorphism of $\g$ affects the aforementioned algebraic data.  In
general, suppose that we are given a representation
$\rho\colon \mathsf{G}\to\operatorname{Inn}(\g)$, where $\mathsf{G}$
is a group and $\operatorname{Inn}(\g)$ is the group of inner
automorphisms of $\g$.  Because $\g_1$ and $\g_2$ are ideals of $\g$,
they are invariant under inner automorphisms, so $\rho$ restricts to
two subrepresentations $\rho_i\colon \mathsf{G}\to \Aut(\g_i)$.  Let
$\h$ be a subalgebra of $\g$ and $g\in \mathsf{G}$.  We consider the
conjugate subalgebra $\bar{\h}=\rho(g)\h$.  Then, for their
corresponding $4$-tuples $(\h_1,\h_2,\q,f)$ and
$(\bar{\h}_1,\bar{\h}_2,\bar{\q},\bar{f})$ we have
\begin{equation}\label{eq:equivariant-goursat-relations}
  \begin{aligned}
    \bar{\h}_i={}&\rho_i(g)\h_i, & \bar{\q}={}&\rho_2 (g)\q, & \bar{f}=\kappa(g)\circ f \circ \rho_1(g)^{-1},
  \end{aligned}
\end{equation}
where $\kappa(g)\colon \h_2/\q \to \bar{\h}_2/\bar{\q}$ sends $Y+\q\mapsto\rho_2(g)Y+\bar{\q}$.
As a consequence, the Lie subalgebras of $\g$ up to $\mathsf{G}$-conjugacy are in a bijective correspondence with the orbits of the action given by~\eqref{eq:equivariant-goursat-relations}.

We now focus on the case at hand.
Decompose $\g=\gl(2,\RR)=\fsl(2,\RR)\oplus \RR I$.
First we need to know the subalgebras of $\fsl(2,\RR)$ up to conjugacy.
This can be done directly by hand or by appealing to the structure of maximal solvable subalgebras of real semisimple Lie algebras (see for example~\cite{MR142687}).
We obtain that apart from the obvious subalgebras $0$ and $\fsl(2,\RR)$, these subalgebras are up to conjugacy
\begin{equation}
  \begin{aligned}
    \mathfrak{t}_0 (2,\RR) ={} & \vecspan \left\{
       \begin{pmatrix} 1 & 0 \\ 0 & -1 \end{pmatrix},
       \begin{pmatrix} 0 & 1 \\ 0 & 0 \end{pmatrix}
    \right\}, &
    \so(2)={}&\RR \begin{pmatrix} 0 & -1 \\ 1 & 0 \end{pmatrix}, \\
    \mathfrak{a}={}&\RR \begin{pmatrix} 1 & 0 \\ 0 & -1 \end{pmatrix}, &
    \mathfrak{n}={}& \RR \begin{pmatrix} 0 & 1 \\ 0 & 0 \end{pmatrix}.
  \end{aligned}
\end{equation}
On the other hand, given any $4$-tuple $(\h_1,\h_2,\q,f)$, the pair $(\h_2,\q)$ is either $(0,0)$, $(\RR I,\RR I)$ or $(\RR I,0)$.
Note that in the first two cases we have $\h_2/\q=0$, and thus $f$ is necessarily the zero map.
If $\q=\h_2=0$, then the corresponding subalgebra is $\h=\h_1$, whereas if $\q=\h_2=\RR I$, the corresponding subalgebra is $\h=\h_1\oplus \RR I$.
The subalgebras obtained by adding $\RR I$ to a nonzero proper subalgebra of $\fsl(2,\RR)$ are
\begin{equation}
  \begin{aligned}
      \mathfrak{t}(2,\RR)={}&\left\{
      \begin{pmatrix} x & y \\ 0 & z \end{pmatrix}
      ~\middle|~x,y,z\in \RR
    \right\}, &
    \mathbb{C} ={}& \vecspan \left\{
       \begin{pmatrix} 1 & 0 \\ 0 & 1 \end{pmatrix},
       \begin{pmatrix} 0 & -1 \\ 1 & 0 \end{pmatrix}
    \right\}, \\
       \mathfrak{diag}(2,\RR) = {} & \vecspan \left\{
       \begin{pmatrix} 1 & 0 \\ 0 & 0 \end{pmatrix},
       \begin{pmatrix} 0 & 0 \\ 0 & 1 \end{pmatrix}
    \right\}, &
    \mathfrak{n}\oplus\RR I = {} & \vecspan \left\{
       \begin{pmatrix} 1 & 0 \\ 0 & 1 \end{pmatrix},
       \begin{pmatrix} 0 & 1 \\ 0 & 0 \end{pmatrix}
    \right\}.
  \end{aligned}
\end{equation}

We now consider the case that $\h_2=\RR I$ and $\q=0$.
Because we need the existence of a surjective homomorphism $\h_1\to \h_2/\q=\RR I$, $\h_1$ is either $\mathfrak{t}_0 (2,\RR)$ or any of the one-dimensional subalgebras.

For the case $\h_1=\mathfrak{t}_0(2,\RR)$, a surjective homomorphism $f\colon \mathfrak{t}_0(2,\RR)\to \RR I$ has to vanish in the ideal $\mathfrak{n}\subseteq \mathfrak{t}_0(2,\RR)$, and we have
\begin{equation}
  f \begin{pmatrix} 1 & 0 \\ 0 & -1 \end{pmatrix} = y I, \quad y\in \RR\setminus \{0\}.
\end{equation}
The normaliser of $\mathfrak{t}_0(2,\RR)$ in $\GL(2,\RR)$ is the subgroup of all nonsingular triangular matrices, and it acts trivially on the space of surjective homomorphisms $\mathfrak{t}_0(2,\RR)\to \RR I$, so different values of $y$ yield nonconjugate subalgebras of the form
\begin{equation}
  \vecspan \left\{
     \begin{pmatrix} y+1 & 0 \\ 0 & y-1 \end{pmatrix},
     \begin{pmatrix} 0 & 1 \\ 0 & 0 \end{pmatrix}
  \right\}=
  \vecspan \left\{
     \begin{pmatrix} \cos \theta & 0 \\ 0 & \sin \theta \end{pmatrix},
     \begin{pmatrix} 0 & 1 \\ 0 & 0 \end{pmatrix}
  \right\}=\s_\theta,
\end{equation}
where $(\cos \theta,\sin\theta)=\tfrac{-1}{\sqrt{2(y^2+1)}}(y+1,y-1)$.
This establishes a bijective correspondence between $y\in \RR \setminus \{0\}$ and $\theta\in ( \tfrac{\pi}{4},\tfrac{5\pi}{4} )\setminus \{\tfrac{3\pi}{4}\}$, and because $\theta$ and $\theta+\pi$ yield the same subalgebra we can directly suppose that $\theta\in [0,\pi)\setminus \{\tfrac{\pi}{4},\tfrac{3\pi}{4}\}$.
Note that the case $\theta=\tfrac{\pi}{4}$ gives $\h=\mathfrak{n}\oplus \RR I$ and the case $\theta=\tfrac{3\pi}{4}$ gives $\h=\mathfrak{t}_{0}(2,\RR)$.

Suppose that $\h_1=\so(2)$.
The map $f$ satisfies
\begin{equation}
  f \begin{pmatrix} 0 & -1 \\ 1 & 0 \end{pmatrix}
  =y I, \quad y\in \RR \setminus \{0\}.
\end{equation}
The stabiliser of $\so(2)$ in $\GL(2,\RR)$ is
\begin{equation}
  \mathsf{CO}(2)\rtimes \mathbb{Z}_2=\left\{
     \begin{pmatrix} a & -b \\ b & a \end{pmatrix}~\middle|~ a^2 + b^2 \neq 0
  \right\} \rtimes \left< \begin{pmatrix} 1 & 0 \\ 0 & -1 \end{pmatrix} \right>,
\end{equation}
where $\mathsf{CO}(2)$ leaves $f$ unchanged and the outer $\mathbb{Z}_2$ allows us to change the sign of $f$.
Thus, a complete set of nonconjugate subalgebras is obtained by letting $y > 0$, and these are of the form
\begin{equation}
  \RR \begin{pmatrix} y & -1 \\ 1 & y \end{pmatrix}.
\end{equation}
Note that the algebra above in the limit case $y=0$ is precisely $\so(2)$.

Now, suppose $\h_1=\mathfrak{a}$.
A surjective map $f\colon \mathfrak{a}\to \RR I$ is of the form
\begin{equation}
  f \begin{pmatrix} 1 & 0 \\ 0 & -1 \end{pmatrix}=y I, \quad y\in \RR\setminus \{0\}.
\end{equation}
The stabiliser of $\mathfrak{a}$ is the semidirect product
\begin{equation}
  \mathsf{Diag}(2,\RR)\rtimes \mathbb{Z}_2 =
  \left\{
     \begin{pmatrix} a & 0 \\ 0 & b \end{pmatrix}
     ~\middle|~ab\neq 0
  \right\}\rtimes \left< \begin{pmatrix} 0 & 1 \\ 1 & 0 \end{pmatrix} \right>,
\end{equation}
where only the outer $\mathbb{Z}_2$ acts nontrivially on $f$ by changing its sign.
We conclude that a complete family of nonconjugate subalgebras is obtained by letting $y < 0$, giving us
\begin{equation}
  \RR \begin{pmatrix} y+1 & 0 \\ 0 & y-1 \end{pmatrix}=
  \RR \begin{pmatrix} x & 0 \\ 0 & 1 \end{pmatrix}, \quad x=\frac{y+1}{y-1}\in (-1,1).
\end{equation}
Different values of $x$ give rise to nonconjugate subalgebras.
Observe that the case $x=-1$ corresponds to $\mathfrak{a}$ whereas the case $x=1$ gives $\RR I$.

Finally, suppose that $\h_1=\mathfrak{n}$.
Then the map $f$ is such that
\begin{equation}
  f \begin{pmatrix} 0 & 1 \\ 0 & 0 \end{pmatrix} = y I, \quad y \neq 0.
\end{equation}
The normaliser of $\mathfrak{n}$ is $\mathsf{T}(2,\RR)$ and its action allows us to rescale $f$ freely, which means that we can simply set $y=1$ and obtain the subalgebra
\begin{equation}
  \RR \begin{pmatrix} 1 & 1 \\ 0 & 1 \end{pmatrix}.
\end{equation}

We can now collect our findings in the following result.

\begin{theorem}\label{thm:gl2Classification}
  Any nonzero proper subalgebra of $\gl(2,\RR)$ is conjugate to exactly one of
  the ones given in Table~\ref{tab:gl2-subalgebras}.
\end{theorem}

\begin{table}[h!]
  \centering
  \setlength{\extrarowheight}{3pt}
  \caption{Representatives of the conjugacy classes of proper subalgebras of $\gl(2,\RR)$.}
  \rowcolors{2}{blue!10}{white}
  \begin{tabular}{>{$}l<{$}>{$}c<{$}l}
    \toprule
    \multicolumn{1}{l}{\textbf{Subalgebra}} & \multicolumn{1}{c}{\textbf{Dimension}} & \textbf{Comments}\\
    \midrule[\heavyrulewidth]
    \fsl(2,\RR) & 3 & Simple \\
    \ft(2,\RR) & 3 & Solvable \\
    \fdiag(2,\RR) & 2 & Abelian \\
    \s_{\theta}=\vecspan \left\{ \begin{colorpmatrix} \cos \theta & 0 \\ 0 & \sin \theta \end{colorpmatrix},\quad \begin{colorpmatrix} 0 & 1 \\ 0 & 0 \end{colorpmatrix} \right\} & 2 & \makecell[l]{$\theta\in [0,\pi)$,\\ Abelian $\Leftrightarrow$ $\theta=\frac{\pi}{4}$} \\ 
    \CC \cong \RR I \oplus \so(2)& 2 & Abelian\\
    \RR \begin{colorpmatrix} x & 0 \\ 0 & 1 \end{colorpmatrix} & 1 & $|x|\leq 1$\\
    \RR \begin{colorpmatrix} 1 & 1 \\  0 & 1 \end{colorpmatrix} & 1 &  \\
    \RR \begin{colorpmatrix} 0 & 1 \\  0 & 0 \end{colorpmatrix} & 1 &  \\
    \RR \begin{colorpmatrix} y & -1 \\ 1 & y \end{colorpmatrix}& 1 & $y\geq 0$\\
    \bottomrule
  \end{tabular}
  \caption*{This table collects representatives of the conjugacy
    classes of proper Lie subalgebras of $\gl(2,\RR)$.  Among the
    one-dimensional subalgebras, $\so(1,1)$ corresponds to the first
    with $x=-1$ and $\so(2)$ corresponds to the
    last with $y=0$.}
  \label{tab:gl2-subalgebras}
\end{table}

Due to the structure theorems in Section~\ref{sec:struc-32-klas}, it is
important to determine which subalgebras of $\gl(2,\RR)$ admit
invariant a symmetric bilinear form $\alpha \in \sS^2(\RR^2)^*$ and an
invariant $\RR^2$-valued symmetric bilinear form $\sigma \in \RR^2
\otimes \sS^2(\RR^2)^*$.  From the determination of the relative
isomorphisms in Section~\ref{sec:relative-isomorphisms}, it is
moreover also important to determine which subalgebras of $\gl(2,\RR)$
leave a nonzero covector in $(\RR^2)^*$ invariant.

This latter condition is easy to see by inspection of the conjugacy
classes of subalgebras of $\gl(2,\RR)$ in
Table~\ref{tab:gl2-subalgebras}.  The result is tabulated in
Table~\ref{tab:gl2-subalgebras-with-W-dual-invariants}.

\begin{table}
  \centering
  \setlength{\extrarowheight}{3pt}
  \caption{Lie subalgebras of $\gl(2,\RR)$ that leave a linear form invariant.}
  \rowcolors{2}{blue!10}{white}
  \begin{tabular}{>{$}l<{$}>{$}l<{$}}
    \toprule
    \multicolumn{1}{l}{\textbf{Subalgebra} $\boldsymbol{\h}$} & \multicolumn{1}{c}{$\left((\RR^{2})^*\right)^\h$}\\
    \midrule[\heavyrulewidth]
    \s_{0} = \vecspan \left\{ \begin{colorpmatrix} 1 & 0 \\ 0 & 0 \end{colorpmatrix},\begin{colorpmatrix} 0 & 1 \\ 0 & 0 \end{colorpmatrix} \right\}& \RR\be^{2} \\
    \RR  \begin{colorpmatrix} 0 & 0 \\ 0 & 1 \end{colorpmatrix}& \RR\be^{1}\\
    \RR  \begin{colorpmatrix} 0 & 1 \\ 0 & 0 \end{colorpmatrix} & \RR\be^{2}\\
    \bottomrule
  \end{tabular}
  \caption*{This table lists the proper subalgebras in
    Table~\ref{tab:gl2-subalgebras} which preserve a linear form $\mu
    \in W^*$ invariant, along with the invariant linear forms.}
  \label{tab:gl2-subalgebras-with-W-dual-invariants}
\end{table}

Up to a change of basis and scaling, there are only three nontrivial
symmetric bilinear forms on $\RR ^{2}$, which are
\begin{equation}
  \be^{1}\otimes \be^{1} + \be^{2}\otimes \be^{2},\quad \be^{1}\otimes
  \be^{2} + \be^{2}\otimes \be^{1}, \quad \be^{2}\otimes \be^{2},
\end{equation}
where $\be_1,\be_2$ is a basis for $\RR^2$ and $\be^1,\be^2$ is the
corresponding dual basis.  The corresponding stabiliser subalgebras are
\begin{equation}
  \so(2)=\RR 
    \begin{pmatrix}
      0 & -1 \\
      1 & 0
    \end{pmatrix},\qquad
    \so(1,1)=\RR 
    \begin{pmatrix}
      1 & 0 \\
      0 & -1
    \end{pmatrix},\qquad
  \s_{0}=\vecspan
  \left\{ 
      \begin{pmatrix}
        1 & 0 \\
        0 & 0
      \end{pmatrix},
      \begin{pmatrix}
        0 & 1 \\
        0 & 0
      \end{pmatrix}
  \right\}.
\end{equation}
Therefore, we can also list the Lie subalgebras of $\gl(2,\RR)$ that admit
an invariant quadratic form in Table~\ref{tab:gl2-subalgs-inv-sbf}.

\begin{table}
  \centering
  \setlength{\extrarowheight}{3pt}
  \caption{Lie subalgebras of $\gl(2,\RR)$ that leave a quadratic form invariant.}
  \rowcolors{2}{blue!10}{white}
  \begin{tabular}{>{$}l<{$}>{$}l<{$}}
    \toprule
    \multicolumn{1}{l}{\textbf{Subalgebra} $\boldsymbol{\h}$} & \multicolumn{1}{c}{$\alpha \in \left(\sS^{2}(\RR^{2})^*\right)^\h$}\\
    \midrule[\heavyrulewidth]
    \s_{0} = \vecspan \left\{ \begin{colorpmatrix} 1 & 0 \\ 0 & 0 \end{colorpmatrix},\quad \begin{colorpmatrix} 0 & 1 \\ 0 & 0 \end{colorpmatrix} \right\}& \be^{2}\otimes \be^{2} \\
    \so(2) & \be^{1}\otimes \be^{1}+\be^{2}\otimes \be^{2}\\
    \so(1,1) & \be^{1}\otimes \be^{2} + \be^2 \otimes \be^1 \\
    \RR  \begin{colorpmatrix} 0 & 1 \\ 0 & 0 \end{colorpmatrix} & \be^{2}\otimes \be^{2} \\
    \RR  \begin{colorpmatrix} 0 & 0 \\ 0 & 1 \end{colorpmatrix}& \be^{1}\otimes \be^{1} \\
    \bottomrule
  \end{tabular}
  \caption*{This table lists the subalgebras in
    Table~\ref{tab:gl2-subalgebras} which preserve a quadratic form
    along with the invariant quadratic forms.  In the last line, we
    have had to change the form of the invariant quadratic form so
    that the subalgebra agrees with one in
    Table~\ref{tab:gl2-subalgebras} instead of merely being conjugate
    to one.}
  \label{tab:gl2-subalgs-inv-sbf}
\end{table}

The determination of which Lie subalgebras admit an invariant $\sigma
\colon \RR^2 \times \RR^2 \to \RR^2$ is a routine linear algebra
exercise and we only list the results in Table~\ref{tab:gl2-subalgs-inv-sigma}.

\begin{table}
  \centering
  \setlength{\extrarowheight}{3pt}
  \caption{Lie subalgebras of $\gl(2,\RR)$ that leave invariant a nonzero symmetric bilinear $\sigma\colon \RR^2 \times \RR^2 \to \RR^2$.}
  \resizebox{\linewidth}{!}{
  \rowcolors{2}{blue!10}{white} 
  \begin{tabular}{>{$}l<{$}>{$}l<{$}}
    \toprule
    \multicolumn{1}{l}{\textbf{Subalgebra} $\boldsymbol{\h}$}& \multicolumn{1}{c}{$\sigma \in (\RR^2 \otimes \sS^2 (\RR^2)^*)^\h$}\\
    \midrule[\heavyrulewidth]
    \s_{0} = \vecspan \left\{ \begin{colorpmatrix} 1 & 0 \\ 0 & 0 \end{colorpmatrix},\quad \begin{colorpmatrix} 0 & 1 \\ 0 & 0 \end{colorpmatrix} \right\}& \be_1 \otimes \left( \be^1 \otimes \be^2 + \be^2 \otimes \be^1 \right) + 2 \be_2 \otimes \be^2\otimes \be^2\\
    \s_{\arctan(1/2)} = \vecspan \left\{ \begin{colorpmatrix} 2 & 0 \\ 0 & 1 \end{colorpmatrix},\quad \begin{colorpmatrix} 0 & 1 \\ 0 & 0 \end{colorpmatrix} \right\}& \be_1 \otimes \be^2 \otimes \be^2\\
    \RR \begin{colorpmatrix} 0 & 0 \\ 0 & 1\end{colorpmatrix} & \be_1 \otimes \be^1 \otimes \be^1, \be_2 \otimes \left( \be^1 \otimes \be^2 + \be^2 \otimes \be^1 \right)\\
    \RR \begin{colorpmatrix} \frac12 & 0 \\ 0 & 1\end{colorpmatrix} & \be_2 \otimes \be^1 \otimes \be^1\\
    \RR \begin{colorpmatrix} 0 & 1 \\ 0 & 0\end{colorpmatrix} & \be_1 \otimes \be^2 \otimes \be^2, \be_1 \otimes \left( \be^1 \otimes \be^2 + \be^2 \otimes \be^1 \right) + 2 \be_2 \otimes \be^2\otimes \be^2\\
    \bottomrule
  \end{tabular}
} \label{tab:gl2-subalgs-inv-sigma}
\end{table}

\section{Real two-dimensional representations of Bianchi Lie algebras}
\label{sec:2d-reps-bianchi}

As discussed in Section~\ref{sec:struc-32-klas}, the Lie bracket of a $(3,2)$-KLA
$\g$ is partially encoded by a Bianchi Lie algebra $\fb$ and a two-dimensional
real representation $\rho\colon \fb \to \gl(W)$.  Furthermore, if
$\g'$ is another $(3,2)$-KLA (with associated representation
$\rho'\colon\fb'\to\gl(W')$) which differs from $\g$ via a relative
isomorphism, due to
Propositions~\ref{prop:generic-d-relative-isomorphisms} and~\ref{prop:3d-relative-isomorphisms}
there exists a Lie algebra
isomorphism $\psi\colon\fb\to\fb'$ and a linear isomorphism $h\colon
W\to W'$ such that $\rho'(\psi(B))=h\rho(B)h^{-1}$ for all $B\in\fb$.
This motivates the following definition.

Let $\fb$, $\fb'$ be two Lie algebras and $\rho\colon\fb\to \gl(W)$,
$\rho'\colon \fb'\to\gl(W')$ linear representations of $\fb$ and
$\fb'$ respectively. We say that $\rho$ and $\rho'$ are
\textit{weakly equivalent} or \textit{weakly isomorphic} if there
exists a linear isomorphism $h\colon W\to W'$ and a Lie algebra
isomorphism $\psi\colon \fb\to\fb'$ satisfying
\begin{equation}
  \rho'\circ \psi=\Ad(h)\circ \rho.
\end{equation}

Recall here that $\Ad(h)\colon \gl(W)\to\gl(W')$ is defined by
$\Ad(h)T=h\circ T \circ h^{-1}$.  Thus, for the classification of
$(3,2)$-KLAs it suffices to know the weak equivalence classes of
two-dimensional representations of Bianchi Lie algebras.  That will be
the purpose of this section.

Choosing a basis for $W$, the image $\im \rho$ is a Lie subalgebra of
$\gl(2,\RR)$ and the kernel $\ker \rho$ is an ideal of $\fb$.  For
 a Bianchi Lie algebra $\fb$ we let $\Aut(\fb)$ denote its Lie
group of automorphisms.  These are listed along with the
classification of Bianchi Lie algebras in Appendix~\ref{app:Bianchi}.

Having made these choices, we arrive at a rigidified notion of weak
equivalence. Indeed, having fixed a representative of each isomorphism
class of Bianchi Lie algebras: namely, the Lie algebras which are listed in
Table~\ref{tab:Bianchi}, $\psi$ now becomes an automorphism of a
Bianchi Lie algebra $\fb$ in that table.  Having chosen bases for $W$
and $W'$ and in that way identifying both $\gl(W)$ and $\gl(W')$ with
$\gl(2,\RR)$, now $h$ becomes an invertible linear transformation of
$\RR^2$.  In summary, we now have that two representations $\rho
\colon \fb \to \gl(2,\RR)$ and $\rho' \colon \fb \to \gl(2,\RR)$ are
weakly equivalent if there exist $\psi \in \Aut(\fb)$ and $h \in
\GL(2,\RR)$, so that $\rho'\circ \psi = \Ad(h) \circ \rho$.

In this section we will list all of the two-dimensional
representations $\rho \colon \fb \to \gl(2,\RR)$ up to weak
equivalence and for each such representation we will determine the
possible invariant $\alpha \in \sS^2 W^*$, $\sigma\colon \sS^2W \to W$
and $\varphi\colon \ext{2} W \to \fb$, as well as record the
calculation of the $\fb$-invariant vectors in $W^*$.  We will also
compute the subgroup $\Aut(\rho)\subseteq \GL(W)\times \Aut(\fb)$,
which we call the \textit{group of weak automorphisms} of $\rho$, that
stabilises $\rho$.  More precisely, we define
\begin{equation}
  \Aut(\rho)=\left\{
    (h,\psi)\in\GL(W)\times \Aut(\fb)\colon \rho\circ \psi=\Ad(h)\circ \rho
  \right\}.
\end{equation}

\subsection{Bianchi I}
\label{sec:B-I}

This is the abelian three-dimensional Lie algebra and hence $\Aut(\fb) =
\GL(\fb)$, which acts transitively on one- and two-dimensional vector
subspaces of $\fb$.  The image $\rho(\fb) \subseteq \gl(2,\RR)$ must be
an abelian Lie algebra and they can be read off from
Table~\ref{tab:gl2-subalgebras}: $\fdiag(2,\RR)$,
$\s_{\pi/4}$, $\CC \cong \RR I \oplus \so(2)$, any one-dimensional
subalgebra and the zero subalgebra.

\subsubsection{\emph{\Romanbar{I}.$1$}: $\im\rho = \fdiag(2,\RR)$}
\label{sec:BI-diag}

We may choose a basis $\{B_1,B_2,B_3\}$ for $\fb$ such that $\rho(B_3)=0$, whereas
\begin{equation}
  \rho(B_1) = \begin{pmatrix} 1 & 0 \\ 0 & 0 \end{pmatrix}\quad\text{and}\quad
  \rho(B_2) = \begin{pmatrix} 0 & 0 \\ 0 & 1 \end{pmatrix}.
\end{equation}
For this representation, we have $\alpha$ and $\sigma$ are zero, by
inspection of Tables~\ref{tab:gl2-subalgs-inv-sbf} and
\ref{tab:gl2-subalgs-inv-sigma}.  Similarly, because $\ext{2} W$ is a
nontrivial representation and $\fb$ is the trivial representation,
there cannot be any nonzero equivariant map $\ext{2} W \to \fb$, so
$\varphi=0$.  Finally, there are no nonzero $\fb$-invariants in $W^*$.
The group of weak automorphisms of $\rho$ is
\begin{equation}
  \begin{split}
    \Aut(\rho)=&\left\{
      \left(
        \begin{pmatrix}
          h_1 & 0 \\
          0 & h_2
        \end{pmatrix},
        \begin{pmatrix}
          1 & 0 & 0 \\
          0 & 1 & 0 \\
          a_{31} & a_{32} & a_{33}
        \end{pmatrix}
      \right) ~ \middle | ~ h_1h_2a_{33}\neq 0
    \right\}\\
    &\bigcup
    \left\{
      \left(
        \begin{pmatrix}
          0 & h_2 \\
          h_1 & 0
        \end{pmatrix},
        \begin{pmatrix}
          0 & 1 & 0 \\
          1 & 0 & 0 \\
          a_{31} & a_{32} & a_{33}
        \end{pmatrix}
      \right)~\middle|~ h_1h_2a_{33}\neq 0
    \right\}.
  \end{split}
\end{equation}

\subsubsection{\emph{\Romanbar{I}.$2$}: $\im\rho = \s_{\pi/4}$}
\label{sec:BI-s_pi/4}

We may choose a basis $\{B_1,B_2,B_3\}$ for $\fb$ such that $\rho(B_3)=0$, whereas
\begin{equation}
  \rho(B_1) = \begin{pmatrix} 1 & 0 \\ 0 & 1 \end{pmatrix} \quad\text{and}\quad
  \rho(B_2) = \begin{pmatrix} 0 & 1 \\ 0 & 0 \end{pmatrix}.
\end{equation}
The situation here is just like the previous case: $\alpha$, $\sigma$
and $(W^*)^\fb$ are all zero.  Finally, since $\ext{2} W$ is a
nontrivial representation, $\varphi$ is also zero.  We also have
\begin{equation}
  \Aut(\rho)=\left\{
    \left(
      \begin{pmatrix}
        a_{22}h_1 &  h_2 \\
        0 & h_1
      \end{pmatrix},
      \begin{pmatrix}
        1 & 0 & 0 \\
        0 & a_{22} & 0 \\
        a_{31} & a_{32} & a_{33}
      \end{pmatrix}
    \right)~\middle|~ h_{1}a_{22}a_{33}\neq 0
  \right\}.
\end{equation}

\subsubsection{\emph{\Romanbar{I}.$3$}: $\im\rho = \CC$}
\label{sec:BI-C}

We may choose a basis $\{B_1,B_2,B_3\}$ for $\fb$ such that $\rho(B_3) =0$, whereas
\begin{equation}
  \rho(B_1) = \begin{pmatrix} 1 & 0 \\ 0 & 1 \end{pmatrix} \quad\text{and}\quad
  \rho(B_2) = \begin{pmatrix} 0 & -1 \\ 1 & 0 \end{pmatrix}.
\end{equation}
We have again that $\varphi$ vanishes identically.  From the tables,
it follows that so do $\alpha$, $\sigma$ and $(W^*)^\fb$. The group of
weak automorphisms of $\rho$ is
\begin{equation}
  \Aut(\rho)=\left\{
    \left(
      \begin{pmatrix}
        h_1 & h_2 \\
        -\varepsilon h_2 & \varepsilon h_1
      \end{pmatrix},
      \begin{pmatrix}
        1 & 0 & 0 \\
        0 & \varepsilon & 0 \\
        a_{31} & a_{32} & a_{33}
      \end{pmatrix}
    \right)~\middle|~a_{33}(h_1^2+h_2^2)\neq 0, \varepsilon=\pm 1
  \right\}.
\end{equation}

\subsubsection{\emph{\Romanbar{I}.$4$}: $\im\rho = \RR \begin{pmatrix}x & 0 \\ 0 & 1\end{pmatrix}$ with $\lvert x \rvert \leq 1$}
\label{sec:BI-1d-x}

We may choose a basis $\{B_1,B_2,B_3\}$ for $\fb$ such that $\rho(B_2) =
\rho(B_3) = 0$ and
\begin{equation}
  \rho(B_1) = \begin{pmatrix} x & 0 \\ 0 & 1 \end{pmatrix},
\end{equation}
with $|x| \leq 1$.
One sees that different values of $x$ give rise to different
weak equivalence classes of representations.
Here all $\alpha,\sigma,\varphi, (W^*)^\fb$
are zero except at the following values of $x$:
\begin{itemize}
\item $x= -1$, which corresponds to $\so(1,1)$.  Here
  $\sigma = 0$, but $\alpha \in \RR (\be^1 \otimes \be^2 + \be^2 \otimes
  \be^1)$.  Furthermore, $\varphi(\be_1 \wedge \be_2)$ can be arbitrary.
\item $x=\frac12$.  Here $\sigma \in \RR(\be_2 \otimes \be^1
  \otimes \be^1)$, but $\alpha,\varphi,(W^*)^\fb$ are zero.
\item $x=0$.  Here $\sigma \in \RR ( \be_1 \otimes
    \be^1 \otimes \be^1 ) \oplus \RR ( \be_2 \otimes (\be^1
    \otimes \be^2 + \be^2 \otimes \be^1) )$, and $(W^*)^\fb = \RR \be^1$, $\alpha\in\RR (\be^1
  \otimes \be^1)$, but $\varphi$ is still zero.
\end{itemize}
The stabiliser $\Aut(\rho)$ of $\rho$ depends on $x$.
Indeed, we have:
\begin{itemize}
	\item If $x=1$, then
	\begin{equation}
		\Aut(\rho)=\left\{
			\left(
				h,
				\begin{pmatrix}
					1 & 0 \\
					\boldsymbol{u} & A
				\end{pmatrix}
			\right)~\middle|~ h,A\in\GL(2,\RR), \boldsymbol{u}\in\RR^2
		\right\}.
	\end{equation}
	\item If $x=-1$, then
	\begin{equation}
		\begin{split}
		\Aut(\rho)=&\left\{
		\left(
		\begin{pmatrix}
		h_{1} & 0 \\
		0 & h_{2} \\
		\end{pmatrix},
		\begin{pmatrix}
			1 & 0 \\
			\boldsymbol{u} & A
		\end{pmatrix}
		\right)~\middle|~ \boldsymbol{u}\in\RR^2, h_1h_2(\det A) \neq 0
		\right\}\\
		&\bigcup \left\{
		\left(
		\begin{pmatrix}
		0 & h_2 \\
		h_1 & 0 \\
		\end{pmatrix},
		\begin{pmatrix}
			-1 & 0 \\
			\boldsymbol{u} & A
		\end{pmatrix}
		\right)~\middle|~ \boldsymbol{u}\in\RR^2, h_1h_2(\det A) \neq 0
		\right\}.
		\end{split}
              \end{equation}
	\item If $x\neq\pm 1$, then
	\begin{equation}
		\Aut(\rho)=\left\{
		\left(
		\begin{pmatrix}
		h_{1} & 0 \\
		0 & h_{2} \\
		\end{pmatrix},
		\begin{pmatrix}
			1 & 0 \\
			\boldsymbol{u} & A
		\end{pmatrix}
		\right)~\middle|~ \boldsymbol{u}\in\RR^2, h_1h_2(\det A) \neq 0
		\right\}.
	\end{equation}
\end{itemize}

\subsubsection{\emph{\Romanbar{I}.$5$}: $\im\rho = \RR \begin{pmatrix}1 & 1 \\ 0 & 1 \end{pmatrix}$}
\label{sec:BI-1d-jordan-1}

We may choose a basis $\{B_1,B_2,B_3\}$ for $\fb$ such that
$\rho(B_2)=\rho(B_3)=0$ and
\begin{equation}
  \rho(B_1) = \begin{pmatrix} 1 & 1 \\ 0 & 1 \end{pmatrix}.
\end{equation}
Here all $\alpha,\sigma$ and $(W^*)^\fb$ are zero from the
tables and $\varphi$ is readily seen to be zero.
The group of weak automorphisms of $\rho$ is
\begin{equation}
	\Aut(\rho)=\left\{
	\left(
	\begin{pmatrix}
	h_1 & h_2 \\
	0 & h_1 \\
	\end{pmatrix},
	\begin{pmatrix}
	1 & 0 \\
	\boldsymbol{u} & A
	\end{pmatrix}
	\right)~\middle|~ h_1(\det A)\neq 0
	\right\}.
\end{equation}

\subsubsection{\emph{\Romanbar{I}.$6$}: $\im\rho = \RR \begin{pmatrix}0 & 1 \\ 0 & 0 \end{pmatrix}$}
\label{sec:BI-1d-jordan-0}

We may choose a basis $\{B_1,B_2,B_3\}$ for $\fb$ such that
$\rho(B_2)=\rho(B_3)=0$ and
\begin{equation}
  \rho(B_1) = \begin{pmatrix} 0 & 1 \\ 0 & 0 \end{pmatrix}.
\end{equation}
Here $\alpha \in \RR (\be^2 \otimes \be^2)$ and
\begin{equation}
  \sigma \in \RR \left(\be_1 \otimes \be^2 \otimes \be^2\right) \oplus \RR\left(\be_1 \otimes (\be^1 \otimes \be^2 + \be^2 \otimes \be^1) +
2 \be_2 \otimes \be^2 \otimes \be^2\right).
\end{equation}
Also, $\varphi$ is unconstrained.  From the tables, $(W^*)^\fb = \RR \be^2$.
One sees that the weak automorphism group of $\rho$ is
\begin{equation}
	\Aut(\rho)=\left\{
	\left(
	\begin{pmatrix}
	a h_1 & h_2 \\
	0 & h_1 \\
	\end{pmatrix},
	\begin{pmatrix}
	a & 0 \\
	\boldsymbol{u} & B
	\end{pmatrix}
	\right)~\middle|~ ah_1(\det B)\neq 0
	\right\}.
\end{equation}

\subsubsection{\emph{\Romanbar{I}.$7$}: $\im\rho = \RR \begin{pmatrix}y & -1 \\ 1 &
    y \end{pmatrix}$ with $y\geq 0$}
\label{sec:BI-1d-y}
We may choose a basis $\{B_1,B_2,B_3\}$ for $\fb$ such that
$\rho(B_2)=\rho(B_3)=0$ and
\begin{equation}
  \rho(B_1) = \begin{pmatrix} y & -1 \\ 1 & y\end{pmatrix}.
\end{equation}
Different values of $y$ give rise to nonconjugate subalgebras, and
thus to representations which are not weakly equivalent. Here $\sigma
= 0$ always and so are $\alpha$ and $\varphi$ 
unless $y =0$, which corresponds to $\so(2)$.  In that case, $\alpha
\in \RR (\be^1\otimes \be^1 + \be^2 \otimes \be^2)$ and similarly,
$\varphi(\be_1 \wedge \be_2)$ can be nonzero.  From the tables, we see
that $(W^*)^\fb = 0$. If $y=0$, the group of weak automorphisms of
$\rho$ is
\begin{equation}
  \Aut(\rho)=\left\{
    \left(
      \begin{pmatrix}
	h_{1} & h_{2} \\
	-\varepsilon h_{2} & \varepsilon h_{1} \\
      \end{pmatrix},
      \begin{pmatrix}
	\varepsilon & 0 \\
	\boldsymbol{u} & A 
      \end{pmatrix}
    \right)~\middle|~  (h_1^2+h_2^2)(\det A)\neq 0, \varepsilon=\pm 1
  \right\}.
\end{equation}
If $y \neq 0$, then we must take $\varepsilon = 1$ in the above formula.

\subsubsection{\emph{\Romanbar{I}.$8$}: $\im\rho = 0$}
\label{sec:BI-zero}

Here both the representation $W$ and the adjoint representation $\fb$
are trivial and hence $\alpha$, $\sigma$ and $\varphi$ are unconstrained.
Clearly, $(W^*)^\fb=W^*$.  Trivially, we have
\begin{equation}
	\Aut(\rho)=\GL(W)\times\Aut(\fb).
\end{equation}

\subsubsection{Summary}
\label{sec:BI-summary}

We summarise the results of this section in Table~\ref{tab:BI-2d-reps-summary}.

\begin{table}[h!]
  \centering
  \setlength{\extrarowheight}{3pt}
  \caption{Two-dimensional real representations of Bianchi $\Romanbar{I}$.}
  \rowcolors{2}{blue!10}{white}
  \begin{tabular}{>{$}l<{$}*{3}{>{$}c<{$}}>{$}c<{$}*{4}{>{$}c<{$}}}
    \toprule
    \multicolumn{1}{c}{Label} & \multicolumn{1}{c}{$\rho(B_1)$} & \multicolumn{1}{c}{$\rho(B_2)$} & \multicolumn{1}{c}{$\rho(B_3)$} & \multicolumn{1}{c}{Remarks} & \multicolumn{1}{c}{$\alpha$} & \multicolumn{1}{c}{$\sigma$} & \multicolumn{1}{c}{$\varphi$} & \multicolumn{1}{c}{$(W^*)^\fb$}\\
    \midrule[\heavyrulewidth]
    \hypertarget{I1}{\Romanbar{I}.1}& \begin{colorpmatrix} 1 & 0\\0 & 0\end{colorpmatrix} & \begin{colorpmatrix} 0 & 0\\ 0 & 1\end{colorpmatrix} &  0 & & & & & \\
    \hypertarget{I2}{\Romanbar{I}.2} & \begin{colorpmatrix} 1 & 0\\0 & 1\end{colorpmatrix} & \begin{colorpmatrix} 0 & 1\\ 0 & 0\end{colorpmatrix} &  0 & & & & & \\
    \hypertarget{I3}{\Romanbar{I}.3} & \begin{colorpmatrix} 1 & 0\\0 & 1\end{colorpmatrix} & \begin{colorpmatrix} 0 & -1\\ 1 & 0\end{colorpmatrix} &  0 & & & & & \\
    \hypertarget{I4}{\Romanbar{I}.4_x} & \begin{colorpmatrix} x & 0\\0 & 1\end{colorpmatrix} & 0 & 0 & |x| \leq 1 & x=-1,0 & x=0,\frac12 & x=-1 & \RR\be^1~(x=0)\\
    \hypertarget{I5}{\Romanbar{I}.5} & \begin{colorpmatrix} 1 & 1\\0 & 1\end{colorpmatrix} & 0 & 0 & & & & & \\
    \hypertarget{I6}{\Romanbar{I}.6} & \begin{colorpmatrix} 0 & 1\\0 & 0\end{colorpmatrix} & 0 & 0 & & \checkmark & \checkmark & \checkmark & \RR\be^2\\
    \hypertarget{I7}{\Romanbar{I}.7_y} & \begin{colorpmatrix} y & -1\\1 & y\end{colorpmatrix} & 0 & 0 & y\geq 0 & y=0 & & y=0 & \\
    \hypertarget{I8}{\Romanbar{I}.8} & 0 & 0 & 0 & & \checkmark & \checkmark & \checkmark & W^*\\
    \bottomrule
  \end{tabular}
  \caption*{Using the basis $\{B_1,B_2,B_3\}$ in Table~\ref{tab:Bianchi}
    for the Bianchi $\Romanbar{I}$ Lie algebra $\fb$, each row is a
    representation $\rho\colon \fb \to \gl(2,\RR)$.  The representation
    has some parameters which in the absence of any ``Remarks'' take any
    possible real values.  We also indicate for each representation
    whether there are any invariant $\alpha \in \sS^2 (\RR^2)^*$,
    $\sigma \in \RR^2 \otimes \sS^2 (\RR^2)^*$ and $\varphi \in \fb
    \otimes  \ext{2} (\RR^2)^*$.  A $\checkmark$ indicates existence
    of a nonzero such invariant for any (allowed) value of the
    parameters,  whereas a blank cell indicates its absence.  Any
    other expression indicates the values of any parameter where a
    nonzero invariant exists.  The last column is the $\fb$-invariant
    subspace of $W^*$: it is assumed to be zero unless otherwise noted.}
  \label{tab:BI-2d-reps-summary}
\end{table}

\subsection{Bianchi II}
\label{sec:B-II}

This is the Heisenberg Lie algebra with brackets $[B_2,B_3]=B_1$
relative to a choice of basis $\{B_1,B_2,B_3\}$.
The automorphism group is given by
\begin{equation}
  \Aut(\fb) = \left\{
    \begin{pmatrix}
      \det A & \boldsymbol{a}^T \\ \boldsymbol{0} & A
    \end{pmatrix}~ \middle | ~  A \in \GL(2,\RR),~ \boldsymbol{a} \in \RR^2\right\}.
\end{equation}
The proper ideals are the first derived ideal $[\fb,\fb]= \RR B_1$ and
the span of $B_1$ together with any nonzero vector in the $(B_2,B_3)$-plane.
Since $\Aut(\fb)$ acts transitively on such nonzero vectors, we are
free to only consider the ideal spanned by $\{B_1,B_2\}$.

It follows from Table~\ref{tab:gl2-subalgebras} that there are no
nilpotent nonabelian Lie subalgebras of $\gl(2,\RR)$, hence $\ker \rho
\neq 0$.  We discuss in turn the cases $\ker\rho = \fb$, $\ker\rho =
\RR B_1$ and $\ker\rho = \vecspan\{B_1,B_2\}$.
For $\ker \rho = \RR B_1$, because $\fb/\RR B_1$ is abelian,
$\im\rho$ can be any two-dimensional abelian subalgebra of
$\gl(2,\RR)$, so it is conjugate to one of $\fdiag(2,\RR)$,
$\s_{\pi/4}$ and $\CC$.

\subsubsection{\emph{\Romanbar{II}}.$1$: $\ker\rho = \RR B_1$ and $\im \rho = \fdiag(2,\RR)$}
\label{sec:BII-1d-ker-im-rho-diag}

We can choose bases for $W$ and $\fb$ such that $\rho(B_1)=0$, whereas
\begin{equation}
  \rho(B_2) = \begin{pmatrix} 1 & 0 \\ 0 & 0 \end{pmatrix} \quad\text{and}\quad
  \rho(B_3) = \begin{pmatrix} 0 & 0 \\ 0 & 1 \end{pmatrix}.
\end{equation}

From the tables in Section~\ref{sec:gl2-subalgebras}, it follows that
$\alpha = 0$, $\sigma = 0$ and $(W^*)^\fb=0$.
By direct computation, one sees that $\varphi = 0$.
The group of weak automorphisms of $\rho$ is
\begin{equation}
	\begin{split}
	\Aut(\rho)=&\left\{
	\left(
	\begin{pmatrix}
	h_{1} & 0 \\
	0 & h_{2} \\
	\end{pmatrix},
	\begin{pmatrix}
	1 & a_1 & a_2 \\
	0 & 1 & 0 \\
	0 & 0 & 1 \\
	\end{pmatrix}
	\right)~\middle|~ h_1h_2\neq 0
	\right\}\\
	&\bigcup \left\{
	\left(
	\begin{pmatrix}
	0 & h_2 \\
	h_1 & 0 \\
	\end{pmatrix},
	\begin{pmatrix}
	-1 & a_1 & a_2 \\
	0 & 0 & 1 \\
	0 & 1 & 0 \\
	\end{pmatrix}
	\right)~\middle|~ h_1h_2\neq 0
	\right\}.
	\end{split}
\end{equation}

\subsubsection{\emph{\Romanbar{II}.$2$}: $\ker\rho = \RR B_1$ and $\im \rho = \s_{\pi/4}$}
\label{sec:BII-1d-ker-im-rho-s-pi-by-4}

We can choose bases for $W$ and $\fb$ such that $\rho(B_1)=0$, whereas
\begin{equation}
  \rho(B_2) = \begin{pmatrix} 1 & 0 \\ 0 & 1 \end{pmatrix} \quad\text{and}\quad
  \rho(B_3) = \begin{pmatrix} 0 & 1 \\ 0 & 0 \end{pmatrix}.
\end{equation}

From the tables in Section~\ref{sec:gl2-subalgebras}, it follows that
$\alpha = 0$, $\sigma = 0$ and $(W^*)^\fb=0$.
It is also readily seen that $\varphi = 0$.
We also have
\begin{equation}
	\Aut(\rho)=\left\{
	\left(
	\begin{pmatrix}
	d h_{1} & h_{2} \\
	0 & h_{1} \\
	\end{pmatrix},
	\begin{pmatrix}
	d & a_1 & a_2 \\
	0 & 1 & 0 \\
	0 & 0 & d \\
	\end{pmatrix}
	\right)~\middle|~ dh_1\neq 0
	\right\}.
\end{equation}

\subsubsection{\emph{\Romanbar{II}.$3$}: $\ker\rho = \RR B_1$ and $\im \rho = \CC$}
\label{sec:BII-1d-ker-im-rho-C}

We can choose bases for $W$ and $\fb$ such that $\rho(B_1)=0$, whereas
\begin{equation}
  \rho(B_2) = \begin{pmatrix} 1 & 0 \\ 0 & 1 \end{pmatrix} \quad\text{and}\quad
  \rho(B_3) = \begin{pmatrix} 0 & -1 \\ 1 & 0 \end{pmatrix}.
\end{equation}

This is very similar to the previous case, since $\rho(B_2)$ takes the
same form.  All of $\alpha$, $\sigma$, $\varphi$ and $(W^*)^\fb$ are zero.
We also see that
\begin{equation}
	\Aut(\rho)=\left\{
	\left(
	\begin{pmatrix}
	h_1 & h_2 \\
	-\varepsilon h_2 & \varepsilon h_1 \\
	\end{pmatrix},
	\begin{pmatrix}
	\varepsilon & a_1 & a_2 \\
	0 & 1 & 0 \\
	0 & 0 & \varepsilon \\
	\end{pmatrix}
	\right)~\middle|~ h_1^2+h_2^2\neq 0, \varepsilon=\pm 1
	\right\}.
\end{equation}

\subsubsection{\emph{\Romanbar{II}.$4$}: $\ker\rho = \vecspan\{B_1,B_2\}$ and $\im \rho =  \RR \begin{pmatrix} x & 0 \\ 0 & 1 \end{pmatrix}$  with $|x| \leq 1$}
\label{sec:BII-2d-ker-im-rho-x}

We can choose bases for $W$ and $\fb$ such that
$\rho(B_1)=\rho(B_2)=0$ and
\begin{equation}
  \rho(B_3) = \begin{pmatrix} x & 0 \\ 0 & 1 \end{pmatrix},
\end{equation}
with $|x| \leq 1$.
Different values of $x$ give rise to different weak isoclasses of representations.

All of $\alpha,\sigma,\varphi,(W^*)^\fb$ are zero except for the
following values of $x$:
\begin{itemize}
\item $x=0$.  Here $\alpha\in\RR (\be^{1}\otimes \be^{1})$ and
  $\sigma \in \vecspan \{\be_1 \otimes \be^1 \otimes \be^1, \be_2
  \otimes (\be^1 \otimes \be^2 + \be^2 \otimes \be^1)\}$.  Also
  $\varphi$ is zero.  Finally, $(W^*)^\fb = \RR \be^1$.
\item $x=\frac12$.  Here $\sigma \in \RR (\be_2 \otimes \be^1
  \otimes \be^1)$ and all other are zero.
\item $x=-1$, which is $\so(1,1)$.  Here $\sigma =
  0$ and $(W^*)^\fb=0$, but $\alpha \in \RR(\be^1 \otimes \be^2 +
  \be^2 \otimes \be^1)$ and $\varphi(\be_1 \wedge \be_2)   \in \RR B_1$.
\end{itemize}
As for the group of weak automorphisms, we have
\begin{itemize}
	\item If $x=1$, one sees that
	\begin{equation}
		\Aut(\rho)=\left\{
		\left(
		h,
		\begin{pmatrix}
		a & b_1 & b_2\\
		0 & a & c \\
		0 & 0 & 1 \\
		\end{pmatrix}
		\right)~\middle|~ a(\det h)\neq 0
		\right\}.
	\end{equation}
	\item If $x=-1$, we have
	\begin{equation}
		\begin{split}
			\Aut(\rho)=&\left\{
			\left(
			\begin{pmatrix}
			h_{1} & 0 \\
			0 & h_{2} \\
			\end{pmatrix},
			\begin{pmatrix}
			a & b_1 & b_2 \\
			0 & a & c \\
			0 & 0 & 1 \\
			\end{pmatrix}
			\right)~\middle|~ ah_1h_2\neq 0
			\right\}\\
			&\bigcup
			\left\{
			\left(
			\begin{pmatrix}
			0 & h_{2} \\
			h_{1} & 0 \\
			\end{pmatrix},
			\begin{pmatrix}
			-a & b_1 & b_2 \\
			0 & a & c \\
			0 & 0 & -1 \\
			\end{pmatrix}
			\right)~\middle|~ ah_1h_2\neq 0
			\right\}.
		\end{split}
	\end{equation}
	\item If $x\neq \pm 1$, we obtain
	\begin{equation}
		\Aut(\rho)=\left\{
			\left(
			\begin{pmatrix}
			h_{1} & 0 \\
			0 & h_{2} \\
			\end{pmatrix},
			\begin{pmatrix}
			a & b_1 & b_2 \\
			0 & a & c \\
			0 & 0 & 1 \\
			\end{pmatrix}
			\right)~\middle|~ ah_1h_2\neq 0
			\right\}.
	\end{equation}
\end{itemize}

\subsubsection{\emph{\Romanbar{II}.$5$}: $\ker\rho = \vecspan\{B_1,B_2\}$ and $\im \rho =  \RR \begin{pmatrix} 1 & 1 \\ 0 & 1 \end{pmatrix}$}
\label{sec:BII-2d-ker-im-rho-jordan-1}

We can choose bases for $W$ and $\fb$ such that $\rho(B_1)=\rho(B_2)=0$ and
\begin{equation}
  \rho(B_3) = \begin{pmatrix} 1 & 1 \\ 0 & 1 \end{pmatrix}.
\end{equation}

From the tables, $\alpha$, $\sigma$ and $(W^*)^\fb$ are all
zero.
Moreover, $\varphi$ vanishes as well.
The stabiliser subgroup of $\rho$ is
\begin{equation}
	\Aut(\rho)=\left\{
	\left(
	\begin{pmatrix}
	h_1 & h_2 \\
	0 & h_1 \\
	\end{pmatrix},
	\begin{pmatrix}
	a & b_1 & b_2 \\
	0 & a & c \\
	0 & 0 & 1 \\
	\end{pmatrix}
	\right)~\middle|~ a h_1\neq 0
	\right\}.
\end{equation}

\subsubsection{\emph{\Romanbar{II}.$6$}: $\ker\rho = \vecspan\{B_1,B_2\}$ and $\im \rho = \RR \begin{pmatrix} 0 & 1 \\ 0 & 0 \end{pmatrix}$}
\label{sec:BII-2d-ker-im-rho-jordan-0}

We can choose bases for $W$ and $\fb$ such that $\rho(B_1)=\rho(B_2)=0$ and
\begin{equation}
  \rho(B_3) = \begin{pmatrix} 0 & 1 \\ 0 & 0 \end{pmatrix}.
\end{equation}

From the tables, $\alpha \in \RR(\be^2 \otimes \be^2)$, $\sigma \in
\vecspan\{\be_1 \otimes \be^2 \otimes \be^2, \be_1 \otimes (\be^1
\otimes \be^2 + \be^2 \otimes \be^1) + 2 \be_2 \otimes \be^2 \otimes
\be^2\}$ and $(W^*)^\fb = \RR\be^2$.  It follows that $\varphi(\be_1 \wedge \be_2) \in \RR B_1$.
We also have
\begin{equation}
	\Aut(\rho)=\left\{
	\left(
	\begin{pmatrix}
	d h_1 & h_2 \\
	0 & h_1 \\
	\end{pmatrix},
	\begin{pmatrix}
	a d & b_1 & b_2 \\
	0 & a & c \\
	0 & 0 & d \\
	\end{pmatrix}
	\right)~\middle|~ adh_1\neq 0
	\right\}.
\end{equation}

\subsubsection{\emph{\Romanbar{II}.$7$}: $\ker\rho = \vecspan\{B_1,B_2\}$ and $\im \rho =
  \RR \begin{pmatrix} y & -1 \\ 1 & y \end{pmatrix}$ with $y\geq 0$}
\label{sec:BII-2d-ker-im-rho-y}
Here we can choose bases for $W$ and $\fb$ such that $\rho(B_1)=\rho(B_2)=0$ and
\begin{equation}
  \rho(B_3) = \begin{pmatrix} y & -1 \\ 1 & y \end{pmatrix}
\end{equation}
for some $y\geq 0$. Furthermore, different values of $y$ yield
different weak isoclasses of representations.

Here $\sigma =0$ always and so are $\alpha,\varphi$ unless
$y=0$, when $\alpha \in \RR(\be^1 \otimes \be^1 + \be^2 \otimes
\be^2)$, and $\varphi(\be_1 \wedge \be_2) \in \RR B_1$.  From the
tables we read off that $(W^*)^\fb = 0$. As for the group of weak
automorphisms, we see that if $y=0$, it satisfies
\begin{equation}
	\begin{split}
		\Aut(\rho)=&\left\{
		\left(
		\begin{pmatrix}
		h_1 & h_2 \\
		-\varepsilon h_2 & \varepsilon h_1 \\
		\end{pmatrix},
		\begin{pmatrix}
		\varepsilon a & b_1 & b_2 \\
		0 & a & c \\
		0 & 0 & \varepsilon \\
		\end{pmatrix}
		\right)~\middle|~ a(h_1^2+h_2^2)\neq 0,\varepsilon=\pm 1
		\right\},
	\end{split}
\end{equation}
while for $y\neq 0$ we have
\begin{equation}
	\begin{split}
		\Aut(\rho)=&\left\{
		\left(
		\begin{pmatrix}
		h_1 & h_2 \\
		-h_2 & h_1 \\
		\end{pmatrix},
		\begin{pmatrix}
		a & b_1 & b_2 \\
		0 & a & c \\
		0 & 0 & 1 \\
		\end{pmatrix}
		\right)~\middle|~ a(h_1^2+h_2^2)\neq 0
		\right\}.
	\end{split}
\end{equation}

\subsubsection{\emph{\Romanbar{II}.$8$}: $\ker\rho = \fb$}
\label{sec:BII-rho-zero}

Here $\alpha,\sigma$ are not constrained at all, whereas
$\im\varphi \in \RR B_1$ and all of $W^*$ is invariant.
We trivially have
\begin{equation}
	\Aut(\rho)=\GL(W)\times\Aut(\fb).
\end{equation}
We will analyse this case more carefully when classifying the
corresponding $(3,2)$-KLAs.

\subsubsection{Summary}
\label{sec:BII-summary}

We summarise the results of this section in Table~\ref{tab:BII-2d-reps-summary}, with the same notation as that explained in Table~\ref{tab:BI-2d-reps-summary}.

\begin{table}[h!]
  \centering
  \setlength{\extrarowheight}{3pt}
  \caption{Two-dimensional real representations of Bianchi $\Romanbar{II}$.}
  \rowcolors{2}{blue!10}{white}
  \begin{tabular}{>{$}l<{$}*{3}{>{$}c<{$}}>{$}c<{$}*{4}{>{$}c<{$}}}
    \toprule
    \multicolumn{1}{c}{Label} & \multicolumn{1}{c}{$\rho(B_1)$} & \multicolumn{1}{c}{$\rho(B_2)$} & \multicolumn{1}{c}{$\rho(B_3)$} & \multicolumn{1}{c}{Remarks} & \multicolumn{1}{c}{$\alpha$} & \multicolumn{1}{c}{$\sigma$} & \multicolumn{1}{c}{$\varphi$} & \multicolumn{1}{c}{$(W^*)^\fb$}\\
    \midrule[\heavyrulewidth]
    \hypertarget{II1}{\Romanbar{II}.1} & 0 & \begin{colorpmatrix} 1 & 0\\0 & 0\end{colorpmatrix} & \begin{colorpmatrix} 0 & 0\\ 0 & 1\end{colorpmatrix} & & & & & \\
    \hypertarget{II2}{\Romanbar{II}.2} & 0 & \begin{colorpmatrix} 1 & 0\\0 & 1\end{colorpmatrix} & \begin{colorpmatrix} 0 & 1\\ 0 & 0\end{colorpmatrix} & & & & & \\
    \hypertarget{II3}{\Romanbar{II}.3} & 0 & \begin{colorpmatrix} 1 & 0\\0 & 1\end{colorpmatrix} & \begin{colorpmatrix} 0 & -1\\ 1 & 0\end{colorpmatrix} & & & & & \\
    \hypertarget{II4}{\Romanbar{II}.4_x} & 0 & 0 & \begin{colorpmatrix} x & 0\\0 & 1\end{colorpmatrix} & |x| \leq 1 & x=-1,0 & x=0,\frac12 & x=-1 & \RR\be^1~(x=0)\\
    \hypertarget{II5}{\Romanbar{II}.5} & 0 & 0 & \begin{colorpmatrix} 1 & 1\\0 & 1\end{colorpmatrix} & & & & & \\
    \hypertarget{II6}{\Romanbar{II}.6} & 0 & 0 & \begin{colorpmatrix} 0 & 1\\0 & 0\end{colorpmatrix} & & \checkmark & \checkmark & \checkmark & \RR\be^2 \\
    \hypertarget{II7}{\Romanbar{II}.7_y} & 0 & 0 & \begin{colorpmatrix} y & -1\\1 & y\end{colorpmatrix} & y\geq 0 & y=0 & & y=0 & \\
    \hypertarget{II8}{\Romanbar{II}.8} & 0 & 0 & 0 & & \checkmark & \checkmark & \checkmark & W^*\\
    \bottomrule
  \end{tabular}
  \label{tab:BII-2d-reps-summary}
\end{table}

\subsection{Bianchi III}
\label{sec:B-III}

This is the unique nonabelian decomposable Bianchi Lie algebra with
$B_1$ central and $[B_2,B_3]=2B_2$.  It is the special case of
Bianchi~$\Romanbar{VI}_c$ with $c=1$, but it is convenient to treat it
separately and then simply assume that $c \neq 1$ when we treat
Bianchi~$\Romanbar{VI}_c$ in Section~\ref{sec:B-VI}.  The automorphism
group is given by
\begin{equation}
  \label{eq:aut-BIII}
  \Aut(\fb) = \left\{\begin{pmatrix} a & 0 & u\\ 0 & b & v\\ 0 & 0 &  1 \end{pmatrix} ~\middle |~ a b \neq 0\right\}.
\end{equation}
The proper ideals are $\RR B_1$, $\RR B_2$ and $\vecspan\{B_2, \sin\theta
B_1 + \cos\theta B_3\}$ for $\theta \in [0,\pi)$.  Up to the action of
the automorphisms, there are two cases for the two-dimensional ideal:
if $\theta \neq \frac\pi2$, we may shift
$B_3 \mapsto B_3 - \tan \theta B_1$ to arrive at the ideal
$\vecspan\{B_2,B_3\}$; whereas if $\theta = \frac\pi2$, then we just
have the ideal $\vecspan\{B_1,B_2\}$.

The Lie algebra $\fb$ is actually isomorphic to $\ft(2,\RR)$.  For instance, an
isomorphism $\iota \colon \fb \to \ft(2,\RR)$ is given by
\begin{equation}
  \label{eq:iso-BIII-t2R}
  \iota(B_1) = \begin{pmatrix}1 & 0 \\ 0 & 1 \end{pmatrix}, \qquad
  \iota(B_2) = \begin{pmatrix}0 & 1 \\ 0 & 0 \end{pmatrix} \qquad\text{and}\qquad
  \iota(B_3) = \begin{pmatrix}-1 & 0 \\ 0 & 1 \end{pmatrix}.
\end{equation}
This means that the possible kernels for representations $\rho\colon \fb
\to \gl(W)$ are $0$, $\RR B_1$, $\RR B_2$, $\vecspan\{B_1,B_2\}$,
$\vecspan\{B_2,B_3\}$ and $\fb$.
Note that in the case $\ker \rho = \RR B_2$, the quotient $\fb/\RR
B_2$ is abelian, so the image of $\rho$ is conjugate to either
$\s_{\pi/4}$, $\fdiag(2,\RR)$ or $\CC$.

It follows from the classification of Lie subalgebras of $\gl(2,\RR)$
in Table~\ref{tab:gl2-subalgebras}, that the only other Bianchi Lie
algebra besides Bianchi~$\Romanbar{III}$ which admits a faithful
two-dimensional representation is Bianchi~$\Romanbar{VIII}$, which is
isomorphic to $\fsl(2,\RR)$.

\subsubsection{\emph{\Romanbar{III}.$1$}: $\ker\rho =0$}
\label{sec:BIII-zero-ker}

In this case, $\rho$ is an isomorphism $\fb \to \ft(2,\RR)$.  It
follows from the tables that $\alpha, \sigma$ and $(W^*)^\fb$
are all zero.  Up to an automorphism of $\fb$, we can take
$\rho = \iota$, the isomorphism
$\iota \colon \fb \stackrel{\cong}{\to} \ft(2,\RR)$ in
equation~\eqref{eq:iso-BIII-t2R}.
It follows from direct computation that $\varphi = 0$.
The weak automorphisms of
$\rho$ are given by
\begin{equation}
	\Aut(\rho)=\left\{
	\left(
	\begin{pmatrix}
	b t & \frac{t v}{2} \\
	0 & t \\
	\end{pmatrix},
	\begin{pmatrix}
	1 & 0 & 0 \\
	0 & b & v \\
	0 & 0 & 1 \\
	\end{pmatrix}
	\right)~\middle|~ bt\neq 0
	\right\}.
\end{equation}

\subsubsection{\emph{\Romanbar{III}.$2$}: $\ker\rho =\RR B_1$}
\label{sec:BIII-ker-B1}

In this case we have (with respect to an appropriate basis of $W$) $\im \rho = \s_\theta$ for $\theta \neq \frac\pi4$,
since $\fb/\RR B_1$ is nonabelian.  We have that $\rho(B_1) = 0$ and
\begin{equation}
  \rho(B_2) = \begin{pmatrix} 0 & r \\ 0 & 0 \end{pmatrix}
  \qquad\text{and}\qquad
  \rho(B_3) = \begin{pmatrix} x & s \\ 0 & y \end{pmatrix},
\end{equation}
where, in particular, $r \neq 0$.  Imposing $[\rho(B_2),\rho(B_3)] = 2
\rho(B_2)$ says that $y = x + 2$. We can rescale $B_2 \mapsto \frac1r
B_2$ to put $r=1$ and we may shift $B_3 \mapsto B_3 - s B_2$ to set
$s =0$.  In summary, we may assume
\begin{equation}
  \rho(B_1) =0, \qquad \rho(B_2) = \begin{pmatrix} 0 & 1 \\ 0 & 0 \end{pmatrix}
  \qquad\text{and}\qquad
  \rho(B_3) = \begin{pmatrix} x & 0 \\ 0 & x+2 \end{pmatrix},
\end{equation}
for any $x \in \RR$ (and different values of $x$ give different isoclasses).

From the tables, we have that $\alpha,\sigma,(W^*)^\fb$ are
zero except in the following cases:
\begin{itemize}
\item $x=-2$ ($\theta = 0$) when $\alpha \in \RR (\be^2 \otimes
  \be^2)$, $\sigma \in \RR(\be_1\otimes (\be^1\otimes\be^2+\be^2\otimes\be^1)+2\be_2\otimes\be^2\otimes\be^2)$ and $(W^*)^\fb = \RR\be^2$;
\item $x=-4$ $(\cot\theta = 2)$, when $\sigma \in \RR(\be_1 \otimes \be^2 \otimes \be^2)$.
\end{itemize}

As for the map $\varphi$, we have:
\begin{itemize}
\item $x=-2$, in which case $\varphi(\be_1 \wedge \be_2)
  \in \RR B_2$;
\item or $x=-1$, in which case $\varphi(\be_1 \wedge \be_2)
  \in \RR B_1$.
\end{itemize}

We have
\begin{equation}
	\Aut(\rho)=\left\{
	\left(
	\begin{pmatrix}
	b t & \frac{t v}{2} \\
	0 & t \\
	\end{pmatrix},
	\begin{pmatrix}
	a & 0 & u \\
	0 & b & v \\
	0 & 0 & 1 \\
	\end{pmatrix}
	\right)~\middle|~ abt\neq 0
	\right\}.
\end{equation}

\subsubsection{\emph{\Romanbar{III}.$3$}: $\ker\rho = \RR B_2$ and $\im\rho = \fdiag(2,\RR)$}
\label{sec:BIII-ker-B2-im-diag}

Here $\rho(B_2)=0$ and
\begin{equation}
  \rho(B_1) = \begin{pmatrix} x & 0 \\ 0 & y \end{pmatrix}
  \qquad\text{and}\qquad \rho(B_3) = \begin{pmatrix} s & 0 \\ 0 &
    t \end{pmatrix}.
\end{equation}
Without loss of generality we can assume that $x \neq 0$, since we can
always swap the basis vectors in $W$.

In that case, we may rescale $B_1$ so that $x = 1$ and we may shift
$B_3$ so that $s = 0$, resulting in
\begin{equation}
  \label{eq:BIII-ker-B2-im-diag}
  \rho(B_1) = \begin{pmatrix} 1 & 0 \\ 0 & y \end{pmatrix},\qquad
  \rho(B_2)= 0 \qquad\text{and}\qquad \rho(B_3) = \begin{pmatrix} 0 &
    0 \\ 0 & t \end{pmatrix},
\end{equation}
where $t \neq 0$.
One sees that two pairs of parameters $(t,y),(t',y')\in(\RR\setminus \{0\})\times \RR$ yield weakly equivalent representations if and only if $(t,y)=(t',y')$ or $yy'\neq 0$ and $(t',y')=(-\tfrac{t}{y},\tfrac{1}{y})$, so a complete set of isoclasses is obtained by restricting ourselves to the subset
\begin{equation}
  \big\{
     (t,y)\in(\RR\setminus \{0\})\times [-1,1]\mid y\neq 1 \text{ or } t > 0
  \big\}.
\end{equation}

From the tables, it follows that $\alpha,\sigma,(W^*)^\fb$ are zero.
With $\rho$ given by equation~\eqref{eq:BIII-ker-B2-im-diag}, we see
that $\varphi$ always vanishes except when $t = -2$ and $y=-1$, in
which case $\varphi(\be_1\wedge\be_2)\in\RR B_2$.  As for the weak
automorphism group, we have that if $y=-1$ then
\begin{equation}
	\begin{split}
		\Aut(\rho)=&\left\{
		\left(
		\begin{pmatrix}
		h_{1} & 0 \\
		0 & h_{2} \\
		\end{pmatrix},
		\begin{pmatrix}
		1 & 0 & 0 \\
		0 & b & v \\
		0 & 0 & 1 \\
		\end{pmatrix}
		\right)~\middle|~ bh_1h_2\neq 0
		\right\}\\
		&\bigcup \left\{
		\left(
		\begin{pmatrix}
		0 & h_{2} \\
		h_{1} & 0 \\
		\end{pmatrix},
		\begin{pmatrix}
		-1 & 0 & t \\
		0 & b & v \\
		0 & 0 & 1 \\
		\end{pmatrix}
		\right)~\middle|~ bh_1h_2\neq 0
		\right\},
	\end{split}
\end{equation}
whereas for $y\neq -1$ we get that
\begin{equation}
	\Aut(\rho)=\left\{
	\left(
	\begin{pmatrix}
	h_{1} & 0 \\
	0 & h_{2} \\
	\end{pmatrix},
	\begin{pmatrix}
	1 & 0 & 0 \\
	0 & b & v \\
	0 & 0 & 1 \\
	\end{pmatrix}
	\right)~\middle|~ bh_1h_2\neq 0
	\right\}.
\end{equation}

\subsubsection{\emph{\Romanbar{III}.$4/5$}: $\ker\rho = \RR B_2$ and $\im\rho = \s_{\pi/4}$}
\label{sec:BIII-ker-B2-im-s-pi4}

Here we have
\begin{equation}
  \rho(B_1)=\begin{pmatrix} y & x \\ 0 & y  \end{pmatrix} \qquad\text{and}\qquad
  \rho(B_3)=\begin{pmatrix} r & s \\ 0 & r  \end{pmatrix}.
\end{equation}
We have to distinguish two cases depending on whether or not $y=0$.

If $y\neq 0$, we may use $\Aut(\fb)$ to rescale $y=1$ and to set $r=0$, resulting in
\begin{equation}
  \rho(B_1) = \begin{pmatrix} 1 & x \\ 0 & 1  \end{pmatrix} \qquad\text{and}\qquad
  \rho(B_3)=\begin{pmatrix} 0 & s \\ 0 & 0  \end{pmatrix}.
\end{equation}
The normaliser of $\im \rho$ in $\GL(2,\RR)$ is the group $\T(2,\RR)$
of invertible upper triangular matrices and we may conjugate by
$\T(2,\RR)$ in order to set $s = 1$, resulting in
\begin{equation}
  \rho(B_1) = \begin{pmatrix} 1 & x \\ 0 & 1  \end{pmatrix} \qquad\text{and}\qquad
  \rho(B_3)=\begin{pmatrix} 0 & 1 \\ 0 & 0  \end{pmatrix}.
\end{equation}
Furthermore, different values of $x$ give non-equivalent representations.
From the tables it follows that $\alpha,\sigma$ and $(W^*)^\fb$
are zero, and $\varphi$ is also seen to be zero.
The group of weak automorphisms is
\begin{equation}
	\Aut(\rho)=\left\{
	\left(
	\begin{pmatrix}
	h_1 & h_2 \\
	0 & h_1 \\
	\end{pmatrix},
	\begin{pmatrix}
	1 & 0 & 0 \\
	0 & b & v \\
	0 & 0 & 1 \\
	\end{pmatrix}
	\right)~\middle|~ bh_1\neq 0
	\right\}.
\end{equation}

If $y = 0$, then $x\neq 0$ and we may use $\Aut(\fb)$ to rescale $x=1$ and
set $s= 0$.  Relabelling the parameters for uniformity of notation, we find
\begin{equation}
  \rho(B_1) = \begin{pmatrix} 0 & 1 \\ 0 & 0  \end{pmatrix} \qquad\text{and}\qquad
  \rho(B_3)=\begin{pmatrix} t & 0 \\ 0 & t  \end{pmatrix},
\end{equation}
where $t \neq 0$.
Different values of $t$ correspond to non-isomorphic representations.
From the tables, it follows that
$\alpha,\sigma$ and $(W^*)^\fb$ are zero.
A calculation shows that $\varphi = 0$ unless $t=-1$, in which
case $\varphi(\be_1\wedge \be_2) \in \RR B_2$.  One sees that
\begin{equation}
	\Aut(\rho)=\left\{
	\left(
	\begin{pmatrix}
	a h_1 & h_2 \\
	0 & h_1 \\
	\end{pmatrix},
	\begin{pmatrix}
	a & 0 & 0 \\
	0 & b & v \\
	0 & 0 & 1 \\
	\end{pmatrix}
	\right)~\middle|~ abh_1\neq 0
	\right\}.
\end{equation}

\subsubsection{\emph{\Romanbar{III}.$6/7$}: $\ker\rho = \RR B_2$ and $\im\rho = \CC$}
\label{sec:BIII-ker-B2-im-C}

Here we have
\begin{equation}
  \rho(B_1)= \begin{pmatrix} x & -y \\ y & x \end{pmatrix} \qquad\text{and}\qquad
  \rho(B_3)= \begin{pmatrix} t & -s \\ s & t \end{pmatrix}.
\end{equation}
We must again distinguish between two cases, depending on whether or
not $x = 0$.

If $x=0$, then $y \neq 0$.  We may use $\Aut(\fb)$ to set $y= 1$ and
$s = 0$, resulting in
\begin{equation}
  \rho(B_1)= \begin{pmatrix} 0 & -1 \\ 1 & 0 \end{pmatrix} \qquad\text{and}\qquad
  \rho(B_3)= \begin{pmatrix} t & 0 \\ 0 & t \end{pmatrix},
\end{equation}
with $t \neq 0$.
One sees that each value of $t$ corresponds to a different weak isoclass of representations.
From the tables, it follows that
$\alpha,\sigma$ and $(W^*)^\fb$ are zero.
Similarly to the previous case, we see that
$\varphi = 0$ unless $t = -1$, in which case
$\varphi(\be_1\wedge \be_2) \in \RR B_2$.  The group of weak automorphisms is
\begin{equation}
	\Aut(\rho)=\left\{
	\left(
	\begin{pmatrix}
	h_1 & h_2 \\
	-\varepsilon h_2 & \varepsilon h_1 \\
	\end{pmatrix},
	\begin{pmatrix}
	\varepsilon & 0 & 0 \\
	0 & b & v \\
	0 & 0 & 1 \\
	\end{pmatrix}
	\right)~\middle|~ b(h_1^2+h_2^2)\neq 0,\varepsilon=\pm 1
	\right\}.
\end{equation}

If $x\neq 0$, we may use $\Aut(\fb)$ to rescale $x = 1$ and eliminate the diagonal terms of $\rho(B_3)$, resulting in
\begin{equation}
  \rho(B_1)= \begin{pmatrix} 1 & -y \\ y & 1 \end{pmatrix} \qquad\text{and}\qquad
  \rho(B_3)= \begin{pmatrix} 0 & -t \\ t & 0 \end{pmatrix},  
\end{equation}
where $t\neq 0$ and where we have relabelled parameters for
uniformity.
It turns out that two pairs $(t,y),(t',y')\in
(\RR\setminus\{0\})\times \RR$ are weakly equivalent if and only if
$(t',y')=\pm (t,y)$, so we obtain a complete set of weak isoclasses by
restricting ourselves to $t>0$.  From the tables, $\alpha,\sigma$ and
$(W^*)^\fb$ are zero. Moreover, a direct computation yields $\varphi =
0$.

The group of weak automorphisms is given by
\begin{equation}
	\Aut(\rho)=\left\{
	\left(
	\begin{pmatrix}
	h_1 & h_2 \\
	-h_2 & h_1 \\
	\end{pmatrix},
	\begin{pmatrix}
	1 & 0 & 0 \\
	0 & b & v \\
	0 & 0 & 1 \\
	\end{pmatrix}
	\right)~\middle|~ b(h_1^2+h_2^2)\neq 0
	\right\}.
\end{equation}

\subsubsection{\emph{\Romanbar{III}.$8$}: $\ker\rho = \vecspan\{B_1,B_2\}$ and $\im\rho = \RR \begin{pmatrix} x & 0 \\ 0 & 1\end{pmatrix}$ with $\left|x\right| \leq 1$}
\label{sec:BIII-ker-B1B2-im-x}

In this case we have $\rho(B_1)= \rho(B_2)= 0$ and
\begin{equation}
  \rho(B_3) =
  \begin{pmatrix}
    t x & 0 \\ 0 & t
  \end{pmatrix}\qquad\text{with $t\neq 0$.}
\end{equation}
We see that two pairs $(t,x),(t',x')\in(\RR\setminus\{0\})\times[-1,1]$ give weakly equivalent representations if and only if $(t,x)=(t',x')$ or $x=-1$ and $t'=-t$, so we obtain a complete set of representatives by restricting ourselves to the subset
\begin{equation}
	\left\{
		 (t,x)\in(\RR\setminus \{0\})\times [-1,1]\mid x\neq -1 \text{ or } t>0 
	\right\}.
\end{equation}
From the tables, $(W^*)^\fb = 0$ unless $x=0$, in which case
$(W^*)^\fb = \RR\be^1$; $\alpha = 0$ unless $x=-1$, in which case $\alpha \in
\RR(\be^1 \otimes \be^2 + \be^2 \otimes \be^1)$, or $x=0$, in which
case $\alpha\in \RR \be^1\otimes\be^1$; and $\sigma=0$ unless
\begin{itemize}
\item $x=0$, in which case $\sigma \in
\vecspan\{\be_1 \otimes \be^1 \otimes \be^1, \be_2 \otimes (\be^1
\otimes \be^2 + \be^2 \otimes \be^1)\}$;
\item or  $x=\frac12$, in which case $\sigma \in \RR \be_2 \otimes
  \be^1 \otimes \be^1$.
\end{itemize}
As for $\varphi$, we see that it is zero unless
\begin{itemize}
\item $x=-1$, in which case $\varphi(\be_1\wedge \be_2) \in \RR B_1$;
\item or $t(x+1)=-2$, in which case $\varphi(\be_1\wedge \be_2) \in \RR B_2$.
\end{itemize}

Finally, one sees that if $x=1$, we have
\begin{equation}
	\Aut(\rho)=\GL(W)\times\Aut(\fb),
\end{equation}
whereas for $x< 1$ we see that
\begin{equation}
	\Aut(\rho)=\left\{
	\left(
	\begin{pmatrix}
	h_1 & 0 \\
	0 & h_2 \\
	\end{pmatrix},
	\begin{pmatrix}
	a & 0 & u \\
	0 & b & v \\
	0 & 0 & 1 \\
	\end{pmatrix}
	\right)~\middle|~ abh_1h_2\neq 0
	\right\}.
\end{equation}

\subsubsection{\emph{\Romanbar{III}.$9$}: $\ker\rho = \vecspan\{B_1,B_2\}$ and $\im\rho = \RR \begin{pmatrix} 1 & 1 \\ 0 & 1\end{pmatrix}$}
\label{sec:BIII-ker-B1B2-im-jordan-1}

In this case we take $\rho(B_1)=\rho(B_2)=0$ and
\begin{equation}
  \rho(B_3) =
  \begin{pmatrix}
    t & t \\ 0 & t
  \end{pmatrix} \qquad \text{where $t \neq 0$.}
\end{equation}
It is clear that different values of $t$ yield different weak isoclasses of representations.
The map $\varphi = 0$ unless $t =-1$,
in which case $\varphi(\be_1 \wedge \be_2) \in \RR B_2$.  From the tables, we have
that $(W^*)^\fb$, $\alpha$ and $\sigma$ are all zero.

The group of weak automorphisms is
\begin{equation}
	\Aut(\rho)=\left\{
	\left(
	\begin{pmatrix}
	h_1 & h_2 \\
	0 & h_1 \\
	\end{pmatrix},
	\begin{pmatrix}
	a & 0 & u \\
	0 & b & v \\
	0 & 0 & 1 \\
	\end{pmatrix}
	\right)~\middle|~ abh_1\neq 0
	\right\}.
\end{equation}

\subsubsection{\emph{\Romanbar{III}.$10$}: $\ker\rho = \vecspan\{B_1,B_2\}$ and $\im\rho = \RR \begin{pmatrix} 0 & 1 \\ 0 & 0\end{pmatrix}$}
\label{sec:BIII-ker-B1B2-im-jordan-0}

Here $\rho(B_1)=\rho(B_2)=0$ and we may put $\rho(B_3)$ in Jordan normal form to obtain
\begin{equation}
  \rho(B_3) =
  \begin{pmatrix}
    0 & 1\\ 0 & 0 
  \end{pmatrix}
\end{equation}
From the tables, $\alpha \in \RR \be^2 \otimes \be^2$,
$\sigma \in \vecspan\{\be_1 \otimes \be^2 \otimes \be^2, \be_1 \otimes
(\be^1 \otimes \be^2 + \be^2 \otimes \be^1) + 2 \be_2 \otimes \be^2
\otimes \be^2\}$ and $(W^*)^\fb = \RR\be^2$.
We also see that $\im \varphi \subseteq \RR B_1$.  We have
\begin{equation}
	\Aut(\rho)=\left\{
	\left(
	\begin{pmatrix}
	h_1 & h_2 \\
	0 & h_1 \\
	\end{pmatrix},
	\begin{pmatrix}
	a & 0 & u \\
	0 & b & v \\
	0 & 0 & 1 \\
	\end{pmatrix}
	\right)~\middle|~ abh_1\neq 0
	\right\}.
\end{equation}

\subsubsection{\emph{\Romanbar{III}.$11$}: $\ker\rho = \vecspan\{B_1,B_2\}$ and $\im\rho = \RR \begin{pmatrix} y & -1\\1  &  y\end{pmatrix} $}
\label{sec:BIII-ker-B1B2-im-y}

Here $\rho(B_1)= \rho(B_2) = 0$ and
\begin{equation}
  \rho(B_3) = \begin{pmatrix} ty & -t\\t  &  ty\end{pmatrix} \qquad
  \text{with $t \neq 0$ and $y \geq 0$.}
\end{equation}
A calculation yields that two pairs $(t,y)$ and $(t',y')$ give
weakly equivalent representations if and only if $(t,y)=(t',y')$
or $y=y'=0$ and $t'=-t$, so we may obtain a complete set of weak
isoclasses of representations by restricting our attention to
\begin{equation}
    \big\{(t,y)\in (\RR\setminus\{0\})\times [0,\infty)\mid
    t>0 \text{ or } y>0 \big\}.
\end{equation}
From the tables, $\sigma$ and $(W^*)^\fb$ are zero always, whereas
$\alpha$ is zero unless $y = 0$, in which case $\alpha
\in \RR(\be^1 \otimes \be^1 + \be^2 \otimes \be^2)$.
We see that $\varphi = 0$ unless
\begin{itemize}
\item $y=0$, in which case $\im \varphi \subseteq \RR B_1$;
\item or $ty = -1$, in which  case $\im \varphi \subseteq \RR B_2$.
\end{itemize}
The group of weak automorphisms is
\begin{equation}
	\Aut(\rho)=\left\{
	\left(
	\begin{pmatrix}
	h_1 & h_2 \\
	-h_2 & h_1 \\
	\end{pmatrix},
	\begin{pmatrix}
	a & 0 & u \\
	0 & b & v \\
	0 & 0 & 1 \\
	\end{pmatrix}
	\right)~\middle|~ ab(h_1^2+h_2^2)\neq 0
	\right\}.
\end{equation}

\subsubsection{\emph{\Romanbar{III}.$12$}: $\ker\rho = \vecspan\{B_2,B_3\}$ and $\im\rho = \RR \begin{pmatrix} x & 0 \\ 0 & 1\end{pmatrix}$ with $|x| \leq 1$}
\label{sec:BIII-ker-B2B3-im-x}

Here $\rho(B_2)=\rho(B_3) = 0$ and, since rescaling $B_1$ is an
automorphism of $\fb$, we can just take
\begin{equation}
  \rho(B_1) =
  \begin{pmatrix}
    x & 0 \\ 0 & 1 \end{pmatrix},
\end{equation}
with $|x|\leq 1$.
Different values of $x$ yield different weak isoclasses.
The conditions on $\alpha,\sigma,(W^*)^\fb$
are as in the case treated in Section~\ref{sec:BIII-ker-B1B2-im-x} and shall not be repeated here.
One sees that $\varphi=0$ unless
$x=-1$, in which case $\varphi(\be_1\wedge\be_2)\in \RR B_1$.
As for the weak automorphism group, we distinguish three cases:
\begin{itemize}
	\item If $x=1$ then
	\begin{equation}
		\Aut(\rho)=\left\{
			\left(
				h,
				\begin{pmatrix}
				1 & 0 & 0 \\
				0 & b & v \\
				0 & 0 & 1 \\
				\end{pmatrix}
			\right)~\middle|~ b(\det h)\neq 0
		\right\}.	
	\end{equation}

	\item If $x=-1$, then
	\begin{equation}
		\begin{split}
		\Aut(\rho)=&\left\{
			\left(
			\begin{pmatrix}
			h_1 & 0 \\
			0 & h_2 \\
			\end{pmatrix},
			\begin{pmatrix}
			1 & 0 & 0 \\
			0 & b & v \\
			0 & 0 & 1 \\
			\end{pmatrix}
			\right)~\middle|~ bh_1h_2\neq 0
		\right\}\\
		&\bigcup \left\{
			\left(
			\begin{pmatrix}
			0 & h_2 \\
			h_1 & 0 \\
			\end{pmatrix},
			\begin{pmatrix}
			-1 & 0 & 0 \\
			0 & b & v \\
			0 & 0 & 1 \\
			\end{pmatrix}
			\right)~\middle|~ bh_1h_2\neq 0
		\right\}.
		\end{split}
	\end{equation}

	\item For $x\in (-1,1)$, we have
	\begin{equation}
		\Aut(\rho)=\left\{
			\left(
			\begin{pmatrix}
			h_1 & 0 \\
			0 & h_2 \\
			\end{pmatrix},
			\begin{pmatrix}
			1 & 0 & 0 \\
			0 & b & v \\
			0 & 0 & 1 \\
			\end{pmatrix}
			\right)~\middle|~ bh_1h_2\neq 0
		\right\}.
	\end{equation}
\end{itemize}

\subsubsection{\emph{\Romanbar{III}.$13$}: $\ker\rho = \vecspan\{B_2,B_3\}$ and $\im\rho = \RR \begin{pmatrix} 1 & 1 \\ 0 & 1\end{pmatrix}$}
\label{sec:BIII-ker-B2B3-im-jordan-1}

Here $\rho(B_2)=\rho(B_3) = 0$ and
\begin{equation}
  \rho(B_1) =
  \begin{pmatrix}
    1 & 1 \\ 0 & 1
  \end{pmatrix}.
\end{equation}
The conditions on $\alpha,\sigma,(W^*)^\fb$ are as in the case treated in
Section~\ref{sec:BIII-ker-B1B2-im-jordan-1} and shall not be repeated here.
The map $\varphi$ is necessarily zero.
The group of weak automorphisms is
\begin{equation}
	\Aut(\rho)=\left\{
		\left(
		\begin{pmatrix}
		h_{1} & h_{2} \\
		0 & h_{1} \\
		\end{pmatrix},
		\begin{pmatrix}
		1 & 0 & 0 \\
		0 & b & v \\
		0 & 0 & 1 \\
		\end{pmatrix}
		\right)~\middle|~ bh_1\neq 0
	\right\}.
\end{equation}

\subsubsection{\emph{\Romanbar{III}.$14$}: $\ker\rho = \vecspan\{B_2,B_3\}$ and $\im\rho = \RR \begin{pmatrix} 0 & 1 \\ 0 & 0\end{pmatrix}$}
\label{sec:BIII-ker-B2B3-im-jordan-0}

Here $\rho(B_2)=\rho(B_3) = 0$ and we can rescale $B_1$ to obtain
\begin{equation}
  \rho(B_1) =
  \begin{pmatrix}
    0 & 1 \\ 0 & 0
  \end{pmatrix}.
\end{equation}
The conditions on $\alpha,\sigma,(W^*)^\fb$ are as in the case treated in
Section~\ref{sec:BIII-ker-B1B2-im-jordan-0} and shall not be repeated
here.  Now $\rho(B) (\be_1 \wedge \be_2) = 0$ for all $B \in \fb$, hence
$\im \varphi$ is central.  A calculation shows that
\begin{equation}
	\Aut(\rho)=\left\{
		\left(
		\begin{pmatrix}
		a h_{1} & h_{2} \\
		0 & h_{1} \\
		\end{pmatrix},
		\begin{pmatrix}
		a & 0 & 0 \\
		0 & b & v \\
		0 & 0 & 1 \\
		\end{pmatrix}
		\right)~\middle|~ abh_1\neq 0
	\right\}.
\end{equation}

\subsubsection{\emph{\Romanbar{III}.$15$}: $\ker\rho = \vecspan\{B_2,B_3\}$ and $\im\rho =
  \RR \begin{pmatrix} y & -1\\1  &  y\end{pmatrix}$ with $y\geq 0$}
\label{sec:BIII-ker-B2B3-im-y}

Here $\rho(B_2)=\rho(B_3) = 0$ and, after possibly rescaling $B_1$, we may suppose
\begin{equation}
  \rho(B_1) =
  \begin{pmatrix}
    y & -1 \\ 1  & y
  \end{pmatrix}, \quad y\geq 0.
\end{equation}
One sees that different values of $y$ correspond to different weak isoclasses for $\rho$.
The conditions on $\alpha,\sigma,(W^*)^\fb$ are as in the case
treated in Section~\ref{sec:BIII-ker-B1B2-im-y} and shall not be
repeated here.  We have that $\varphi = 0$ unless $y = 0$, in which
case $\im\varphi$ is central.  If $y=0$, then the weak automorphism
group is
\begin{equation}
	\Aut(\rho)=\left\{
		\left(
		\begin{pmatrix}
		h_1 & h_2 \\
		-\varepsilon h_2 & \varepsilon h_1 \\
		\end{pmatrix},
		\begin{pmatrix}
		\varepsilon & 0 & 0 \\
		0 & b & v \\
		0 & 0 & 1 \\
		\end{pmatrix}
		\right)~\middle|~ b(h_1^2+h_2^2)\neq 0,\varepsilon = \pm 1
	\right\},
\end{equation}
while for $y>0$ we have
\begin{equation}
	\Aut(\rho)=\left\{
		\left(
		\begin{pmatrix}
		h_1 & h_2 \\
		-h_2 & h_1 \\
		\end{pmatrix},
		\begin{pmatrix}
		1 & 0 & 0 \\
		0 & b & v \\
		0 & 0 & 1 \\
		\end{pmatrix}
		\right)~\middle|~ b(h_1^2+h_2^2)\neq 0
	\right\}.
\end{equation}

\subsubsection{\emph{\Romanbar{III}.$16$}: $\rho = 0$}
\label{sec:BIII-rho-0}

In this case, $\alpha, \sigma$ are unconstrained and
$(W^*)^\fb = W^*$; whereas $\varphi$ takes values in the
centre.  Trivially,
\begin{equation}
	\Aut(\rho)=\GL(W)\times\Aut(\fb).
\end{equation}

\subsubsection{Summary}
\label{sec:BIII-summary}

We summarise the results of this section in Table~\ref{tab:BIII-2d-reps-summary}, with the same notation as that explained in Table~\ref{tab:BI-2d-reps-summary}.

\begin{table}[h!]
  \centering
  \setlength{\extrarowheight}{2pt}
  \caption{Two-dimensional real representations of Bianchi $\Romanbar{III}$.}
  \rowcolors{2}{blue!10}{white}
  \resizebox{\linewidth}{!}{
    \begin{tabular}{>{$}l<{$}*{3}{>{$}c<{$}}>{$}c<{$}*{4}{>{$}c<{$}}}
      \toprule
      \multicolumn{1}{c}{Label} & \multicolumn{1}{c}{$\rho(B_1)$}& \multicolumn{1}{c}{$\rho(B_2)$} & \multicolumn{1}{c}{$\rho(B_3)$} & \multicolumn{1}{c}{Remarks} & \multicolumn{1}{c}{$\alpha$} & \multicolumn{1}{c}{$\sigma$} & \multicolumn{1}{c}{$\varphi$} & \multicolumn{1}{c}{$(W^*)^\fb$}\\
      \midrule[\heavyrulewidth]
      \hypertarget{III1}{\Romanbar{III}.1} & \begin{colorpmatrix}1 & 0 \\ 0 & 1 \end{colorpmatrix} & \begin{colorpmatrix}0 & 1 \\ 0 & 0 \end{colorpmatrix} & \begin{colorpmatrix}-1 & 0 \\ 0 & 1 \end{colorpmatrix} & & & & & \\
      \hypertarget{III2}{\Romanbar{III}.2_x} & 0 & \begin{colorpmatrix}0 & 1 \\ 0 & 0 \end{colorpmatrix} & \begin{colorpmatrix}x & 0 \\ 0 & x+2 \end{colorpmatrix} & & x=-2 & x=-2,-4 & x=-1,-2 & \RR\be^2~(x=-2)\\
      \hypertarget{III3}{\Romanbar{III}.3_{t,y}} & \begin{colorpmatrix}1 & 0 \\ 0 & y \end{colorpmatrix} & 0 & \begin{colorpmatrix}0 & 0 \\ 0 & t \end{colorpmatrix} & \makecell{t\neq 0\\\lvert y\rvert \leq 1 \\ y\neq 1 \text{ or } t>0} & & & (t,y)=(-2,-1) & \\
      \hypertarget{III4}{\Romanbar{III}.4_x} & \begin{colorpmatrix}1 & x \\ 0 & 1 \end{colorpmatrix} & 0 & \begin{colorpmatrix}0 & 1 \\ 0 & 0 \end{colorpmatrix} & & & & & \\
      \hypertarget{III5}{\Romanbar{III}.5_t} & \begin{colorpmatrix}0 & 1 \\ 0 & 0 \end{colorpmatrix} & 0 & \begin{colorpmatrix}t & 0 \\ 0 & t \end{colorpmatrix} & t \neq 0 & & & t=-1 & \\
      \hypertarget{III6}{\Romanbar{III}.6_t} & \begin{colorpmatrix}0 & -1 \\ 1 & 0 \end{colorpmatrix} & 0 & \begin{colorpmatrix}t & 0 \\ 0 & t \end{colorpmatrix} & t\neq 0 & & & t=-1 & \\
      \hypertarget{III7}{\Romanbar{III}.7_{t,y}} & \begin{colorpmatrix}1 & -y \\ y & 1 \end{colorpmatrix} & 0 & \begin{colorpmatrix}0 & -t \\ t & 0 \end{colorpmatrix} & t> 0 & & & & \\
      \hypertarget{III8}{\Romanbar{III}.8_{t,x}} &  0 & 0 & \begin{colorpmatrix}t x & 0 \\ 0 & t \end{colorpmatrix} &\makecell{t\neq 0 \\ \lvert x\rvert \leq 1\\x\neq -1\text{ or }t > 0} & x=-1,0 & x=0,\tfrac12 & \left\{\begin{array}{l}x=-1\\t(x+1)=-2\end{array}\right. & \RR\be^1~(x=0)\\
      \hypertarget{III9}{\Romanbar{III}.9_t} &  0 & 0 & \begin{colorpmatrix}t & t \\ 0 & t \end{colorpmatrix} & t \neq 0 & & & t=-1 & \\
      \hypertarget{III10}{\Romanbar{III}.10} &  0 & 0 & \begin{colorpmatrix}0 & 1 \\ 0 & 0 \end{colorpmatrix} &  & \checkmark & \checkmark & \checkmark & \RR\be^2\\
      \hypertarget{III11}{\Romanbar{III}.11_{t,y}} &  0 & 0 & \begin{colorpmatrix}t y & -t \\ t & t y \end{colorpmatrix} & \makecell{t \neq 0 \\ y\geq 0\\t>0\text{ or }y>0} & y=0 & & \left\{\begin{array}{l}y=0\\t y = -1\end{array}\right. & \\
      \hypertarget{III12}{\Romanbar{III}.12_x} & \begin{colorpmatrix} x & 0 \\ 0 & 1 \end{colorpmatrix} & 0 & 0 & |x|\leq 1 & x=-1,0 & x=0,\tfrac12 & x=-1 & \RR\be^1~(x=0)\\
      \hypertarget{III13}{\Romanbar{III}.13} & \begin{colorpmatrix}1 & 1 \\ 0 & 1 \end{colorpmatrix} & 0 & 0 & & & & & \\
      \hypertarget{III14}{\Romanbar{III}.14} & \begin{colorpmatrix}0 & 1 \\ 0 &  0 \end{colorpmatrix} & 0 & 0 & & \checkmark & \checkmark & \checkmark & \RR\be^2\\
      \hypertarget{III15}{\Romanbar{III}.15_y} & \begin{colorpmatrix} y & -1 \\ 1 & y \end{colorpmatrix} & 0 & 0 & y\geq 0 & y=0 & & y=0 & \\
      \hypertarget{III16}{\Romanbar{III}.16} &  0 & 0 & 0 & & \checkmark & \checkmark & \checkmark & W^*\\
      \bottomrule
    \end{tabular}
  }
  \label{tab:BIII-2d-reps-summary}
\end{table}

\subsection{Bianchi IV}
\label{sec:B-IV}

This is the solvable Lie algebra with brackets $[B_1,B_3]=B_1$ and
$[B_2,B_3]=B_1 + B_2$.  The automorphism group is given by
\begin{equation}
  \Aut(\fb) = \left\{
    \begin{pmatrix} a & b & r \\ 0 & a & s \\ 0 & 0 & 1 \end{pmatrix}~ \middle | ~  a\neq 0\right\}.
\end{equation}
The proper ideals are the first derived ideal $[\fb,\fb]=
\vecspan\{B_1,B_2\}$ and the one-dimensional ideal $\RR B_1$.

\subsubsection{\emph{\Romanbar{IV}.$1$}: $\ker\rho = \RR B_1$}
\label{sec:BIV-1d-ker-rho}

The quotient Lie algebra $\fb/\RR B_1$ is isomorphic to the affine
Lie algebra $\aff(1,\RR)$, so that $\im\rho$ is conjugate to $\s_\theta$
for $\theta \neq \frac\pi4$ and $0\leq \theta <\pi$.

We can choose bases such that $\rho(B_1)=0$, whereas
\begin{equation}
  \rho(B_2) = \begin{pmatrix} 0 & u \\ 0 & 0 \end{pmatrix} \quad\text{and}\quad
  \rho(B_3) = \begin{pmatrix} x & v \\ 0 & x+1 \end{pmatrix}.
\end{equation}
Using $\Aut(\fb)$, we can set $u= 1$ and $v = 0$, so that in the end
\begin{equation}
  \rho(B_1) = 0, \quad \rho(B_2) = \begin{pmatrix} 0 & 1 \\ 0 &
    0 \end{pmatrix} \quad\text{and}\quad \rho(B_3) = \begin{pmatrix} x
    & 0 \\ 0 & x +1 \end{pmatrix},
\end{equation}
where $x \in \RR$ determines $\rho$ up to weak equivalence.

From the tables and from calculations we see that
$\alpha,\sigma,\varphi,(W^*)^\fb$ are all zero except at the
following values of $x$:
\begin{itemize}
\item $x=-2$: $\sigma \in \RR(\be_1 \otimes \be^2 \otimes \be^2)$, but
    $\alpha,\,\varphi$ and $(W^*)^{\fb}$ are zero;
\item and $x=-1$: $\alpha \in \RR(\be^2 \otimes \be^2)$, $\sigma \in
  \RR(\be_1 \otimes (\be^1 \otimes \be^2 + \be^2 \otimes \be^1) + 2
  \be_2 \otimes \be^2\otimes \be^2)$, $(W^*)^\fb = \RR\be^2$ and
  $\varphi(\be_1 \wedge \be_2) \in \RR B_1$.
\end{itemize}

The group of weak automorphisms is
\begin{equation}
  \Aut(\rho)=\left\{
    \left(
      \begin{pmatrix}
        a h & s h\\
        0 & h \\
      \end{pmatrix},
      \begin{pmatrix}
        a & b & r \\
        0 & a & s \\
        0 & 0 & 1 \\
      \end{pmatrix}
    \right)~\middle|~ ah\neq 0
  \right\}.
\end{equation}

\subsubsection{\emph{\Romanbar{IV}.$2$}: $\ker\rho = \vecspan\{B_1,B_2\}$ and $\im\rho = \RR \begin{pmatrix}x & 0 \\ 0 & 1\end{pmatrix}$ with $\lvert x\rvert \leq 1$}
\label{sec:BIV-2d-ker-im-x}

We can find bases such that $\rho(B_1)=\rho(B_2)=0$ and
\begin{equation}
  \rho(B_3) = \begin{pmatrix} t x & 0 \\ 0 & t \end{pmatrix},
\end{equation}
where $t \neq 0$.
It turns out that two pairs $(t,x)$ and $(t',x')$ yield weakly equivalent representations if and only if $(t,x)=(t',x')$ or $x=x'=-1$ and $t'=-t$, so we obtain a complete set of weak isoclasses by restricting ourselves to the subset
\begin{equation}
	\big\{
		(t,x)\in (\RR\setminus \{0\} ) \times [-1,1] \mid x\neq -1 \text{ or } t>0
	\big\}.
\end{equation}

The quadratic form $\alpha \in \sS^2 W^*$ is zero unless $x=0$, in which case $\alpha \in \RR(\be^1 \otimes \be^1)$ or
$x=-1$, which corresponds to $\so(1,1)$, in which case $\alpha \in \RR(\be^1 \otimes \be^2 + \be^2 \otimes
\be^1)$.  There are no nonzero invariants in $W^*$ unless $x=0$, in which case
$(W^*)^\fb = \RR\be^1$.

The symmetric bilinear form $\sigma \in W \otimes \sS^2W^*$ is zero
except when $x=0$, in which case $\sigma \in \vecspan\{\be_1 \otimes
\be^1 \otimes \be^1, \be_2 \otimes (\be^1\otimes \be^2 + \be^2 \otimes
\be^1)\}$, or when $x=\frac12$, in which case $\sigma \in \RR(\be_2
\otimes \be^1 \otimes \be^1)$.

The 2-form $\varphi \in \fb \otimes \ext{2} W^*$ is zero, unless
$t(x+1) =-1$, in which case $\varphi(\be_1 \wedge \be_2) \in \RR B_1$.

Finally, if $x=1$ we have
\begin{equation}
	\Aut(\rho)=\GL(W)\times \Aut(\fb),
\end{equation}
whereas for $x\neq 1$ we see that
\begin{equation}
	\Aut(\rho)=\left\{
		\left(
		\begin{pmatrix}
		h_{1} & 0 \\
		0 & h_{2} \\
		\end{pmatrix},
		\begin{pmatrix}
		a & b & r \\
		0 & a & s \\
		0 & 0 & 1 \\
		\end{pmatrix}
		\right)~\middle|~ ah_1h_2\neq 0
	\right\}.
\end{equation}

\subsubsection{\emph{\Romanbar{IV}.$3$}: $\ker\rho = \vecspan\{B_1,B_2\}$ and $\im\rho = \RR\begin{pmatrix}1 & 1 \\ 0 & 1\end{pmatrix}$}
\label{sec:BIV-2d-ker-im-jordan-1}

We can find bases such that 
\begin{equation}
  \rho(B_1) = \begin{pmatrix} 0 & 0 \\ 0 & 0 \end{pmatrix}, \quad
  \rho(B_2) = \begin{pmatrix} 0 & 0 \\ 0 & 0 \end{pmatrix} \quad\text{and}\quad
  \rho(B_3) = \begin{pmatrix} t & t \\ 0 & t \end{pmatrix},
\end{equation}
where $t \neq 0$ determines $\rho$ up to weak isomorphism.

From the tables, it follows that $\alpha$, $\sigma$ and
$(W^*)^\fb$ are all zero.
A calculation shows that $\varphi = 0$ unless
$t=-\frac12$, in which case $\varphi(\be_1 \wedge \be_2) \in \RR B_1$.
We see that
\begin{equation}
	\Aut(\rho)=\left\{
		\left(
		\begin{pmatrix}
		h_{1} & h_{2} \\
		0 & h_{1} \\
		\end{pmatrix},
		\begin{pmatrix}
		a & b & r \\
		0 & a & s \\
		0 & 0 & 1 \\
		\end{pmatrix}
		\right)~\middle|~ ah_1\neq 0
	\right\}.
\end{equation}

\subsubsection{\emph{\Romanbar{IV}.$4$}: $\ker\rho = \vecspan\{B_1,B_2\}$ and $\im\rho = \RR\begin{pmatrix}0 & 1 \\ 0 & 0\end{pmatrix}$}
\label{sec:BIV-2d-ker-im-jordan-0}

We can find bases such that $\rho(B_1)=\rho(B_2)=0$ and by putting $\rho(B_3)$ in Jordan normal form we can assume that
\begin{equation}
  \rho(B_3) = \begin{pmatrix} 0 & 1 \\ 0 & 0 \end{pmatrix},
\end{equation}

From the tables, it follows that $\alpha \in \RR(\be^2 \otimes
\be^2)$, $(W^*)^\fb = \RR\be^2$ and
\begin{equation}
  \sigma \in \vecspan \left\{\be_1 \otimes \left( \be^1 \otimes \be^2  + \be^2 \otimes \be^1  \right) + 2 \be_2 \otimes \be^2 \otimes \be^2, \be_1 \otimes \be^2 \otimes \be^2\right\}.
\end{equation}
Since $\fb$ has zero centre, it follows that $\varphi = 0$.
The group of weak automorphisms is once again
\begin{equation}
	\Aut(\rho)=\left\{
		\left(
		\begin{pmatrix}
		h_{1} & h_{2} \\
		0 & h_{1} \\
		\end{pmatrix},
		\begin{pmatrix}
		a & b & r \\
		0 & a & s \\
		0 & 0 & 1 \\
		\end{pmatrix}
		\right)~\middle|~ ah_1\neq 0
	\right\}.
\end{equation}

\subsubsection{\emph{\Romanbar{IV}.$5$}: $\ker\rho = \vecspan\{B_1,B_2\}$ and $\im\rho
  =\RR \begin{pmatrix}y  & -1 \\  1 & y\end{pmatrix}$ with $y\geq 0$}
\label{sec:BIV-2d-ker-im-y}

We can find bases such that $\rho(B_1)=\rho(B_2)=0$ and
\begin{equation}
  \rho(B_3) = \begin{pmatrix} t y & -t \\ t & t y \end{pmatrix},
\end{equation}
where $t \neq 0$ and $y\geq 0$.
One sees that two different pairs $(t,y)$ and $(t',y')$ give rise to weakly
equivalent representations if and only if $(t,y)=(t',y')$ or $y=y'=0$ and $t'=-t$,
so we get a complete set of weak isoclasses by restricting ourselves to the
subset
\begin{equation}
	\big\{
		(t,y)\in (\RR\setminus \{0\} ) \times [0,\infty) \mid t>0 \text{ or } y>0
	\big\}.
\end{equation}
From the tables, we see that $\alpha =0$ unless $y=0$, in which case $\alpha
\in \RR(\be^1 \otimes \be^1 + \be^2 \otimes \be^2)$.  It follows that
$\sigma$ and $(W^*)^\fb$ are always zero.
The map $\varphi$ vanishes unless $t y = - \frac12$ in which case
$\varphi(\be_1 \wedge \be_2) \in \RR B_1$.

A calculation yields
\begin{equation}
	\Aut(\rho)=\left\{
		\left(
		\begin{pmatrix}
		h_{1} & h_{2} \\
		-h_{2} & h_{1} \\
		\end{pmatrix},
		\begin{pmatrix}
		a & b & r \\
		0 & a & s \\
		0 & 0 & 1 \\
		\end{pmatrix}
		\right) ~\middle|~ a(h_1^2+h_2^2)\neq 0
	\right\}.
\end{equation}

\subsubsection{\emph{\Romanbar{IV}.$6$}: $\rho = 0$}
\label{sec:BIV-rho-zero}

Since the centre of $\fb$ is zero, it follows that $\varphi$ is zero.  The other data $\alpha,\sigma$ are
unconstrained, but we are free to use $\GL(W)$ to bring them to normal
forms.  Also, of course, $(W^*)^\fb = W^*$.
We have obviously that
\begin{equation}
	\Aut(\rho)=\GL(W)\times \Aut(\fb).
\end{equation}

\subsubsection{Summary}
\label{sec:BIV-summary}

We summarise the results of this section in Table~\ref{tab:BIV-2d-reps-summary}, with the same notation as that explained in Table~\ref{tab:BI-2d-reps-summary}.

\begin{table}[h!]
  \centering
  \setlength{\extrarowheight}{3pt}
  \caption{Two-dimensional real representations of Bianchi $\Romanbar{IV}$.}
  \rowcolors{2}{blue!10}{white}
  \resizebox{\linewidth}{!}{
  \begin{tabular}{>{$}l<{$}*{3}{>{$}c<{$}}>{$}c<{$}*{4}{>{$}c<{$}}}
    \toprule
    \multicolumn{1}{c}{Label} & \multicolumn{1}{c}{$\rho(B_1)$} & \multicolumn{1}{c}{$\rho(B_2)$} & \multicolumn{1}{c}{$\rho(B_3)$} & \multicolumn{1}{c}{Remarks} & \multicolumn{1}{c}{$\alpha$} & \multicolumn{1}{c}{$\sigma$} & \multicolumn{1}{c}{$\varphi$} & \multicolumn{1}{c}{$(W^*)^\fb$}\\
    \midrule[\heavyrulewidth]
    \hypertarget{IV1}{\Romanbar{IV}.1_x} & 0 & \begin{colorpmatrix} 0 & 1\\0 & 0\end{colorpmatrix} & \begin{colorpmatrix} x & 0\\ 0 & x+1\end{colorpmatrix} & & x=-1 & x=-1,-2 & x=-1 & \RR\be^2~ (x=-1)\\
    \hypertarget{IV2}{\Romanbar{IV}.2_{t,x}} &     0 & 0 & \begin{colorpmatrix} t x & 0\\0 & t\end{colorpmatrix} & \makecell{\lvert x \rvert\leq1\\ t\neq 0\\x\neq -1 \text{ or } t>0} & x=-1,0 & x=0,\frac12 & t(x+1)=-1 & \RR\be^1~(x=0)\\
    \hypertarget{IV3}{\Romanbar{IV}.3_t} &     0 & 0 & \begin{colorpmatrix} t & t\\0 & t\end{colorpmatrix} & t \neq 0 & & & t=-\frac12 & \\
    \hypertarget{IV4}{\Romanbar{IV}.4} &     0 & 0 & \begin{colorpmatrix} 0 & 1\\0 & 0\end{colorpmatrix} &  & \checkmark & \checkmark & & \RR\be^2\\
    \hypertarget{IV5}{\Romanbar{IV}.5_{t,y}} &     0 & 0 & \begin{colorpmatrix} t y & -t\\ t
      & t y\end{colorpmatrix} & \makecell{t\neq 0\\y\geq 0\\t>0\text{ or }y>0} & y=0 & & t y=-\frac12 & \\
    \hypertarget{IV6}{\Romanbar{IV}.6} &     0 & 0 & 0 & & \checkmark & \checkmark & & W^*\\
    \bottomrule
  \end{tabular}
}
  \label{tab:BIV-2d-reps-summary}
\end{table}

\subsection{Bianchi V}
\label{sec:B-V}

This is the solvable Lie algebra with nonzero brackets $[B_1,B_3]=B_1$
and $[B_2,B_3]=B_2$.  The automorphism group is
\begin{equation}
  \Aut(\fb) = \left\{\begin{pmatrix}A & u \\ 0 & 1\end{pmatrix} ~
    \middle | ~  A \in \GL(2,\RR),~u \in \RR^2\right\} \cong \Aff(2,\RR).
\end{equation}
The proper ideals are $[\fb,\fb] = \vecspan\{B_1,B_2\}$ and any line
therein.  Since the affine group acts transitively on lines, we may
choose $\RR B_1$ without loss of generality.  Hence up to
weak automorphisms, the possible values of $\ker\rho$ can be taken to be
$\RR B_1$, $\vecspan\{B_1,B_2\}$ and $\fb$.

\subsubsection{\emph{\Romanbar{V}.$1$}: $\ker\rho = \RR B_1$}
\label{sec:BV-ker-B1}

In this case $\im \rho \cong \fb/\RR B_1$ is nonabelian and hence
isomorphic to $\s_\theta$ for $\theta \neq \frac\pi4$.  We can
therefore choose bases such that $\rho(B_1)=0$, whereas
\begin{equation}
  \rho(B_2) =  \begin{pmatrix}0 & s \\ 0 & 0 \end{pmatrix} \qquad\text{and}\qquad
  \rho(B_3) =  \begin{pmatrix}x & r \\ 0 & t \end{pmatrix}.
\end{equation}
Demanding that $[\rho(B_2),\rho(B_3)]=\rho(B_2)$, we see that $t = x +
1$.  Using $\Aut(\fb)$ we may rescale $B_2$ (so we can set $s=1$) and shift $B_3$ by any
multiple of $B_2$ (so we can set $r=0$), resulting in
\begin{equation}
  \rho(B_1) = 0, \qquad \rho(B_2) =  \begin{pmatrix}0 & 1 \\ 0 & 0 \end{pmatrix} \qquad\text{and}\qquad
  \rho(B_3) =  \begin{pmatrix}x & 0 \\ 0 & x+1 \end{pmatrix},
\end{equation}
where $x$ can be any real number (and it determines $\rho$ up to weak equivalence).

From the tables, we have that
\begin{itemize}
\item $(W^*)^\fb = 0$ unless $x=-1$, in which case $(W^*)^\fb = \RR\be^2$;
\item $\alpha = 0$ unless $x=-1$, in which case $\alpha \in \RR(\be^2
  \otimes \be^2)$;
\item $\sigma = 0$ unless
  \begin{itemize}[label=$\circ$]
  \item $x=-1$, in which case $\sigma \in \RR(\be_1 \otimes (\be^1 \otimes \be^2 + \be^2 \otimes \be^1) + 2 \be_2\otimes \be^2 \otimes \be^2)$; 
  \item or $x=-2$, in which case $\sigma \in \RR \be_1 \otimes \be^2 \otimes \be^2$.
  \end{itemize}
\end{itemize}
Moreover, we see that $\varphi = 0$ unless $x=-1$, in which case
$\im\varphi \subseteq \vecspan\{B_1,B_2\}$.

In this case, the group of weak automorphisms is
\begin{equation}
	\Aut(\rho)=\left\{
		\left(
		\begin{pmatrix}
		a_{22} h & u_2 h \\
		0 & h \\
		\end{pmatrix},
		\begin{pmatrix}
		a_{11} & a_{12} & u_1 \\
		0 & a_{22} & u_2 \\
		0 & 0 & 1 \\
		\end{pmatrix}
		\right)~\middle|~ a_{11}a_{22}h \neq 0
	\right\}.
\end{equation}

\subsubsection{\emph{\Romanbar{V}.$2$}: $\ker\rho=\vecspan\{B_1,B_2\}$ and $\im\rho =  \RR \begin{pmatrix}x & 0 \\ 0 & 1\end{pmatrix}$ with $\lvert x\rvert\leq 1$}
\label{sec:BV-2d-ker-im-x}

We can find bases such that $\rho(B_1)=\rho(B_2)=0$ and
\begin{equation}
  \rho(B_3) = \begin{pmatrix} t x & 0 \\ 0 & t\end{pmatrix},
\end{equation}
where $t \neq 0$.
As with previous cases, we obtain a complete set of isoclasses by restricting ourselves to the subset
\begin{equation}
	\big\{
		(t,x)\in(\RR\setminus\{0\})\times [-1,1] \mid x\neq -1 \text{ or } t>0
	\big\}.
\end{equation}

The situation for $\alpha,\sigma,(W^*)^\fb$ is as for
any other case with such $\im\rho$ (see, e.g.,
Section~\ref{sec:BIV-2d-ker-im-x}) and shall not be repeated here.
As for $\varphi$, we have that it must be zero
unless $t(x+1) = -1$, in which case $\im \varphi \subseteq
\vecspan\{B_1,B_2\}$.
As for the weak automorphism group, we see that
\begin{itemize}
	\item If $x=1$, this group is
	\begin{equation}
		\Aut(\rho)=\GL(W)\times \Aut(\fb).
	\end{equation}
	\item Otherwise, we have
	\begin{equation}
		\Aut(\rho)=\left\{
			\left(
				\begin{pmatrix}
					h_1 & 0 \\
					0 & h_2
				\end{pmatrix},
				\begin{pmatrix}
					A & \boldsymbol{u} \\
					0 & 1
				\end{pmatrix}
			\right)~\middle|~ h_1h_2(\det A)\neq 0
		\right\}.
	\end{equation}
\end{itemize}

\subsubsection{\emph{\Romanbar{V}.$3$}: $\ker\rho=\vecspan\{B_1,B_2\}$ and $\im\rho =  \RR \begin{pmatrix}1 & 1 \\ 0 & 1\end{pmatrix}$}
\label{sec:BV-2d-ker-im-jordan-1}

We can find bases such that $\rho(B_1)=\rho(B_2)=0$ and
\begin{equation}
  \rho(B_3) = \begin{pmatrix} t  & t \\ 0 & t\end{pmatrix},
\end{equation}
where $t \neq 0$ determines $\rho$ up to weak isomorphism.  The situation for $\alpha,\sigma,(W^*)^\fb$
is as for any other case with such $\im\rho$ (see, e.g.,
Section~\ref{sec:BIV-2d-ker-im-jordan-1}) and shall not be repeated
here.  The situation with $\varphi$ is very similar to the previous
two cases: $\varphi = 0$ unless $t = -\frac12$, in which
case $\im\varphi\subseteq \vecspan\{B_1,B_2\}$.
The group of weak automorphisms is
\begin{equation}
	\Aut(\rho)=\left\{
		\left(
			\begin{pmatrix}
				h_1 & h_2 \\
				0 & h_1
			\end{pmatrix},
			\begin{pmatrix}
				A & \boldsymbol{u} \\
				0 & 1
			\end{pmatrix}
		\right)~\middle|~ h_1(\det A)\neq 0
	\right\}.
\end{equation}

\subsubsection{\emph{\Romanbar{V}.$4$}: $\ker\rho=\vecspan\{B_1,B_2\}$ and $\im\rho =  \RR \begin{pmatrix}0 & 1 \\ 0 & 0\end{pmatrix}$}
\label{sec:BV-2d-ker-im-jordan-0}

We can find bases such that $\rho(B_1)=\rho(B_2)=0$ and putting $\rho(B_3)$ in Jordan normal form we get
\begin{equation}
  \rho(B_3) = \begin{pmatrix} 0  & 1 \\ 0 & 0\end{pmatrix},
\end{equation}
The situation for $\alpha,\sigma,(W^*)^\fb$
is as for any other case with such $\im\rho$ (see, e.g.,
Section~\ref{sec:BIV-2d-ker-im-jordan-0}) and shall not be repeated
here.  Here $\varphi=0$ since $\fb$ has trivial centre.
We also obtain
\begin{equation}
	\Aut(\rho)=\left\{
		\left(
			\begin{pmatrix}
				h_1 & h_2 \\
				0 & h_1
			\end{pmatrix},
			\begin{pmatrix}
				A & \boldsymbol{u} \\
				0 & 1
			\end{pmatrix}
		\right)~\middle|~ h_1(\det A)\neq 0
	\right\}.
\end{equation}

\subsubsection{\emph{\Romanbar{V}.$5$}: $\ker\rho=\vecspan\{B_1,B_2\}$ and $\im\rho =
  \RR \begin{pmatrix}y & -1 \\ 1 & y\end{pmatrix}$ with $y\geq0$}
\label{sec:BV-2d-ker-im-y}

We can find bases such that $\rho(B_1)=\rho(B_2)=0$ and
\begin{equation}
  \rho(B_3) = \begin{pmatrix} t y  & -t \\ t & t y\end{pmatrix}.
\end{equation}
where $t \neq 0$ and $y\geq 0$.
As with previous cases, we get a complete set of isoclasses by
taking $(t,y)$ to be in
\begin{equation}
	\big\{
		(t,y)\in (\RR\setminus \{0\} ) \times [0,\infty) \mid t>0 \text{ or } y>0
	\big\}.
\end{equation}
The situation for $\alpha,\sigma,(W^*)^\fb$
is as for any other case with such $\im\rho$ (see, e.g.,
Section~\ref{sec:BIV-2d-ker-im-y}) and shall not be repeated here.
A nonzero $\varphi$ requires $2 t y = -1$, in which
case $\im \varphi \subseteq \vecspan\{B_1, B_2\}$.
The group of weak automorphisms is
\begin{equation}
	\Aut(\rho)=\left\{
		\left(
			\begin{pmatrix}
				h_1 & h_2 \\
				-h_2 & h_1
			\end{pmatrix},
			\begin{pmatrix}
				A & \boldsymbol{u} \\
				0 & 1
			\end{pmatrix}
		\right)~\middle|~ (\det A)(h_1^2+h_2^2)\neq 0
	\right\}.
\end{equation}

\subsubsection{\emph{\Romanbar{V}.$6$}: $\rho=0$}
\label{sec:BV-rho-zero}

Here $\alpha,\sigma$ are unconstrained and all of $W^*$ is
invariant, but since $\fb$ has trivial centre, $\varphi$ is zero.
The group of weak automorphisms is trivially
\begin{equation}
	\Aut(\rho)=\GL(W)\times \Aut(\fb).
\end{equation}

\subsubsection{Summary}
\label{sec:BV-summary}

We summarise the results of this section in Table~\ref{tab:BV-2d-reps-summary}, with the same notation as that explained in Table~\ref{tab:BI-2d-reps-summary}.

\begin{table}[h!]
  \centering
  \setlength{\extrarowheight}{3pt}
  \caption{Two-dimensional real representations of Bianchi $\Romanbar{V}$.}
  \rowcolors{2}{blue!10}{white}
  \resizebox{\linewidth}{!}{
  \begin{tabular}{>{$}l<{$}*{3}{>{$}c<{$}}>{$}c<{$}*{4}{>{$}c<{$}}}
    \toprule
    \multicolumn{1}{c}{Label} & \multicolumn{1}{c}{$\rho(B_1)$} & \multicolumn{1}{c}{$\rho(B_2)$} & \multicolumn{1}{c}{$\rho(B_3)$} & \multicolumn{1}{c}{Remarks} & \multicolumn{1}{c}{$\alpha$} & \multicolumn{1}{c}{$\sigma$} & \multicolumn{1}{c}{$\varphi$} & \multicolumn{1}{c}{$(W^*)^\fb$}\\
    \midrule[\heavyrulewidth]
    \hypertarget{V1}{\Romanbar{V}.1_x} & 0 & \begin{colorpmatrix} 0 & 1\\0 & 0\end{colorpmatrix} & \begin{colorpmatrix} x & 0\\ 0 & x+1 \end{colorpmatrix} & & x=-1 & x = -1,-2 & x=-1 & \RR\be^2~(x=-1)\\
    \hypertarget{V2}{\Romanbar{V}.2_{t,x}} & 0 & 0 & \begin{colorpmatrix} t x & 0\\0 & t\end{colorpmatrix} & \makecell{\lvert x\rvert\leq1\\ t\neq 0 \\x\neq -1 \text{ or }t>0} & x= -1,0 & x= 0,\frac12 & t(x+1)=-1 & \RR\be^1~(x=0)\\
    \hypertarget{V3}{\Romanbar{V}.3_t} & 0 & 0 & \begin{colorpmatrix} t & t\\0 & t\end{colorpmatrix} & t \neq 0 & & & t=-\frac12 & \\
    \hypertarget{V4}{\Romanbar{V}.4} & 0 & 0 & \begin{colorpmatrix} 0 & 1\\0 & 0\end{colorpmatrix} &  & \checkmark & \checkmark & & \RR\be^2\\
    \hypertarget{V5}{\Romanbar{V}.5_{t,y}} & 0 & 0 & \begin{colorpmatrix} t y & -t\\t  & t y\end{colorpmatrix} & \makecell{t\neq 0 \\ y\geq 0\\ t>0\text{ or }y>0} & y=0 & & 2 t y=-1 & \\
    \hypertarget{V6}{\Romanbar{V}.6} & 0 & 0 & 0 & & \checkmark & \checkmark & & W^* \\
    \bottomrule
  \end{tabular}
}
  \label{tab:BV-2d-reps-summary}
\end{table}

\subsection{Bianchi VI$_{c\neq 1}$}
\label{sec:B-VI}

This is the Lie algebra with brackets $[B_1,B_3]=(c-1)B_1$ and
$[B_2,B_3]=(c+1)B_2$ for $c\geq 0$ and $c\neq 1$, since $c=1$ is
Bianchi~\Romanbar{III}, which was already treated in
Section~\ref{sec:B-III}.

The automorphism group depends on whether $c>0$ or $c=0$.  If $c>0$,
\begin{equation}
  \Aut(\fb) = \left\{\begin{pmatrix} a & 0 & u\\0 & b & v\\ 0 & 0 & 1 \end{pmatrix} ~ \middle | ~  ab \neq 0\right\}.
\end{equation}
If $c=0$ there are additional automorphisms
\begin{equation}
  \Aut(\fb)_- = \left\{\begin{pmatrix} 0 & a & u\\b & 0 & v\\ 0 & 0 & - 1 \end{pmatrix} ~ \middle | ~  ab \neq 0\right\}.
\end{equation}

The proper ideals are $\RR B_1$, $\RR B_2$ and $\vecspan\{B_1,B_2\}$,
hence the possible $\ker \rho$ are $\RR B_1$, $\RR B_2$,
$\vecspan\{B_1,B_2\}$ and $\fb$.
Note that the ideals $\RR B_{1}$ and $\RR B_{2}$ are conjugate if and only if $c=0$.

\subsubsection{\emph{\Romanbar{VI}.$1$}: $\ker\rho = \RR B_1$}
\label{sec:BVI-ker-B1}

If $\ker \rho = \RR B_1$, then $\fb/\ker\rho \cong \aff(1,\RR)$ and
hence we may take $\im\rho = \s_\theta$ for $\theta \neq \frac\pi4$.
This means that we can choose bases where $\rho(B_1)=0$, but
\begin{equation}
	\rho(B_2) = \begin{pmatrix} 0 & s \\ 0 & 0 \end{pmatrix}
	\qquad\text{and}\qquad
	\rho(B_3) = \begin{pmatrix} x & r \\ 0 &
		y \end{pmatrix},
\end{equation}
where $s \neq 0$ and $x,y\in\RR$, $(x,y)\neq (0,0)$.  Using
$\Aut(\fb)$ we can rescale $B_2$ to make $s=1$ and translate $B_3$ so
that $r=0$.  Demanding that $[\rho(B_2),\rho(B_3)] = (c+1)\rho(B_2)$,
we find that $y=x+c+1$. In
summary, we have
\begin{equation}
  \rho(B_1)=0, \qquad \rho(B_2) = \begin{pmatrix} 0 & 1 \\ 0 & 0 \end{pmatrix}
  \qquad\text{and}\qquad
  \rho(B_3) = \begin{pmatrix} x & 0 \\ 0 & x+c+1 \end{pmatrix},
\end{equation}
where $\cot\theta = \frac{x}{x+c+1}$.
The parameter $x\in\RR$ determines $\rho$ up to weak isomorphism.
From the tables, $\alpha = 0$
and $(W^*)^\fb=0$ unless $x+c+1=0$, in which case $\alpha \in \RR \be^2
\otimes \be^2$ and $(W^*)^\fb = \RR \be^2$ and $\sigma=0$ except in the following cases:
\begin{itemize}
\item $x+c+1=0$, in which case $\sigma\in\RR(\be_{1}\otimes(\be^{1}\otimes\be^{2}+\be^{2}\otimes\be^{1})+2\be_{2}\otimes\be^{2}\otimes\be^{2})$;
\item or $x+2c+2=0$, in which case $\sigma \in \RR(\be_1\otimes \be^2 \otimes \be^2)$.
\end{itemize}
Concerning $\varphi$, this map has to vanish except in the
following two possibilities:
\begin{itemize}
\item $x+c+1=0$, in which case $\im\varphi \subseteq \RR B_2$;
\item or $x=-c$, in which case $\im \varphi \subseteq \RR B_1$.
\end{itemize}
As for the group of weak automorphisms, we obtain
\begin{equation}
	\Aut(\rho)=\left\{
		\left(
		\begin{pmatrix}
		b h_1 & h_2 \\
		0 & h_1 \\
		\end{pmatrix},
		\begin{pmatrix}
		a & 0 & u \\
		0 & b & \frac{(c+1) h_2}{h_1} \\
		0 & 0 & 1 \\
		\end{pmatrix}
		\right)~\middle|~ abh_1\neq 0
	\right\}.
\end{equation}

\subsubsection{\emph{\Romanbar{VI}.$2$}: $\ker\rho=\RR B_2$}
\label{sec:BVI-ker-B2}

First of all, note that $\rho$ is not weakly equivalent to a representation in the family $\Romanbar{VI}.1$ if and only if $c>0$, so we will assume that this is the case.
Again, $\im \rho$ can be taken to be $\s_\theta$ for $\theta \neq \frac\pi4$ in $[0,\pi)$.
Using $\Aut(\fb)$ and demanding that we get a representation, we find
\begin{equation}
  \rho(B_1) = \begin{pmatrix} 0 & 1 \\ 0 & 0 \end{pmatrix}, \qquad \rho(B_2)=0
  \qquad\text{and}\qquad
  \rho(B_3) = \begin{pmatrix} x & 0 \\ 0 & x+c-1 \end{pmatrix},
\end{equation}
for $x\in\RR$ (which determines $\rho$ up to weak equivalence) and where now $\cot\theta = \frac{x}{x+c-1}$. The
situation with $\alpha,\sigma$ is just as in the previous case
in terms of $\theta$, but the relation between $x$ and $\theta$ is
different.  From the tables, $\alpha = 0$ and $(W^*)^\fb= 0$
unless $x= 1-c$, in which case $\alpha \in \RR \be^2 \otimes \be^2$
and $(W^*)^\fb = \RR\be^2$ and $\sigma=0$ except in the following cases:
\begin{itemize}
\item $x= 1-c$, in which case $\sigma\in\RR(\be_{1}\otimes(\be^{1}\otimes\be^{2}+\be^{2}\otimes\be^{1})+2\be_{2}\otimes\be^{2}\otimes\be^{2})$;
\item or $x=2(1-c)$, in which case $\sigma \in \RR(\be_1\otimes \be^2 \otimes \be^2)$.
\end{itemize}
The analysis of $\varphi$ is \emph{mutatis mutandis} as in the
previous case.  It is only nonzero in the following two cases:
\begin{itemize}
\item $x=1-c$, in which case $\im \varphi \subseteq \RR B_1$;
\item or $x=-c$, in which case $\im \varphi \subseteq \RR	B_2$.
\end{itemize}
In this case, the group of weak automorphisms is
\begin{equation}
	\Aut(\rho)=\left\{
		\left(
		\begin{pmatrix}
		a h_1 & h_2 \\
		0 & h_1 \\
		\end{pmatrix},
		\begin{pmatrix}
		a & 0 & \frac{(c-1) h_2}{h_1} \\
		0 & b & v \\
		0 & 0 & 1 \\
		\end{pmatrix}
		\right)~\middle|~ abh_1\neq 0
	\right\}.
\end{equation}

\subsubsection{\emph{\Romanbar{VI}.$3$}: $\ker\rho=\vecspan\{B_1,B_2\}$ and $\im\rho =  \RR\begin{pmatrix}x & 0 \\ 0 & 1\end{pmatrix}$ with $\lvert x\rvert\leq 1$}
\label{sec:BVI-ker-B1B2-im-x}

We can choose bases where $\rho(B_1)=\rho(B_2)=0$ and
\begin{equation}
  \rho(B_3) =
  \begin{pmatrix}
    t x & 0 \\ 0 & t
  \end{pmatrix} \qquad \text{where $t\neq 0$ and $\lvert x\rvert \leq 1$.}
\end{equation}
We briefly discuss isoclasses.
For $c=0$, one sees that $(t,x)$ and $(t',x')$ yield weakly isomorphic representations if and only if $x=x'$ and $t'=\pm t$, so a complete set of isoclasses is taken by restricting ourselves to $(t,x)\in \RR^+\times [-1,1]$.
Meanwhile, for $c\neq 0$ we get a complete set of isoclasses by restricting ourselves to
\begin{equation}
	\big\{
		(t,x)\in (\RR\setminus \{0\})\times[-1,1]\mid x\neq -1 \text{ or }t>0
	\big\}.
\end{equation}

The situation with $\alpha,\sigma,(W^*)^\fb$ can be read off
from other cases with the same $\im\rho$ (see, e.g.,
Section~\ref{sec:BIV-2d-ker-im-x}) and shall not be repeated.

The map $\varphi$ always vanishes except in the following two cases:
\begin{itemize}
\item $t(x+1) = 1 - c$, in which case $\im \varphi \subseteq \RR B_1$;
\item or $t(x+1) = - 1 - c$, in which case $\im \varphi \subseteq \RR B_2$.
\end{itemize}

As for the weak automorphisms, we have:
\begin{itemize}
	\item If $x=1$, then
	\begin{equation}
		\Aut(\rho)=\GL(W)\times (\Aut(\fb)\setminus \Aut(\fb)_{-}),
	\end{equation}
	where the second factor is the full automorphism group of $\fb$ for $c>0$.
	\item If $x=-1$ and $c=0$, then
	\begin{equation}
		\begin{split}
		\Aut(\rho)=&\left\{
			\left(
			\begin{pmatrix}
			h_{1} & 0 \\
			0 & h_{2} \\
			\end{pmatrix},
			\begin{pmatrix}
			a & 0 & u \\
			0 & b & v \\
			0 & 0 & 1 \\
			\end{pmatrix}
			\right)~\middle|~ abh_1h_2\neq 0
		\right\}\\
		&\bigcup \left\{
			\left(
			\begin{pmatrix}
			0 & h_2 \\
			h_1 & 0 \\
			\end{pmatrix},
			\begin{pmatrix}
			0 & a & u \\
			b & 0 & v \\
			0 & 0 & -1 \\
			\end{pmatrix}
			\right)~\middle|~ abh_1h_2\neq 0
		\right\}
		\end{split}
	\end{equation}
	\item If $x\neq \pm 1$ or $x=-1$ and $c>0$, then
	\begin{equation}
		\Aut(\rho)=\left\{
			\left(
			\begin{pmatrix}
			h_{1} & 0 \\
			0 & h_{2} \\
			\end{pmatrix},
			\begin{pmatrix}
			a & 0 & u \\
			0 & b & v \\
			0 & 0 & 1 \\
			\end{pmatrix}
			\right)~\middle|~ abh_1h_2\neq 0
		\right\}.
	\end{equation}
\end{itemize}

\subsubsection{\emph{\Romanbar{VI}.$4$}: $\ker\rho=\vecspan\{B_1,B_2\}$ and $\im\rho \in  \RR\begin{pmatrix}1 & 1 \\ 0 & 1\end{pmatrix}$}
\label{sec:BVI-ker-B1B2-im-jordan-1}

We can choose bases where $\rho(B_1)=\rho(B_2)=0$ and
\begin{equation}
  \rho(B_3) =
  \begin{pmatrix}
    t & t \\ 0 & t
  \end{pmatrix}\qquad \text{where $t \neq 0$.}
\end{equation}
One sees that for $c>0$ each $t$ determines the isoclass of $\rho$, while for $c=0$ we have to restrict ourselves to $t>0$ to get a complete set of isoclasses. 
The situation with $\alpha,\sigma,(W^*)^\fb$ can be read off
from other cases with the same $\im\rho$ (see, e.g.,
Section~\ref{sec:BIV-2d-ker-im-jordan-1}) and shall not be repeated.
In addition,
we see that $\varphi = 0$ except in two cases:
\begin{itemize}
\item $t = \frac{1-c}2$, in which case $\im \varphi \subseteq \RR B_1$;
\item or $t = -\frac{1+c}2$, in which case $\im \varphi \subseteq \RR B_2$.
\end{itemize}

The weak automorphism group is
\begin{equation}
	\Aut(\rho)=\left\{
		\left(
		\begin{pmatrix}
		h_1 & h_2 \\
		0 & h_1 \\
		\end{pmatrix},
		\begin{pmatrix}
		a & 0 & u \\
		0 & b & v \\
		0 & 0 & 1 \\
		\end{pmatrix}
		\right)~\middle|~ abh_1\neq 0
	\right\}.
\end{equation}

\subsubsection{\emph{\Romanbar{VI}.$5$}: $\ker\rho=\vecspan\{B_1,B_2\}$ and $\im\rho \in  \RR\begin{pmatrix}0 & 1 \\ 0 & 0\end{pmatrix}$}
\label{sec:BVI-ker-B1B2-im-jordan-0}

We can choose bases where $\rho(B_1)=\rho(B_2)=0$ and we may put $\rho(B_3)$ in Jordan normal form to get
\begin{equation}
  \rho(B_3) =
  \begin{pmatrix}
    0 & 1 \\ 0 & 0
  \end{pmatrix}.
\end{equation}
The situation with $\alpha,\sigma,(W^*)^\fb$ can be read off
from other cases with the same $\im\rho$ (see, e.g.,
Section~\ref{sec:BIV-2d-ker-im-jordan-0}) and shall not be repeated.
The map $\varphi$ is readily seen to be zero.

The weak automorphism group is
\begin{equation}
	\Aut(\rho)=\left\{
		\left(
		\begin{pmatrix}
		h_1 & h_2 \\
		0 & h_1 \\
		\end{pmatrix},
		\begin{pmatrix}
		a & 0 & u \\
		0 & b & v \\
		0 & 0 & 1 \\
		\end{pmatrix}
		\right)~\middle|~ abh_1\neq 0
	\right\},
\end{equation}
and for $c=0$ we get the additional automorphisms
\begin{equation}
	\Aut(\rho)_{-}=\left\{
		\left(
		\begin{pmatrix}
		h_1 & h_2 \\
		0 & -h_1 \\
		\end{pmatrix},
		\begin{pmatrix}
		0 & a & u \\
		b & 0 & v \\
		0 & 0 & -1 \\
		\end{pmatrix}
		\right)~\middle|~ abh_1\neq 0
	\right\}.
\end{equation}

\subsubsection{\emph{\Romanbar{VI}.$6$}: $\ker\rho=\vecspan\{B_1,B_2\}$ and $\im\rho \in  \RR\begin{pmatrix}y & -1 \\ 1 & y\end{pmatrix}$ with $y \geq0$}
\label{sec:BVI-ker-B1B2-im-y}

We can choose bases where $\rho(B_1)=\rho(B_2)=0$ and
\begin{equation}
  \rho(B_3) =
  \begin{pmatrix}
    t y & -t \\ t & t y
  \end{pmatrix}\qquad \text{where $t \neq 0$ and $y\geq 0$.}
\end{equation}
We first comment on weak isoclasses.
If $c>0$, then two pairs $(t,y)$ and $(t',y')$ give weakly
equivalent representations if and only if $(t,y)=(t',y')$ or
$y=y'=0$ and $t'=\pm t$, so we obtain a complete set of isoclasses
by taking $(t,y)$ to be in the set
\begin{equation}
    \big\{
            (t,y)\in (\RR\setminus \{0\} ) \times [0,\infty) \mid t>0 \text{ or } y>0
    \big\}.
\end{equation}
However, for $c=0$ we see that $(t,y)$ and $(t',y')$ give weakly
equivalent representations if and only if $\vert t'\vert=\vert t\vert$
and $y'=y$, so we get a complete set of isoclasses by only choosing
pairs $(t,y)$ with $t>0$.  The situation with
$\alpha,\sigma,(W^*)^\fb$ can be read off from other cases with the
same $\im\rho$ (see, e.g., Section~\ref{sec:BIV-2d-ker-im-y}) and
shall not be repeated.

A computation yields that $\varphi = 0$ except
in two cases:
\begin{itemize}
\item $t y = \frac{1-c}2$, in which case $\im\varphi\subseteq \RR B_1$;
\item or $t y = - \frac{1+c}2$, in which case $\im\varphi\subseteq \RR B_2$.
\end{itemize}

Finally, for $c>0$ or $c=0$ and $y\neq 0$ we have
\begin{equation}
	\Aut(\rho)=\left\{
		\left(
		\begin{pmatrix}
		h_1 & h_2 \\
		-h_2 & h_1 \\
		\end{pmatrix},
		\begin{pmatrix}
		a & 0 & u \\
		0 & b & v \\
		0 & 0 & 1 \\
		\end{pmatrix}
		\right)~\middle|~ ab(h_1^2+h_2^2)\neq 0
	\right\},
\end{equation}
while for $(c,y)=(0,0)$ we obtain
\begin{equation}
	\begin{split}
		\Aut(\rho)=&\left\{
			\left(
			\begin{pmatrix}
			h_1 & h_2 \\
			-h_2 & h_1 \\
			\end{pmatrix},
			\begin{pmatrix}
			a & 0 & u \\
			0 & b & v \\
			0 & 0 & 1 \\
			\end{pmatrix}
			\right)~\middle|~ ab(h_1^2+h_2^2)\neq 0
		\right\} \\
		&\bigcup \left\{
			\left(
			\begin{pmatrix}
			h_1 & h_2 \\
			h_2 & -h_1 \\
			\end{pmatrix},
			\begin{pmatrix}
			0 & a & u \\
			b & 0 & v \\
			0 & 0 & -1 \\
			\end{pmatrix}
			\right)~\middle|~ ab(h_1^2+h_2^2)\neq 0
		\right\}.
	\end{split}
\end{equation}

\subsubsection{\emph{\Romanbar{VI}.$7$}: $\rho=0$}
\label{sec:BVI-rho-zero}

As usual $\alpha,\sigma$ are unconstrained and all of $W^*$ is
invariant, but $\im\varphi$ is central and hence zero, since
$c\neq 1$.
The group of weak automorphisms is trivially
\begin{equation}
	\Aut(\rho)=\GL(W)\times \Aut(\fb).
\end{equation}

\subsubsection{Summary}
\label{sec:BVI-summary}

We summarise the results of this section in Table~\ref{tab:BVI-2d-reps-summary}, with the same notation as that explained in Table~\ref{tab:BI-2d-reps-summary}.

\begin{table}
  \centering
  \setlength{\extrarowheight}{5pt}
  \caption{Two-dimensional real representations of Bianchi $\Romanbar{VI}_{1 \neq c \geq 0}$.}
  \rowcolors{2}{blue!10}{white}
  \resizebox{\linewidth}{!}{%
      \begin{tabular}{>{$}l<{$}*{3}{>{$}c<{$}}>{$}c<{$}*{4}{>{$}c<{$}}}
        \toprule
        \multicolumn{1}{c}{Label} & \multicolumn{1}{c}{$\rho(B_1)$} & \multicolumn{1}{c}{$\rho(B_2)$} & \multicolumn{1}{c}{$\rho(B_3)$} & \multicolumn{1}{c}{Remarks} & \multicolumn{1}{c}{$\alpha$} & \multicolumn{1}{c}{$\sigma$} & \multicolumn{1}{c}{$\varphi$} & \multicolumn{1}{c}{$(W^*)^\fb$}\\
        \midrule[\heavyrulewidth]
        \hypertarget{VI1}{\Romanbar{VI}.1_{c,x}} & 0 & \begin{colorpmatrix} 0 & 1\\0 & 0\end{colorpmatrix} & \begin{colorpmatrix} x & 0\\ 0 & x+c+1 \end{colorpmatrix} & & x=-(1+c) & \left\{\begin{array}{l}x=-(1+c)\\ x=-2(1+c)\end{array}\right. & \left\{\begin{array}{l}x=-(1+c)\\x = -c\end{array}\right. & \RR\be^2~(x=-1-c)\\
        \hypertarget{VI2}{\Romanbar{VI}.2_{c,x}} & \begin{colorpmatrix} 0 & 1\\0 & 0\end{colorpmatrix} & 0 & \begin{colorpmatrix} x & 0\\ 0 & x+c-1 \end{colorpmatrix} & c>0 & x=1-c & \left\{\begin{array}{l}x=1-c\\ x=2(1-c)\end{array}\right. & \left\{\begin{array}{l}x=1-c\\x = -c\end{array}\right. & \RR\be^2~(x=1-c)\\
        \hypertarget{VI3}{\Romanbar{VI}.3_{c,t,x}} & 0 & 0 & \begin{colorpmatrix} t x & 0\\0 & t\end{colorpmatrix} & \begin{tabular}{@{}>{$}c<{$}@{}}c=0, \lvert x \rvert \leq 1, t>0 \\\text{or} \\\left\{\begin{array}{c}c>0, t\neq 0, \lvert x \rvert \leq 1 \\ x\neq -1 \text{ or } t>0\end{array}\right.\end{tabular}& x=-1,0 & x=0,\frac12 & t(x+1)=\pm 1 - c & \RR\be^1~(x=0)\\
        \hypertarget{VI4}{\Romanbar{VI}.4_{c,t}} & 0 & 0 & \begin{colorpmatrix} t & t\\0 & t\end{colorpmatrix} & \begin{tabular}{@{}>{$}c<{$}@{}}t\neq 0 \\ c>0 \text{ or }t>0\end{tabular} & & & t=\frac12 (\pm 1 - c) & \\
        \hypertarget{VI5}{\Romanbar{VI}.5_c} & 0 & 0 & \begin{colorpmatrix} 0 & 1\\0 & 0\end{colorpmatrix} &  & \checkmark & \checkmark & & \RR\be^2\\
        \hypertarget{VI6}{\Romanbar{VI}.6_{c,t,y}} & 0 & 0 & \begin{colorpmatrix} t y & -t\\t  & t y\end{colorpmatrix} & \begin{tabular}{@{}>{$}c<{$}@{}}c=0,y\geq 0,t>0\\\text{or}\\\left\{\begin{array}{c} c>0,y\geq 0,t\neq 0 \\ t>0 \text{ or } y> 0 \end{array}\right. \end{tabular} & y=0 & & t y=\frac12(\pm 1- c) & \\
        \hypertarget{VI7}{\Romanbar{VI}.7_c} & 0 & 0 & 0 & & \checkmark & \checkmark & & W^*\\
        \bottomrule
      \end{tabular}
    }
\label{tab:BVI-2d-reps-summary}
\end{table}

\subsection{Bianchi VII$_{c\geq 0}$}
\label{sec:B-VII}

This is the solvable Lie algebra with nonzero brackets $[B_1,B_3]=c
B_1 - B_2$ and $[B_2,B_3]=B_1 + c B_2$, where $c \geq 0$.  The
automorphism group depends on whether or not $c=0$.  If $c>0$, then
\begin{equation}
  \Aut(\fb) = \left\{\begin{pmatrix} a & b & u \\ -b & a & v \\ 0 & 0
      & 1 \end{pmatrix} ~ \middle | ~ a^2 + b^2 > 0\right\}.
\end{equation}
If $c=0$ then there are additional automorphisms
\begin{equation}
  \Aut(\fb)_- = \left\{\begin{pmatrix} a & b & u \\ b & -a & v \\ 0 & 0
      & -1 \end{pmatrix} ~ \middle | ~ a^2 + b^2 > 0\right\}.
\end{equation}
There is a unique proper ideal $[\fb,\fb] = \vecspan\{B_1,B_2\}$ and
this means that for every two-dimensional representation $\rho\colon
\fb\to\gl(W)$, the space $\Hom_\fb (\ext{2} W,\fb)$ is trivial
because otherwise the image of a nonzero equivariant map $\ext{2}
W\to \fb$ would be a one-dimensional ideal.  This means that $\varphi
= 0$ in every case.

The possible $\ker\rho$ are $[\fb,\fb]$ and $\fb$ itself. Since the
quotient $\fb/[\fb,\fb]$ is a one-dimensional Lie algebra, we have
four cases to consider.

\subsubsection{\emph{\Romanbar{VII}.$1$}: $\ker\rho=\vecspan\{B_1,B_2\}$ and $\im\rho =  \RR \begin{pmatrix}x & 0 \\ 0 & 1\end{pmatrix}$ with $\lvert x\rvert \leq 1$}
\label{sec:BVII-2d-ker-im-x}

We may choose bases such that $\rho(B_1)=\rho(B_2)=0$ and
\begin{equation}
  \rho(B_3) = \begin{pmatrix} t x & 0 \\ 0 & t \end{pmatrix} \qquad\text{with $t \neq 0$.}
\end{equation}
If $c=0$, a complete set of isoclasses is obtained by assuming $t>0$, whereas for $c>0$ we get a complete set of isoclasses by restricting ourselves to the set
\begin{equation}
	\big\{
		 (t,x)\in(\RR\setminus\{0\})\times [-1,1]\mid x\neq -1 \text{ or }t>0
	\big\}.
\end{equation}
We may read off the results for $\alpha,\sigma,(W^*)^\fb$ from
any other case with the same $\im\rho$ (see, e.g.,
Section~\ref{sec:BIV-2d-ker-im-x}), so we shall not repeat them here. The
map $\varphi = 0$ automatically.
If $x=1$, then the weak automorphism group is
\begin{equation}
	\Aut(\rho)=\GL(W)\times(\Aut(\fb)\setminus \Aut(\fb)_{-}).
\end{equation}
For $x=-1$ and $c=0$, we obtain
\begin{equation}
	\begin{split}
		\Aut(\rho)=&\left\{
		\left(
		\begin{pmatrix}
		h_1 & 0 \\
		0 & h_2 \\
		\end{pmatrix},
		\begin{pmatrix}
		a & b & u \\
		-b & a & v \\
		0 & 0 & 1 \\
		\end{pmatrix}
		\right)~\middle|~ h_1h_2(a^2+b^2)\neq 0
		\right\} \\
		&\bigcup \left\{
			\left(
			\begin{pmatrix}
			0 & h_2 \\
			h_1 & 0 \\
			\end{pmatrix},
			\begin{pmatrix}
			a & b & u \\
			b & -a & v \\
			0 & 0 & -1 \\
			\end{pmatrix}
			\right)~\middle|~ h_1h_2(a^2+b^2)\neq 0
		\right\}.
	\end{split}
\end{equation}
Finally, for $x\neq \pm 1$ or $c>0$ and $x\neq 1$, we get
\begin{equation}
	\Aut(\rho)=\left\{
	\left(
	\begin{pmatrix}
	h_1 & 0 \\
	0 & h_2 \\
	\end{pmatrix},
	\begin{pmatrix}
	a & b & u \\
	-b & a & v \\
	0 & 0 & 1 \\
	\end{pmatrix}
	\right)~\middle|~ h_1h_2(a^2+b^2)\neq 0
	\right\}.
\end{equation}

\subsubsection{\emph{\Romanbar{VII}.$2$}: $\ker\rho=\vecspan\{B_1,B_2\}$ and $\im\rho =  \RR \begin{pmatrix}1 & 1 \\ 0 & 1\end{pmatrix}$}
\label{sec:BVII-2d-ker-im-jordan-1}

We may choose bases such that $\rho(B_1)=\rho(B_2)=0$ and
\begin{equation}
  \rho(B_3) = \begin{pmatrix} t & t \\ 0 & t \end{pmatrix} \qquad\text{with $t \neq 0$.}
\end{equation}
If $c>0$, then $t$ determines $\rho$ up to isomorphism, whereas for $c=0$ we have to assume $t>0$ to get a complete set of isoclasses.
This is very similar to the previous case: we may read off
$\alpha,\sigma,(W^*)^\fb$ from, e.g.,
Section~\ref{sec:BIV-2d-ker-im-jordan-1}.  As explained earlier,
$\varphi=0$.  The weak automorphism group is
\begin{equation}
	\Aut(\rho)=\left\{
		\left(
		\begin{pmatrix}
		h_1 & h_2 \\
		0 & h_1 \\
		\end{pmatrix},
		\begin{pmatrix}
		a & b & u \\
		-b & a & v \\
		0 & 0 & 1 \\
		\end{pmatrix}
		\right)~\middle|~ h_1(a^2+b^2)\neq 0
	\right\}.
\end{equation}

\subsubsection{\emph{\Romanbar{VII}.$3$}: $\ker\rho=\vecspan\{B_1,B_2\}$ and $\im\rho =  \RR \begin{pmatrix}0 & 1 \\ 0 & 0\end{pmatrix}$}
\label{sec:BVII-2d-ker-im-jordan-0}

We may choose bases such that $\rho(B_1)=\rho(B_2)=0$ and putting $\rho(B_3)$ in Jordan form we obtain
\begin{equation}
  \rho(B_3) = \begin{pmatrix} 0 & 1 \\ 0 & 0 \end{pmatrix}.
\end{equation}
We can read off
$\alpha,\sigma,(W^*)^\fb$ from, e.g.,
Section~\ref{sec:BIV-2d-ker-im-jordan-0}. As before, $\varphi=0$.
The group of weak automorphisms when $c>0$ is
\begin{equation}
	\Aut(\rho)=\left\{
	\left(
	\begin{pmatrix}
	h_1 & h_2 \\
	0 & h_1 \\
	\end{pmatrix},
	\begin{pmatrix}
	a & b & u \\
	-b & a & v \\
	0 & 0 & 1 \\
	\end{pmatrix}
	\right)~\middle|~ h_1(a^2+b^2)\neq 0
	\right\},
\end{equation}
while for $c=0$ we get extra automorphisms of the form
\begin{equation}
	\Aut(\rho)_{-}=\left\{
	\left(
	\begin{pmatrix}
	h_1 & h_2 \\
	0 & -h_1 \\
	\end{pmatrix},
	\begin{pmatrix}
	a & b & u \\
	b & -a & v \\
	0 & 0 & -1 \\
	\end{pmatrix}
	\right)~\middle|~ h_1(a^2+b^2)\neq 0 
	\right\}.
\end{equation}

\subsubsection{\emph{\Romanbar{VII}.$4$}: $\ker\rho=\vecspan\{B_1,B_2\}$ and $\im\rho =   \RR \begin{pmatrix}y & -1 \\ 1 & y\end{pmatrix}$ with $y\geq 0$}
\label{sec:BVII-2d-ker-im-y}

We may choose bases such that $\rho(B_1)=\rho(B_2)=0$ and
\begin{equation}
  \rho(B_3) = \begin{pmatrix} t y & -t \\ t & t y \end{pmatrix} \qquad\text{with $t \neq 0$.}
\end{equation}
If $c>0$, then a complete set of weak isoclasses of representations is obtained
by taking $(t,y)$ in the set
\begin{equation}
    \big\{
            (t,y)\in (\RR\setminus \{0\} ) \times [0,\infty) \mid t>0 \text{ or } y>0
    \big\}.
\end{equation}
However, for $c=0$ we have to assume $t>0$ to get a complete set of isoclasses. 
We may read off $\alpha,\sigma,(W^*)^\fb$ from, e.g.,
Section~\ref{sec:BIV-2d-ker-im-y}. Again $\varphi =0$ automatically.
As for $\Aut(\rho)$, we have:
\begin{itemize}
	\item For $y=0$ and $c=0$, we see that
	\begin{equation}
		\Aut(\rho)=\left\{
		\left(
		\begin{pmatrix}
		h_1 & h_2 \\
		-\varepsilon h_2 & \varepsilon h_1 \\
		\end{pmatrix},
		\begin{pmatrix}
		a & b & u \\
		-\varepsilon b & \varepsilon a & v \\
		0 & 0 & \varepsilon \\
		\end{pmatrix}
		\right)~\middle|~ (h_1^2+h_2^2)(a^2+b^2)\neq 0,\varepsilon=\pm 1
		\right\}.
	\end{equation}
	\item Otherwise, we obtain
	\begin{equation}
		\Aut(\rho)=\left\{
		\left(
		\begin{pmatrix}
		h_1 & h_2 \\
		-h_2 & h_1 \\
		\end{pmatrix},
		\begin{pmatrix}
		a & b & u \\
		-b & a & v \\
		0 & 0 & 1 \\
		\end{pmatrix}
		\right)~\middle|~ (h_1^2+h_2^2)(a^2+b^2)\neq 0
		\right\}.
	\end{equation}
\end{itemize}

\subsubsection{\emph{\Romanbar{VII}.$5$}: $\rho=0$}
\label{sec:BVII-rho-zero}

Here $\alpha,\sigma$ are unconstrained and all of $W^*$ is
invariant, whereas $\varphi$ is zero since $\fb$ has
trivial centre.
We obviously have
\begin{equation}
  \Aut(\rho)=\GL(W)\times \Aut(\fb).
\end{equation}

\subsubsection{Summary}
\label{sec:BVII-summary}

We summarise the results of this section in Table~\ref{tab:BVII-2d-reps-summary}, with the same notation as that explained in Table~\ref{tab:BI-2d-reps-summary}.

\begin{table}
  \centering
  \setlength{\extrarowheight}{3pt}
  \caption{Two-dimensional real representations of Bianchi $\Romanbar{VII}_{c \geq 0}$.}
  \rowcolors{2}{blue!10}{white}
  \resizebox{\linewidth}{!}{
  \begin{tabular}{>{$}l<{$}*{3}{>{$}c<{$}}>{$}c<{$}*{6}{>{$}c<{$}}}
    \toprule
    \multicolumn{1}{c}{Label} & \multicolumn{1}{c}{$\rho(B_1)$} & \multicolumn{1}{c}{$\rho(B_2)$} & \multicolumn{1}{c}{$\rho(B_3)$} & \multicolumn{1}{c}{Remarks} & \multicolumn{1}{c}{$\alpha$} & \multicolumn{1}{c}{$\sigma$} & \multicolumn{1}{c}{$\varphi$} & \multicolumn{1}{c}{$(W^*)^\fb$}\\
    \midrule[\heavyrulewidth]
    \hypertarget{VII1}{\Romanbar{VII}.1_{c,t,x}} & 0 & 0 & \begin{colorpmatrix} t x & 0\\0 & t\end{colorpmatrix} & \begin{tabular}{@{}>{$}c<{$}@{}}c=0, t>0, \lvert x \rvert \leq 1 \\ \text{or} \\ \left\{\begin{array}{c}c>0, t\neq 0, \lvert x \rvert \leq 1 \\ x\neq -1 \text{ or } t>0\end{array}\right.\end{tabular} & x=-1,0 & x=0,\frac12 & & \RR\be^1~(x=0)\\
    \hypertarget{VII2}{\Romanbar{VII}.2_{c,t}} & 0 & 0 & \begin{colorpmatrix} t & t\\0 & t\end{colorpmatrix} & \begin{tabular}{@{}>{$}c<{$}@{}}t \neq 0 \\ c>0 \text{ or } t>0\end{tabular} & & & & \\
    \hypertarget{VII3}{\Romanbar{VII}.3_{c}} & 0 & 0 & \begin{colorpmatrix} 0 & 1\\0 & 0\end{colorpmatrix} &  & \checkmark & \checkmark & & \RR\be^2\\
    \hypertarget{VII4}{\Romanbar{VII}.4_{c,t,y}} & 0 & 0 & \begin{colorpmatrix} t y & -t\\t  & t y\end{colorpmatrix} & \begin{tabular}{@{}>{$}c<{$}@{}}c=0, t>0, y\geq 0 \\ \text{or} \\ \left\{\begin{array}{c}c>0, t\neq 0, y\geq 0 \\ y>0 \text{ or } t>0\end{array}\right.\end{tabular} & y=0 & & & \\
    \hypertarget{VII5}{\Romanbar{VII}.5_c} & 0 & 0 & 0 & & \checkmark & \checkmark & & W^*\\
    \bottomrule
  \end{tabular}
}
  \label{tab:BVII-2d-reps-summary}
\end{table}

\subsection{Bianchi VIII}
\label{sec:B-VIII}

This is $\fb = \fsl(2,\RR)$, which is simple.  The Lie brackets are
$[B_1,B_2]=-B_3$, $[B_2,B_3]=B_1$ and $[B_1,B_3]=-B_2$.  There are no
proper ideals.  So we have to consider two cases: either $\rho =0$ or
else $\rho$ is an isomorphism.

\subsubsection{\emph{$\Romanbar{VIII}.1$: $\ker \rho = 0$}}\label{sec:VIII1}
If $\rho$ is an isomorphism, then $\rho$ is equivalent to the standard representation.
We know from the tables that $\alpha = 0$,
$(W^*)^\fb =0$ and $\sigma\colon \sS^2 W \to W$ has to be zero, since
$\sS^2 W$ is the adjoint representation which is irreducible and
inequivalent to the 2-dimensional defining representation.
Because $\fsl(2,\RR)$ has no one-dimensional ideals, it follows directly that any equivariant map $\varphi\colon \ext{2} W\to \fb$ is zero.
The group of weak automorphisms is
\begin{equation}
	\Aut(\rho)=\big\{
		 (g,\Ad(g)) \mid g\in \GL(2,\RR)
	\big\}.
\end{equation}

\subsubsection{\emph{\Romanbar{VIII}.$2$}: $\rho = 0$}\label{sec:VIII2}
Here, $\alpha,\sigma$ are unconstrained and all of
$W^*$ is invariant, but because $\fb$ has zero centre, $\varphi$ is
zero.
The group of weak automorphisms is clearly
\begin{equation}
	\Aut(\rho)=\GL(W)\times \Aut(\fb).
\end{equation}

\subsection{Bianchi IX}
\label{sec:B-IX}

This is $\fb = \su(2)$, which is simple and has no nontrivial
two-dimensional real representations.  Hence the only possibility here is
$\rho = 0$.  Hence $\alpha,\sigma$ are unconstrained and all of
$W^*$ is invariant, but since $\fb$ has trivial centre, $\varphi$ is
zero.
Clearly, the group of weak automorphisms is
\begin{equation}
	\Aut(\rho)=\GL(W)\times\Aut(\fb).
\end{equation}

\section{Classification of $(3,2)$-KLAs}
\label{sec:classification-3-2}

In this section we classify $(3,2)$-kinematical algebras with spatial
isotropy of dimension larger than $2$. In particular, we provide the
proofs of Theorems~\ref{thm:generic-d-32-KLAs-classification}
and~\ref{thm:3d-32-KLAs-classification} stated in the Introduction.

As we have seen in Section~\ref{sec:struc-32-klas}, for every
$(3,2)$-KLA $\g=\so(V)\oplus (V\otimes W)\oplus \fb$ (with
$\dim V = d \geq 3$) the subspace $\fb$ is a Lie subalgebra, and
therefore the direct sum $\h=\so(V)\oplus \fb$ is also a Lie
subalgebra. Moreover, the adjoint representation $\ad_{\h}$ decomposes
$\g$ into a direct sum of modules $\g = \h \oplus (V\otimes W)$, where
$V$ is the vector representation of $\so(V)$ and $W$ is a (real)
two-dimensional representation of the Bianchi algebra $\fb$.

Now, suppose that we are given a Bianchi algebra $\fb$ and a
two-dimensional representation $W$ of $\fb$, so that $V\otimes W$
becomes a representation of $\h=\so(V)\oplus \fb$. Then we can
construct a $(3,2)$-KLA structure on $\g=\h\oplus (V\otimes W)$ by
viewing $V\otimes W$ as the $2d$-dimensional abelian algebra and
setting $\g=\h\ltimes (V\otimes W)$. We say in this case that $\g$ is
the \textit{split abelian extension} of $\h$ by $V\otimes W$:
\begin{equation}
  \label{eq:split-abelian-extension}
  \begin{tikzcd}
    0 \ar[r] & V \otimes W \ar[r] & \g \ar[r] & \h \ar[r] \arrow[l, bend left=45] & 0.
  \end{tikzcd}
\end{equation}

\begin{example}\label{xm:galilei}
  As mentioned in the introduction, minimally breaking the spatial
  isotropy of Galilei spacetime gives rise to a homogeneous spacetime
  of a $(3,2)$-KLA.  In this case, it is not difficult to see that the
  $(3,2)$-KLA is a split abelian extension associated to the
  representation $\hyperlink{II6}{\Romanbar{II}.6}$ of
  Bianchi~$\Romanbar{II}$.
\end{example}

We shall take the existence of these $(3,2)$-KLAs as accounted for and
not mention them explicitly for the remainder of this section.  We
concentrate instead on those $(3,2)$-KLAs $\g$ with nontrivial Lie
bracket $\ext{2}(V\otimes W) \to \g$.

\subsection{$(3,2)$-KLAs with $d>3$}
\label{sec:d-neq-3-32KLAs}

Let $\g$ be a $(3,2)$-KLA for which $d> 3$.
Due to Theorem~\ref{th:generic-d-32-structure-theorem} the Lie bracket of $\g$ is determined by a two-dimensional real representation $\rho\colon \fb\to \gl(W)$ of a Bianchi Lie algebra $\fb$ and a $\fb$-equivariant map $\varphi\colon \ext{2} W\to \fb$.
In order to reconstruct the Lie bracket from this data, we also need the $\fb$-invariant symmetric bilinear form $\alpha\colon W\times W\to \RR$ defined by relation~\eqref{eq:generic-jacobi-relation}.

Moreover, let $\rho \colon \fb\to \gl(W)$ be a representation.  By
Proposition~\ref{prop:generic-d-relative-isomorphisms} two equivariant
maps $\varphi,\varphi'\colon \ext{2} W \to \fb$ induce relatively
isomorphic algebras (together with $\fb$ and $\rho$) if and only if
there exists an element $(h,\psi)\in\Aut(\rho)$ such that
$\psi\circ \varphi=h^*\varphi'=(\det h)\varphi'$.  As a consequence,
we obtain the following restatement of
Proposition~\ref{prop:generic-d-relative-isomorphisms}:
\begin{proposition}\label{prop:d-neq-3-orbits}
  Let $\fb$ be a Bianchi algebra and $\rho \colon \fb\to \gl(W)$ a
  two-dimensional representation.  Then, the relative isomorphism
  classes of $(3,2)$-KLAs with $d>3$ associated with $\fb$ and $\rho$
  are in a one-to-one correspondence with the orbits of the natural
  action $\Aut(\rho) \curvearrowright  \Hom_\fb(\ext{2} W,\fb)$ given
  by
  \begin{equation}\label{eq:generic-d-klas-action}
    (h,\psi)\cdot \varphi = \frac{1}{\det h} (\psi \circ \varphi).
  \end{equation}
\end{proposition}

Therefore, we need to calculate the orbits of the action~\eqref{eq:generic-d-klas-action} for each of the representations constructed in Section~\ref{sec:2d-reps-bianchi} to determine the relative isomorphism classes of $(3,2)$-KLAs with $d>3$.
We note that the orbit of $0$ corresponds to the split abelian extension.

We first establish some terminology.
Let $\mathsf{G}$ be a group and $\mathsf{G}\curvearrowright X$ an action of $\mathsf{G}$ on a set $X$.
We say that a subset $S\subseteq X$ is a \emph{complete set of normal forms} for the action if the restriction of the canonical projection $X\to X/\mathsf{G}$ to $S$ is a bijection (in other words, if every element of $X$ is in the same orbit as exactly one element of $S$).

The following definition will be useful in order to compute our orbit
spaces of interest. Let $\mathsf{G}\curvearrowright X$ and
$\mathsf{H}\curvearrowright Y$ be two group actions. We say that the
actions are \emph{orbit equivalent} if there exists a set bijection
$f\colon X\to Y$ such that for each $x\in X$ we have
$f(\mathsf{G}\cdot x)=\mathsf{H}\cdot f(x)$. The map $f$ is said to be
an \emph{orbit equivalence map}.

It is clear that if $f\colon X \to Y$ is an orbit equivalence map and
$S\subseteq X$ is a complete set of normal forms for the action $\mathsf{G}
\curvearrowright X$, then $f(S)$ is a complete set of normal forms for
the action $H \curvearrowright Y$.

Before embarking on the classification, we fix the following notation:
by $\{\be_1,\be_2\}$ we denote a basis of $W$, and by
$\{B_1,B_2,B_3\}$ a basis of $\fb$; for each $i$, $j\in \{1,2\}$ we
write
\begin{equation}
  \varphi(\be_1\wedge \be_2)=\varphi^k B_k, \quad \alpha(\be_i,\be_j)=\alpha_{ij}.
\end{equation}
If $(\rho_{ij})$ is the matrix of $\rho(\varphi(\be_1\wedge\be_2))$
with respect to the basis $\{\be_1,\be_2\}$ of $W$, then the
coefficients $\alpha_{ij}$ can be computed by
\begin{equation}
  \begin{pmatrix}
    \alpha_{11} & \alpha_{12} \\
    \alpha_{12} & \alpha_{22}
  \end{pmatrix} =
  \begin{pmatrix}
    -\rho_{21} & \rho_{11} \\
    \rho_{11} & \rho_{12}
  \end{pmatrix}.
\end{equation}
Because the invariant $\alpha$ is computed directly from $\varphi$, we do not focus on it during the classification process, and we simply give its corresponding value when listing the isomorphism classes.

We now determine the orbit spaces of the action~\eqref{eq:generic-d-klas-action} associated with each two-dimensional representation in Section~\ref{sec:2d-reps-bianchi}.

\subsubsection{Bianchi I}

We see from Table~\ref{tab:BI-2d-reps-summary} that the representations $W$ of $\fb$ having nonzero invariant maps $\varphi \colon\ext{2} W \to \fb$ are $\Romanbar{I}.4_{x=-1}$, $\Romanbar{I}.6$, $\Romanbar{I}.7_{y=0}$ and $\Romanbar{I}.8$.

For the case of $\Romanbar{I}.4_{x=-1}$, the action of $\Aut(\rho)$ on
$\Hom_{\fb}(\ext{2} W,\fb)$ has three orbits labelled by the condition $\varphi(\be_1 \wedge \be_2)\in \{ 0,B_1,B_2 \}$.

For the cases of $\Romanbar{I}.6$ and $\Romanbar{I}.7_{y=0}$, one sees that the orbits of the action $\Aut(\rho) \curvearrowright  \Hom_{\fb}(\ext{2} W,\fb)$ are labelled by $\varphi(\be_1\wedge \be_2)\in \{ 0,\pm B_1, B_2 \}$.

Finally, for the trivial representation $\Romanbar{I}.8$, because the projection of $\Aut(\rho)$ onto $\Aut(\fb)=\GL(\fb)$ is surjective, there are two possible normal forms of $\varphi$, obtained by letting $\varphi = 0$ or $\varphi(\be_1\wedge \be_2)=B_1$.

\subsubsection{Bianchi II}

From Table~\ref{tab:BII-2d-reps-summary} we obtain that the representations $W$ of $\fb$ for which $\Hom_\fb(\ext{2} W,\fb)$ is nontrivial are $\Romanbar{II}.4_{x=-1}$, $\Romanbar{II}.6$, $\Romanbar{II}.7_{y=0}$ and $\Romanbar{II}.8$.
In all cases we have $\Hom_{\fb}(\ext{2} W,\fb)=\ext{2} W^* \otimes \RR B_1$.
Because $\RR B_1$ is the only one-dimensional ideal of $\fb$ and it is contained in the kernel of any two-dimensional representation $\rho$ of $\fb$, we deduce immediately that $\alpha=0$ for any choice of $\varphi$.

It is easy to check that in all cases the map $\varphi \in \Hom_{\fb}(\ext{2} W,\fb)\mapsto \varphi^1 \in \RR$ is an orbit equivalence map between the action of $\Aut(\rho)$ on $ \Hom_{\fb}(\ext{2} W,\fb) $ and the standard action of $(\RR \setminus \{ 0 \})$ on $\RR$.
Consequently, every nonzero $\varphi \in \Hom_{\fb}(\ext{2} W,\fb)$ lies in the same orbit as $(\be^1 \wedge \be^2 )\otimes B_1$.

\subsubsection{Bianchi III}

In this setting, we see from Table~\ref{tab:BIII-2d-reps-summary} that the representations $W$ of $\fb$ for which the space $\Hom_{\fb}(\ext{2} W,\fb)$ is nontrivial are the following ones: $\Romanbar{III}.2_{x=-1}$, $\Romanbar{III}.2_{x=-2}$, $\Romanbar{III}.3_{(t,y)=(-2,-1)}$, $\Romanbar{III}.5_{t=-1}$, $\Romanbar{III}.6_{t=-1}$, $\Romanbar{III}.8_{x=-1}$, $\Romanbar{III}.8_{t(x+1)=-2}$, $\Romanbar{III}.9_{t=-1}$, $\Romanbar{III}.10$, $\Romanbar{III}.11_{y=0}$, $\Romanbar{III}.11_{ty=-1}$, $\Romanbar{III}.12_{x=-1}$, $\Romanbar{III}.14$, $\Romanbar{III}.15_{y=0}$ and $\Romanbar{III}.16$.

In all of the above cases, the space $\Hom_{\fb}(\ext{2} W,\fb)$ is
one-dimensional---in fact, spanned by $B_i$, where $i\in \{ 1,2 \}$
depending on the case. Furthermore, except for the representations
$\Romanbar{III}.2_{x=-2}$, $\Romanbar{III}.14$ and
$\Romanbar{III}.15_{y=0}$, the action of $\Aut(\rho)$ on
$\Hom_{\fb}(\ext{2} W,\fb)$ is orbit equivalent to the standard action
of $\RR\setminus \{ 0 \}$ on $\RR$, meaning that any nonzero $\varphi$
can be renormalised to $(\be^1 \wedge \be^2)\otimes B_i$ by an
adequate element of $\Aut(\rho)$. In the remaining cases, the action
is orbit equivalent to the standard action of $\RR^+$ on $\RR$, so we
obtain a complete set of normal forms for $\varphi$ by letting
$\varphi(\be_1 \wedge \be_2)\in \{ 0,\pm B_i \}$.

\subsubsection{Bianchi IV}
For this Bianchi algebra $\fb$, we see from Table~\ref{tab:BIV-2d-reps-summary} that the representations $W$ admitting nonzero invariant elements in $\Hom_{\fb}(\ext{2} W,\fb)$ are $\Romanbar{IV}.1_{x=-1}$, $\Romanbar{IV}.2_{t(x+1)=-1}$, $\Romanbar{IV}.3_{t=-1/2}$ and $\Romanbar{IV}.5_{ty=-1/2}$.
We have $\Hom_{\fb}(\ext{2} W,\fb)=\ext{2} W^* \otimes \RR B_1$ in all cases and, except for the case of $\Romanbar{IV}.1_{x=-1}$, the action of $\Aut(\rho)$ on this space is orbit equivalent to the action of $\RR \setminus \{ 0 \}$ on $\RR$, so we have a unique nonzero normal form for $\varphi$ given by letting $\varphi(\be_1\wedge \be_2)=B_1$.
In the remaining case, the action is orbit equivalent to the standard action $\RR^+ \curvearrowright \RR$, and thus the three possible normal forms for $\varphi$ are given by $\varphi = 0$ and $\varphi(\be_1 \wedge \be_2) = \pm B_1$.

\subsubsection{Bianchi V}
From Table~\ref{tab:BV-2d-reps-summary} we see that the two-dimensional representations $W$ of $\fb$ that admit nonzero invariant $2$-forms in $ \Hom_{\fb}(\ext{2} W,\fb) $ are $\Romanbar{V}.1_{x=-1}$, $\Romanbar{V}.2_{t(x+1)=-1}$, $\Romanbar{V}.3_{t=-1/2}$ and $\Romanbar{V}.5_{ty=-1/2}$.
In all cases we have that the image of $\varphi$ lies in $\vecspan \{ B_1,B_2 \}$.

For the case of $\Romanbar{V}.1_{x=-1}$, it can be shown that the map
$
    \varphi \in \Hom_{\fb}(\ext{2} W,\fb) \mapsto (\varphi^1,\varphi^2)\in \RR^2
$
establishes an orbit equivalence between the action of $\Aut(\rho)$ on $\Hom_{\fb}(\ext{2} W,\fb)$ and the standard representation of the group
\begin{equation}
    \mathsf{L}= \left\{
        \begin{pmatrix}
            a & b \\
            0 & d
        \end{pmatrix} ~\middle|~ a \neq 0,\, d > 0
    \right\}
\end{equation}
on the plane.
Now, a complete set of normal forms for the action $\mathsf{L} \curvearrowright \RR^2$ is $\{ (0,0),(1,0),(0,\pm 1) \}$, and therefore we have four orbits of the $\Aut(\rho)$-action labelled by the possibilities $\varphi(\be_1 \wedge \be_2)\in \{ 0,B_1, \pm B_2 \}$.

The remaining cases can be tackled at the same time.
Indeed, in all of these the map $\varphi \in \Hom_{\fb}(\ext{2} W,\fb)\mapsto (\varphi^1, \varphi^2)\in \RR^2$ is an orbit equivalence map between the action of $\Aut(\rho)$ on $\Hom_{\fb}(\ext{2} W,\fb)$ and the standard representation of $\GL(2,\RR)$ on $\RR^2$.
Because of this, every nonzero $\varphi$ can be renormalised to the map $(\be^1 \wedge \be^2)\otimes B_1$.

\subsubsection{Bianchi VI}

From Table~\ref{tab:BVI-2d-reps-summary}, we see that the
two-dimensional representations $W$ of $\fb$ that admit nonzero
invariant elements in $\Hom_{\fb}(\ext{2} W,\fb)$ are the following:
$\Romanbar{VI}.1_{x=-1-c}$, $\Romanbar{VI}.1_{x=-c}$,
$\Romanbar{VI}.2_{x=1-c}$, $\Romanbar{VI}.2_{x=-c}$,
$\Romanbar{VI}.3_{t(x+1)=1-c}$, $\Romanbar{VI}.3_{t(x+1)=-1-c}$,
$\Romanbar{VI}.4_{2t=1-c}$, $\Romanbar{VI}.4_{2t=-1-c}$,
$\Romanbar{VI}.6_{2ty=1-c}$ and $\Romanbar{VI}.6_{2ty=-1-c}$.  Note
that in all of these cases we have
$\Hom_{\fb}(\ext{2} W,\fb)=\ext{2} W^* \otimes \RR B_i$ for $i=1$ or
$i=2$ depending on the representation.

In all cases but $\Romanbar{VI}.1_{x=-1-c}$ and
$\Romanbar{VI}.2_{x=1-c}$, the action of $\Aut (\rho)$ on the space $
\Hom_{\fb}(\ext{2} W,\fb)$ is orbit equivalent to the standard action
$(\RR \setminus \{ 0 \})\curvearrowright \RR$, and therefore every
nonzero $\varphi$ can be renormalised to $(\be^1\wedge \be^2)\otimes
B_i$.

For the two remaining representations, the action of $\Aut(\rho)$ on $\Hom_{\fb}(\ext{2} W,\fb)$ is orbit equivalent to the usual action of $\RR^+$ on $\RR$, and thus there are three possible normal forms of $\varphi$ given by letting $\varphi(\be_1 \wedge \be_2)\in \{ 0,\pm B_i \}$.

\subsubsection{Bianchi VII, VIII and IX}
For these cases there is nothing to do, because if $\fb$ is any of these Bianchi algebras and $W$ is a two-dimensional representation of $\fb$, we get $\Hom_{\fb}(\ext{2} W,\fb)=0$ automatically.

\subsubsection{Summary}

We now gather the results of our calculations.
Tables~\ref{tab:generic-d-class-II} and~\ref{tab:generic-d-class-III} contain for each two-dimensional representation $\rho$ the nonzero orbit representatives of the standard action of $\Aut(\rho)$ on $\Hom_\fb(\ext{2} W, \fb)$, which thus give $(3,2)$-KLA structures valid for all $d>3$.
On the one hand, the normal forms of $\varphi$ in Table~\ref{tab:generic-d-class-II} are such that $\rho \circ \varphi = 0$, and thus the corresponding bilinear form $\alpha$ is zero.
These normal forms correspond to relative isomorphism classes of algebras in
class~(\ref{th-item:generic-d-phi}).
On the other hand, the normal forms in Table~\ref{tab:generic-d-class-III} have a nonvanishing corresponding $\alpha$, which we have also described.
Any form in this table corresponds to an algebra in
class~(\ref{th-item:generic-d-phi-alpha}).
Due to Proposition~\ref{prop:generic-d-relative-isomorphisms}, it is clear that
a $(3,2)$-KLA that is in either class~(\ref{th-item:generic-d-phi}) or~(\ref{th-item:generic-d-phi-alpha})
cannot be relatively isomorphic to one in class~(\ref{th-item:generic-d-split-abelian-extension}).
In addition, given two relatively isomorphic $(3,2)$-KLAs $\g$ and $\g'$ with $d>3$, their corresponding
bilinear forms $\alpha$ and $\alpha'$ must have the same rank, so an algebra
in class~(\ref{th-item:generic-d-phi}) cannot be relatively isomorphic to an algebra
in class~(\ref{th-item:generic-d-phi-alpha}).
Because of this, we conclude that Theorem~\ref{thm:generic-d-32-KLAs-classification} holds.

\begin{table}[h!]
    \centering
    \captionsetup{justification=centering}
    \renewcommand{\arraystretch}{1.5}
    \caption{Relative isomorphism classes of generic $(3,2)$-KLAs in class (\ref{th-item:generic-d-phi}).}
    \label{tab:generic-d-class-II}
    \rowcolors{2}{blue!10}{white}
    \begin{tabular}[t]{>{$}l<{$}>{$}c<{$}}
        \toprule
        \multicolumn{1}{l}{\textbf{Representation}} & \multicolumn{1}{c}{$B=\varphi(\be_1 \wedge \be_2)$} \\
        \midrule[\heavyrulewidth] 
        \hyperlink{I4}{\Romanbar{I}.4_{x=-1}} & B_2 \\
        \hyperlink{I6}{\Romanbar{I}.6} & B_2 \\
        \hyperlink{I7}{\Romanbar{I}.7_{y=0}} & B_2 \\
        \hyperlink{I8}{\Romanbar{I}.8} & B_1 \\
        \hyperlink{II4}{\Romanbar{II}.4_{x=-1}} & B_1 \\
        \hyperlink{II6}{\Romanbar{II}.6} & B_1 \\
        \hyperlink{II7}{\Romanbar{II}.7_{y=0}} & B_1 \\
        \hyperlink{II8}{\Romanbar{II}.8} & B_1 \\
        \hyperlink{III2}{\Romanbar{III}.2_{x=-1}} & B_1 \\
        \hyperlink{III3}{\Romanbar{III}.3_{(t,y)=(-2,-1)}} & B_2 \\
        \hyperlink{III5}{\Romanbar{III}.5_{t=-1}} & B_2 \\
        \hyperlink{III6}{\Romanbar{III}.6_{t=-1}} & B_2 \\
        \hyperlink{III8}{\Romanbar{III}.8_{t,x=-1}} & B_1 \\
        \hyperlink{III8}{\Romanbar{III}.8_{t(x+1)=-2}} & B_2 \\
        \hyperlink{III9}{\Romanbar{III}.9_{t=-1}} & B_2 \\
        \hyperlink{III10}{\Romanbar{III}.10} & B_1 \\
        \hyperlink{III11}{\Romanbar{III}.11_{y=0}} & B_1 \\
        \hyperlink{III11}{\Romanbar{III}.11_{ty=-1}} & B_2 \\ 
        \hyperlink{III16}{\Romanbar{III}.16} & B_1 \\
        \bottomrule
    \end{tabular}
    \hspace{2cm}
    \rowcolors{2}{blue!10}{white}
    \begin{tabular}[t]{>{$}l<{$}>{$}c<{$}}
        \toprule
        \multicolumn{1}{l}{\textbf{Representation}} & \multicolumn{1}{c}{$B=\varphi(\be_1 \wedge \be_2)$} \\
        \midrule[\heavyrulewidth] 
        \hyperlink{IV1}{\Romanbar{IV}.1_{x=-1}} & \pm B_1 \\
        \hyperlink{IV2}{\Romanbar{IV}.2_{t(x+1)=-1}} & B_1 \\
        \hyperlink{IV3}{\Romanbar{IV}.3_{t=-1/2}} & B_1 \\
        \hyperlink{IV5}{\Romanbar{IV}.5_{ty=-1/2}} & B_1 \\
        \hyperlink{V1}{\Romanbar{V}.1_{x=-1}} & B_1 \\
        \hyperlink{V2}{\Romanbar{V}.2_{t(x+1)=-1}} & B_1 \\
        \hyperlink{V3}{\Romanbar{V}.3_{t=-1/2}} & B_1 \\
        \hyperlink{V5}{\Romanbar{V}.5_{ty=-1/2}} & B_1 \\
        \hyperlink{VI1}{\Romanbar{VI}.1_{x=-c}} & B_1 \\
        \hyperlink{VI2}{\Romanbar{VI}.2_{x=-c}} & B_2 \\
        \hyperlink{VI3}{\Romanbar{VI}.3_{t(x+1)=1-c}} & B_1 \\
        \hyperlink{VI3}{\Romanbar{VI}.3_{t(x+1)=-1-c}} & B_2 \\
        \hyperlink{VI4}{\Romanbar{VI}.4_{2t=1-c}} & B_1 \\
        \hyperlink{VI4}{\Romanbar{VI}.4_{2t=-1-c}} & B_2 \\
        \hyperlink{VI6}{\Romanbar{VI}.6_{2ty=1-c}} & B_1 \\
        \hyperlink{VI6}{\Romanbar{VI}.6_{2ty=-1-c}} & B_2 \\
        \bottomrule
    \end{tabular}
    \caption*{The Lie bracket in $V\otimes W$ is given by \\ $[u \otimes \be_i,v\otimes \be_j]=\epsilon_{ij}\left< u,v\right>B$.}
\end{table}

\begin{table}[h!]
    \centering
    \captionsetup{justification=centering}
    \caption{Relative isomorphism classes of generic $(3,2)$-KLAs in class (\ref{th-item:generic-d-phi-alpha}).}
    \label{tab:generic-d-class-III}
    \rowcolors{2}{blue!10}{white}
    \begin{tabular}[h!]{>{$}l<{$}>{$}c<{$}>{$}c<{$}}
        \toprule
        \multicolumn{1}{l}{\textbf{Representation}} & \multicolumn{1}{c}{$B=\varphi(\be_1\wedge\be_2)$} & \multicolumn{1}{c}{$(\alpha_{ij})$} \\
        \midrule[\heavyrulewidth]
        \hyperlink{I4}{\Romanbar{I}.4_{x=-1}} & B_1 & \begin{colorpmatrix}0 & -1 \\ -1 & 0\end{colorpmatrix} \\[2.5ex]
        \hyperlink{I6}{\Romanbar{I}.6} & \pm B_1 & \begin{colorpmatrix}0 & 0 \\ 0 & \pm 1\end{colorpmatrix} \\[2.5ex]
        \hyperlink{I7}{\Romanbar{I}.7_{y=0}} & \pm B_1 & \begin{colorpmatrix}\mp 1 & 0 \\ 0 & \mp 1\end{colorpmatrix} \\[2.5ex]
        \hyperlink{III2}{\Romanbar{III}.2_{x=-2}} & \pm B_2 & \begin{colorpmatrix} 0 & 0 \\ 0 & \pm 1 \end{colorpmatrix} \\[2.5ex]
        \hyperlink{III12}{\Romanbar{III}.12_{x=-1}} & B_1 & \begin{colorpmatrix} 0 & -1 \\ -1 & 0 \end{colorpmatrix} \\[2.5ex]
        \hyperlink{III14}{\Romanbar{III}.14} & \pm B_1 & \begin{colorpmatrix} 0 & 0 \\ 0 & \pm1 \end{colorpmatrix} \\[2.5ex]
        \hyperlink{III15}{\Romanbar{III}.15_{y=0}} & \pm B_1 & \begin{colorpmatrix} \mp 1 & 0 \\ 0 & \mp 1 \end{colorpmatrix} \\[2.5ex]
        \hyperlink{V1}{\Romanbar{V}.1_{x=-1}} & \pm B_2 & \begin{colorpmatrix} 0 & 0 \\ 0 & \pm 1 \end{colorpmatrix} \\[2.5ex]
        \hyperlink{VI1}{\Romanbar{VI}.1_{x=-1-c}} & \pm B_2 & \begin{colorpmatrix} 0 & 0 \\ 0 & \pm 1 \end{colorpmatrix} \\[2.5ex]
        \hyperlink{VI2}{\Romanbar{VI}.2_{x=1-c}} & \pm B_1 & \begin{colorpmatrix} 0 & 0 \\ 0 & \pm 1 \end{colorpmatrix} \\
        \bottomrule
    \end{tabular}
    \caption*{The Lie bracket in $V\otimes W$ is given by \\ $[u \otimes \be_i,v\otimes \be_j]=\alpha_{ij}u\curlywedge v+\epsilon_{ij}\left< u,v\right>B$.}
\end{table}

\subsection{$(3,2)$-KLAs with $d=3$}
\label{sec:3d-32-KLAs}

This section is devoted to proving
Theorem~\ref{thm:3d-32-KLAs-classification}, stated in the Introduction.

According to Theorem~\ref{thm:3d-32StructureTheorem-simplified}, a
$(3,2)$-KLA $\g=\so(V)\oplus (V\otimes W) \oplus \fb$ with $\dim V=3$
is specified (after performing an adequate relative isomorphism) by
the $4$-tuple $(\rho, \varphi,\alpha,\sigma)$, where
$\rho \colon \fb \to \gl(W)$ is a real two-dimensional representation
of a Bianchi Lie algebra $\fb$, $\varphi \colon \ext{2} W \to \fb$,
$\alpha \colon W \times W \to \RR$ and
$\sigma \colon W \times W \to W$ are $\fb$-equivariant bilinear maps
subject to the equations~\eqref{eq:3d-algebraic-data-relations}. In
terms of this data, the Lie brackets of $\g$ are such that
$\so(V) \oplus \fb$ is a Lie subalgebra acting on $V \otimes W$ in the
standard way, while the restriction of the bracket to
$\ext{2} (V\otimes W)$ is given by~\eqref{eq:3d-lie-bracket}.
Furthermore, from Proposition~\ref{prop:3d-relative-isomorphisms} it
follows that after fixing $\fb$ and $\rho$, the relative isomorphism
classes of $(3,2)$-KLAs associated with these are in a bijective
correspondence with the orbits of the group
\begin{equation}
  \Aut(\rho) \ltimes (W^*)^\fb = \left\{ (\psi,h,\mu) \in \Aut(\fb)
    \times \GL(W) \times (W^*)^\fb \colon \rho \circ \psi =
    \Ad(h) \circ \rho \right\}
\end{equation}
acting on the set of solutions $(\varphi,\alpha,\sigma)$ of
equations~\eqref{eq:3d-algebraic-data-relations}, where the action of $(\psi,h,\mu)$
on $(\varphi,\alpha,\sigma)$ is $(\varphi',\alpha',\sigma')$ with
\begin{equation}
  h^*\varphi' = \psi \circ \varphi, \quad h^*\alpha' = \alpha - \mu \circ \sigma + \mu \otimes \mu, \quad h^*\sigma' = h \circ \sigma - \mu \otimes h - h \otimes \mu.
\end{equation}

Our first remark is that if $\sigma = 0$ then we have the same data
$(\varphi, \alpha)$ and the same equations as in the generic case
$\dim V > 3$.  Moreover we have the same notion of relative
isomorphism, since if two solutions $(\varphi,\alpha,\sigma=0)$ and
$(\varphi',\alpha',\sigma'=0)$ differ by an element of the form
$(\psi,h,\mu)$, then it follows from the previous equation that
$\mu=0$.  The discussion in Section~\ref{sec:d-neq-3-32KLAs} therefore
applies and we may restrict ourselves to representations $\rho \colon
\fb \to \gl(W)$ with nonzero $\sigma$.

A second observation is that a closer look at
equations~\eqref{eq:3d-algebraic-data-relations} shows that the
Bianchi Lie algebra $\fb$ is directly involved only in determining
$\Hom_{\fb}(\ext{2} W, \fb)$, while all the other conditions depend on
the image $\im\rho \subseteq \gl(W)$.  This results in many solutions
$(\varphi,\alpha,\sigma)$ being identical for different $\rho$.

We choose bases $\{\be_1,\be_2\}$ for $W$ and $\{B_1,B_2,B_3\}$ for $\fb$.
Moreover, we define scalars $\alpha_{ij}$, $\sigma^k_{ij}$ and $\varphi^a$ by
\begin{equation}
    \alpha_{ij} = \alpha(\be_i,\be_j), \qquad \sigma(\be_i,\be_j) = 
  \sigma^k_{ij} \be_k, \qquad \varphi(\be_i \wedge
  \be_j) = \epsilon_{ij} \varphi^a B_a,
\end{equation}
where $\epsilon_{ij} = - \epsilon_{ji}$ is the Levi-Civita symbol normalised so that
$\epsilon_{12}=1$.  We also let $\be^i$, $i \in \{1,2\}$, be the canonical dual basis for $W^*$.

We now summarise our results in terms of the subalgebras $\im \rho$ of
$\gl(2,\RR)$ listed in Table~\ref{tab:gl2-subalgs-inv-sigma} and, of
course, the trivial representation $\rho = 0$.

\subsubsection{Case $\im \rho = \s_0$}
\label{sec:im-rho-s0}

This scenario corresponds to the representations $\Romanbar{III}.2_{x=-2}$,
$\Romanbar{IV}.1_{x=-1}$, $\Romanbar{V}.1_{x=-1}$,
$\Romanbar{VI}.1_{x=-(1+c)}$ and $\Romanbar{VI}.2_{x=1-c}$.  For these
representations, we have
\begin{equation}
	\Hom_\fb(\sS^2 W,W)=\RR (\be_1\otimes(\be^1\otimes\be^2+\be^2\otimes\be^1)+2\be_2\otimes\be^2\otimes\be^2).
\end{equation}

Suppose that $\sigma\in \Hom_\fb(\sS^2 W,W)$ and $\varphi\in \Hom_\fb(\ext{2} W,\fb)$ satisfy~\eqref{eq:3d-algebraic-data-relations}.
Because $(W^*)^\fb=\RR \be^2$ and the covector $\mu_2 \be^2$ acts on $\sigma$ by sending $\sigma_{12}^1$ to $\sigma_{12}^1-\mu_2$, we see that the pair $(\sigma,\varphi)$ admits a renormalisation for which $\sigma=0$, which means that none of these representations admit a $(3,2)$-KLA structure exclusive to $d=3$.

\subsubsection{Case $\im \rho = \s_{\arctan(1/2)}$}
\label{sec:im-rho-s-arctan-half}

This corresponds to $\Romanbar{III}.2_{x=-4}$,
$\Romanbar{IV}.1_{x=-2}$, $\Romanbar{V}.1_{x=-2}$,
$\Romanbar{VI}.1_{x=-2(1+c)}$ and $\Romanbar{VI}.2_{x=2(1-c)}$.  Let
$(\varphi,\sigma,\alpha)$ be a solution of
equations~\eqref{eq:3d-algebraic-data-relations} for which
$\sigma\neq 0$.  Then, the only nonzero component of $\sigma$ is
$\sigma_{22}^1$, and because $(W^*)^\fb=0$ we see that this solution
is not in the same orbit as one with $\sigma =0$.  Moreover, since
$\sigma(\be_2,\be_2)\wedge\be_2=\sigma_{22}^1\be_1\wedge\be_2\neq 0$,
the second equation in~\eqref{eq:3d-algebraic-data-relations} forces
$\varphi=0$.  In particular, we get by direct calculation that
$\alpha=0$.  For all representations $\rho$ in this class, the
projection of the group $\Aut(\rho)$ onto $\GL(W)$ contains matrices
of the form
\begin{equation}
  h =
  \begin{pmatrix}
    a & b \\ 0 & d
  \end{pmatrix}
\end{equation}
with $a/d^2 \neq 0$ but otherwise arbitrary (perhaps after compensating with a
suitable automorphism of $\fb$).  The action of such $h$ on $\be_1
\otimes \be^2 \otimes \be^2$ is to rescale it by $a/d^2$ and hence if
$\sigma^1_{22}$ is nonzero, we may assume that $\sigma^1_{22}=1$.
We conclude that these representations admit a unique $(3,2)$-KLA structure exclusive to $d=3$.

\subsubsection{Case $\im \rho = \RR \begin{pmatrix}0&0\\0&1 \end{pmatrix}$}
\label{sec:im-rho-0001}

This corresponds to $\Romanbar{I}.4_{x=0}$,
$\Romanbar{II}.4_{x=0}$, $\Romanbar{III}.8_{x=0}$,
$\Romanbar{III}.12_{x=0}$, $\Romanbar{IV}.2_{x=0}$,
$\Romanbar{V}.2_{x=0}$, $\Romanbar{VI}.3_{x=0}$ and
$\Romanbar{VII}.1_{x=0}$.
We have in these cases that $(W^*)^\fb=\RR\be^1$ and the possible nonzero components of $\sigma$ are $\sigma_{11}^1$ and $\sigma_{12}^2$, without any restrictions.
Suppose that the triple $(\sigma,\varphi,\alpha)$ satisfies~\eqref{eq:3d-algebraic-data-relations}.
Because the action of $\mu_1\be^1$ sends $\sigma_{11}^1$ to $\sigma_{11}^1-2\mu_1$ and $\sigma_{12}^2$ to $\sigma_{12}^2-\mu_1$, we can renormalise $\sigma$ so that $\sigma_{12}^2=0$. 
If $\sigma\neq 0$, then we see that
\begin{equation}
	\sigma(\be_1,\be_1)\wedge\be_2+\sigma(\be_1,\be_2)\wedge\be_1+\sigma(\be_2,\be_1)\wedge\be_1=\sigma_{11}^1(\be_1\wedge\be_2)
\end{equation}
does not vanish, so the second equation in~\eqref{eq:3d-algebraic-data-relations} forces $\varphi=0$.
By direct computation, one also sees that $\alpha=0$.
Because the projection of $\Aut(\rho)$ onto $\GL(W)$ contains the subgroup of all nonsingular diagonal matrices in all cases, we can rescale $\sigma_{11}^1$ freely, so we conclude that these representations admit a unique $(3,2)$-KLA structure exclusive to $d=3$ given by letting $\sigma_{11}^1=1$ and $\varphi = 0$.

\subsubsection{Case $\im \rho = \RR \begin{pmatrix}\frac12&0\\0&1 \end{pmatrix}$}
\label{sec:im-rho-half001}

This case corresponds to $\Romanbar{I}.4_{x=1/2}$,
$\Romanbar{II}.4_{x=1/2}$, $\Romanbar{III}.8_{x=1/2}$,
$\Romanbar{III}.12_{x=1/2}$, $\Romanbar{IV}.2_{x=1/2}$,
$\Romanbar{V}.2_{x=1/2}$, $\Romanbar{VI}.3_{x=1/2}$ and
$\Romanbar{VII}.1_{x=1/2}$.
In all of these cases we have $(W^*)^\fb=0$ and the only nonzero component of $\sigma$ is $\sigma_{11}^2$.
Once again the projection of $\Aut(\rho)$ onto $\GL(W)$ contains the subgroup of all nonsingular diagonal matrices, so its action allows us to rescale $\sigma_{11}^2$ freely, so if $\sigma\neq 0$ we can renormalise $\sigma$ so that $\sigma_{11}^2=1$.
Observe that
\begin{equation}
	\sigma(\be_1,\be_1)\wedge\be_1=-\sigma_{11}^2\be_1\wedge\be_2,
\end{equation}
so the second equation in~\eqref{eq:3d-algebraic-data-relations} implies that $\varphi= 0$.
As a consequence, $\alpha=0$, and we conclude that all of these representations admit a unique $(3,2)$-KLA structure exclusive to $d=3$ obtained by letting $\sigma_{11}^2=1$ and $\varphi = 0$.

\subsubsection{Case $\im \rho = \RR \begin{pmatrix}0&1\\0&0 \end{pmatrix}$}
\label{sec:im-rho-jordan-0}

This case corresponds to representations $\Romanbar{I}.6$,
$\Romanbar{II}.6$, $\Romanbar{III}.10$, $\Romanbar{III}.14$,
$\Romanbar{IV}.4$, $\Romanbar{V}.4$, $\Romanbar{VI}.5$  and
$\Romanbar{VII}.3$.
The possible nonzero components of $\sigma$ are $\sigma_{12}^1$, $\sigma_{22}^1$ and $\sigma_{22}^2=2\sigma_{12}^1$ (with $\sigma_{12}^1$ and $\sigma_{22}^1$ free).
We also have $(W^*)^\fb=\RR \be^2$ and the action of $\mu_2\be^2$ on $\sigma$ sends
\begin{equation}
	\sigma_{12}^1\mapsto\sigma_{12}^1-\mu_2 \quad \text{and}\quad\sigma_{22}^1\mapsto \sigma_{22}^1.
\end{equation}
Because of this, we can renormalise any solution $(\varphi,\alpha,\sigma)$ so that $\sigma_{12}^1=0$.
For such a solution, assuming that $\sigma\neq 0$, we have $\sigma(\be_2,\be_2)\wedge\be_2=\sigma_{22}^1\be_1\wedge\be_2$, so the second equation in~\eqref{eq:3d-algebraic-data-relations} yields $\varphi=0$.
A calculation shows that $\alpha=0$, and in all cases the action of $\Aut(\rho)$ on $\RR \be_1\otimes\be^2\otimes\be^2$ is orbit equivalent to the standard action of $\RR\setminus \{0\}$ on $\RR$, so we can renormalise $\sigma$ so that $\sigma_{22}^1=1$.
Thus, these representations admit exactly one $(3,2)$-KLA structure exclusive to $d=3$ given by $\sigma_{22}^1=1$ and $\varphi = 0$.

\subsubsection{Case $\rho = 0$}
\label{sec:im-rho-0}

This corresponds to representations $\Romanbar{I}.8$,
$\Romanbar{II}.8$, $\Romanbar{III}.16$,$\Romanbar{IV}.6$,
$\Romanbar{V}.6$, $\Romanbar{VI}.7$, $\Romanbar{VII}.5$,
$\Romanbar{VIII}.2$ and the only real two-dimensional representation
of Bianchi~$\Romanbar{IX}$.

We show that any solution $(\varphi,\alpha,\sigma)$ of~\eqref{eq:3d-algebraic-data-relations} admits a renormalisation with either $\sigma=0$ or $\varphi=0$.
Indeed, from the results in Appendix~\ref{app:binary-cubics} we see that $\sigma$ may be assumed to have trace zero, or in other words, that for every $w_1\in W$ the operator $\hat{\sigma}_{w_1}\in\End(W)$ given by $\hat{\sigma}_{w_1}(w_2)=\sigma(w_1,w_2)$ is in $\fsl(W)$.
In particular, this means that for every $w_2,w_3\in W$ we have
\begin{equation}
  0=\hat{\sigma}_{w_1}(w_2\wedge w_3)=\sigma(w_1,w_2)\wedge w_3+w_2\wedge \sigma(w_1,w_3).
\end{equation}
As a consequence, the expression $\sigma(w_1,w_2)\wedge w_3$ is totally symmetric in the three variables, so the second equation in~\eqref{eq:3d-algebraic-data-relations} becomes
\begin{equation}
  3\varphi(\sigma(w_1,w_2)\wedge w_3)=0, \quad w_1,w_2,w_3\in W.
\end{equation}
Thus, because $\varphi$ is either zero or injective, we deduce
from the equation above that either $\sigma=0$ (so the corresponding
$(3,2)$-KLA structure is generic) or $\varphi=0$.  Note that if
$\varphi=0$, then the coefficients of $\alpha$ are
\begin{equation}\label{eq:alpha-in-terms-of-sigma}
  \begin{split}
    \alpha_{11} &= (\sigma^2_{22} - \sigma^1_{12})\sigma^2_{11} +  (\sigma^1_{11}-\sigma^2_{12})\sigma^2_{12},\\
    \alpha_{12} &= \sigma^1_{12}\sigma^2_{12} - \sigma^1_{22} \sigma^2_{11},\\
    \alpha_{22} &=  (\sigma^1_{11}-\sigma^2_{12})\sigma^1_{22}  + (\sigma^2_{22} - \sigma^1_{12})\sigma^1_{12}.
  \end{split}
\end{equation}
We conclude that in order to determine the $(3,2)$-KLA structures
exclusive to $d=3$, it is necessary and sufficient to understand the
orbit space of the action of $\GL(W)\ltimes W^*$ on
$\Hom(\sS^2 W,W)$ studied in Appendix~\ref{app:binary-cubics}, where
we relate this classification to the classical problem of determining
the normal forms for binary cubic forms of weight $-1$ under the
action of $\GL(2,\RR)$.  Summarising the results in that section, we
have the following normal forms $\sigma_{s}$ (where $s=1,\dots,4$) for
nonzero $\sigma$, where we list only the nonzero coefficients
$\sigma_{ij}^k$ and $\alpha_{ij}$.
\begin{itemize}
\item $(\sigma_1)_{11}^1 = - (\sigma_1)_{12}^2 = \tfrac13$ and $(\sigma_1)_{22}^1 = -1$, hence $(\alpha_1)_{11}=-\tfrac29$, $(\alpha_1)_{22}=-\tfrac23$;
\item $(\sigma_2)_{11}^1 = - (\sigma_2)_{12}^2 = \tfrac13$ and $(\sigma_2)_{22}^1  = 1$, hence $(\alpha_2)_{11}=-\tfrac29$, $(\alpha_2)_{22}=\tfrac23$;
\item $(\sigma_3)_{12}^1 = -(\sigma_3)_{22}^2 = \tfrac13$, hence $(\alpha_3)_{22}=-\tfrac29$;
\item and $(\sigma_4)_{22}^1 = 1$, hence $\alpha_4 = 0$.
\end{itemize}
Recall that it is possible that a general nonzero $\sigma$ is in the
same orbit as the zero tensor.

It will be simpler for us to consider different renormalisations of
these $\sigma_s$ so that their corresponding bilinear forms $\alpha_s$
vanish.  Let us define modified bilinear maps $\tau_s$ as follows:
\begin{itemize}
  \item $\tau_1=(-\operatorname{id}_W,\tfrac23 \be^1)\cdot \sigma_1$, with nonzero components $\tau_{11}^1=\tau_{12}^2=\tau_{22}^1=1$;
  \item $\tau_2=(-\operatorname{id}_W,\tfrac23 \be^1)\cdot \sigma_2$, with nonzero components $\tau_{11}^1=\tau_{12}^2=-\tau_{22}^1=1$;
  \item $\tau_3=(h,\tfrac13\be^2) \cdot \sigma_3=\be_1 \otimes \be^1 \otimes \be^1$, where
  \begin{equation}
    h=\begin{pmatrix} 0 & -1 \\ -1 & 0 \end{pmatrix};
  \end{equation}
  \item and $\tau_4= (-h,0)\cdot \sigma_4= \be_2\otimes \be^1
    \otimes \be^1$, where $h$ is the same matrix as above.
\end{itemize}
Then by construction each $\tau_s$ is in the same orbit as $\sigma_s$,
and a calculation shows that the bilinear form $\beta_s$ constructed
from $\tau_s$ via~\eqref{eq:3d-main-alpha-formula} vanishes
identically.  We therefore obtain that every solution
of~\eqref{eq:3d-algebraic-data-relations} taking the form
$(\varphi=0,\alpha,\sigma)$ is in the same orbit as exactly one of
either $(0,0,0)$ or a $(0,0,\tau_s)$ with $s\in\{1,\dots,4\}$.

\subsubsection{Summary}

Let us now gather our results and prove
Theorem~\ref{thm:3d-32-KLAs-classification}. From our calculations, we
have seen that every $(3,2)$-KLA
$\g=\so(V)\oplus (V\otimes W)\oplus \fb$ with $\dim V = d = 3$ is
relatively isomorphic to one that has an associated two-dimensional
representation $\rho \colon \fb\to \gl(W)$ and such that
$[\fb, V \otimes W]\subseteq V \otimes W$. Moreover, the remaining
objects $\varphi$ and $\sigma$ determining the Lie bracket of $\g$ can
be brought to a normal form where either $\varphi=0$ or $\sigma = 0$.

On the one hand, if $\sigma = 0$, then we are back at the problem of
determining the normal forms of $\varphi$ up to the action of
$\Aut(\rho)$, so we conclude that $\g$ is either a split abelian
extension, an algebra in class~(\ref{th-item:3d-class-II})
(classified in Table~\ref{tab:generic-d-class-II})
or an algebra in class~(\ref{th-item:3d-class-III})
(classified in Table~\ref{tab:generic-d-class-III}).

On the other hand, if $\varphi = 0$, then either $\sigma$ can be
brought to zero by an element of $(W^*)^{\fb}$ (so $\g$ is actually a
split abelian extension in disguise) or it may be transformed into
exactly one of the normal forms described in this section. This last
case corresponds to $\g$ being in
class~(\ref{th-item:3d-class-IV}). The canonical forms for
$\sigma$ (and thus all $(3,2)$-KLAs with $d=3$ and in
class~(\ref{th-item:3d-class-IV})) are collected in
Table~\ref{tab:3d-32-KLAs}, and they all have its corresponding
$\alpha$ equal to zero.

All in all, Theorem~\ref{thm:3d-32-KLAs-classification} is proved.

\begin{table}
    \centering
    \captionsetup{justification=centering}
    \caption{Relative isomorphism classes of $(3,2)$-KLAs with $d=3$ in class (\ref{th-item:3d-class-IV}).}\label{tab:3d-32-KLAs}
    \rowcolors{2}{blue!10}{white}
    \begin{tabular}[t]{>{$}l<{$}*{3}{>{$}r<{$}}}
        \toprule
        \multicolumn{1}{l}{\textbf{Representation}} & \multicolumn{1}{c}{$\sigma_{11}$} & \multicolumn{1}{c}{$\sigma_{12}$} & \multicolumn{1}{c}{$\sigma_{22}$} \\
        \midrule[\heavyrulewidth]
        \hyperlink{I4}{\Romanbar{I}.4_{x=0}} & \be_1 & 0 & 0 \\
        \hyperlink{I4}{\Romanbar{I}.4_{x=1/2}} & \be_2 & 0 & 0 \\
        \hyperlink{I6}{\Romanbar{I}.6} & 0 & 0 & \be_1 \\
        \hyperlink{I8}{\Romanbar{I}.8} & \be_1 & \be_2 & \be_1 \\
        \hyperlink{I8}{\Romanbar{I}.8} & \be_1 & \be_2 & -\be_1 \\
        \hyperlink{I8}{\Romanbar{I}.8} & \be_1 & 0 & 0 \\
        \hyperlink{I8}{\Romanbar{I}.8} & \be_2 & 0 & 0 \\
        \hyperlink{II4}{\Romanbar{II}.4_{x=0}} & \be_1 & 0 & 0 \\
        \hyperlink{II4}{\Romanbar{II}.4_{x=1/2}} & \be_2 & 0 & 0 \\
        \hyperlink{II6}{\Romanbar{II}.6} & 0 & 0 & \be_1 \\
        \hyperlink{II8}{\Romanbar{II}.8} & \be_1 & \be_2 & \be_1 \\
        \hyperlink{II8}{\Romanbar{II}.8} & \be_1 & \be_2 & -\be_1 \\
        \hyperlink{II8}{\Romanbar{II}.8} & \be_1 & 0 & 0 \\
        \hyperlink{II8}{\Romanbar{II}.8} & \be_2 & 0 & 0 \\
        \hyperlink{III2}{\Romanbar{III}.2_{x=-4}} & 0 & 0 & \be_1 \\
        \hyperlink{III8}{\Romanbar{III}.8_{t,x=0}} & \be_1 & 0 & 0 \\
        \hyperlink{III8}{\Romanbar{III}.8_{t,x=1/2}} & \be_2 & 0 & 0 \\
        \hyperlink{III10}{\Romanbar{III}.10} & 0 & 0 & \be_1 \\
        \hyperlink{III12}{\Romanbar{III}.12_{x=0}} & \be_1 & 0 & 0 \\
        \hyperlink{III12}{\Romanbar{III}.12_{x=1/2}} & \be_2 & 0 & 0 \\
        \hyperlink{III14}{\Romanbar{III}.14} & 0 & 0 & \be_1 \\
        \hyperlink{III16}{\Romanbar{III}.16} & \be_1 & \be_2 & \be_1 \\
        \hyperlink{III16}{\Romanbar{III}.16} & \be_1 & \be_2 & -\be_1 \\
        \hyperlink{III16}{\Romanbar{III}.16} & \be_1 & 0 & 0 \\
        \hyperlink{III16}{\Romanbar{III}.16} & \be_2 & 0 & 0 \\
        \hyperlink{IV1}{\Romanbar{IV}.1_{x=-2}} & 0 & 0 & \be_1 \\
        \hyperlink{IV2}{\Romanbar{IV}.2_{t,x=0}} & \be_1 & 0 & 0 \\
        \hyperlink{IV2}{\Romanbar{IV}.2_{t,x=1/2}} & \be_2 & 0 & 0 \\
        \hyperlink{IV4}{\Romanbar{IV}.4} & 0 & 0 & \be_1 \\
        \hyperlink{IV6}{\Romanbar{IV}.6} & \be_1 & \be_2 & \be_1 \\
        \hyperlink{IV6}{\Romanbar{IV}.6} & \be_1 & \be_2 & -\be_1 \\
        \hyperlink{IV6}{\Romanbar{IV}.6} & \be_1 & 0 & 0 \\
        \hyperlink{IV6}{\Romanbar{IV}.6} & \be_2 & 0 & 0 \\
        \bottomrule
    \end{tabular}\hspace{2cm}
    \rowcolors{2}{blue!10}{white}
    \begin{tabular}[t]{>{$}l<{$}*{3}{>{$}r<{$}}}
        \toprule
        \multicolumn{1}{l}{\textbf{Representation}} & \multicolumn{1}{c}{$\sigma_{11}$} & \multicolumn{1}{c}{$\sigma_{12}$} & \multicolumn{1}{c}{$\sigma_{22}$} \\
        \midrule[\heavyrulewidth]
        \hyperlink{V1}{\Romanbar{V}.1_{x=-2}} & 0 & 0 & \be_1 \\
        \hyperlink{V2}{\Romanbar{V}.2_{x=0}} & \be_1 & 0 & 0 \\
        \hyperlink{V2}{\Romanbar{V}.2_{x=1/2}} & \be_2 & 0 & 0 \\
        \hyperlink{V4}{\Romanbar{V}.4} & 0 & 0 & \be_1 \\
        \hyperlink{V6}{\Romanbar{V}.6} & \be_1 & \be_2 & \be_1 \\
        \hyperlink{V6}{\Romanbar{V}.6} & \be_1 & \be_2 & -\be_1 \\
        \hyperlink{V6}{\Romanbar{V}.6} & \be_1 & 0 & 0 \\
        \hyperlink{V6}{\Romanbar{V}.6} & \be_2 & 0 & 0 \\
        \hyperlink{VI1}{\Romanbar{VI}.1_{x=-2(1+c)}} & 0 & 0 & \be_1 \\
        \hyperlink{VI2}{\Romanbar{VI}.2_{x=2(1-c)}} & 0 & 0 & \be_1 \\
        \hyperlink{VI3}{\Romanbar{VI}.3_{c,t,x=0}} & \be_1 & 0 & 0 \\
        \hyperlink{VI3}{\Romanbar{VI}.3_{c,t,x=1/2}} & \be_2 & 0 & 0 \\
        \hyperlink{VI5}{\Romanbar{VI}.5} & 0 & 0 & \be_1 \\
        \hyperlink{VI7}{\Romanbar{VI}.7} & \be_1 & \be_2 & \be_1 \\
        \hyperlink{VI7}{\Romanbar{VI}.7} & \be_1 & \be_2 & -\be_1 \\
        \hyperlink{VI7}{\Romanbar{VI}.7} & \be_1 & 0 & 0 \\
        \hyperlink{VI7}{\Romanbar{VI}.7} & \be_2 & 0 & 0 \\
        \hyperlink{VII1}{\Romanbar{VII}.1_{c,t,x=0}} & \be_1 & 0 & 0 \\
        \hyperlink{VII1}{\Romanbar{VII}.1_{c,t,x=1/2}} & \be_2 & 0 & 0 \\
        \hyperlink{VII3}{\Romanbar{VII}.3} & 0 & 0 & \be_1 \\
        \hyperlink{VII5}{\Romanbar{VII}.5} & \be_1 & \be_2 & \be_1 \\
        \hyperlink{VII5}{\Romanbar{VII}.5} & \be_1 & \be_2 & -\be_1 \\
        \hyperlink{VII5}{\Romanbar{VII}.5} & \be_1 & 0 & 0 \\
        \hyperlink{VII5}{\Romanbar{VII}.5} & \be_2 & 0 & 0 \\
        \hyperref[sec:VIII2]{\Romanbar{VIII}.2} & \be_1 & \be_2 & \be_1 \\
        \hyperref[sec:VIII2]{\Romanbar{VIII}.2} & \be_1 & \be_2 & -\be_1 \\
        \hyperref[sec:VIII2]{\Romanbar{VIII}.2} & \be_1 & 0 & 0 \\
        \hyperref[sec:VIII2]{\Romanbar{VIII}.2} & \be_2 & 0 & 0 \\
        \hyperref[sec:B-IX]{\Romanbar{IX}} & \be_1 & \be_2 & \be_1 \\
        \hyperref[sec:B-IX]{\Romanbar{IX}} & \be_1 & \be_2 & -\be_1 \\
        \hyperref[sec:B-IX]{\Romanbar{IX}} & \be_1 & 0 & 0 \\
        \hyperref[sec:B-IX]{\Romanbar{IX}} & \be_2 & 0 & 0 \\
        \bottomrule
\end{tabular}
    \caption*{The Lie bracket in $V\otimes W$ is given by \\ $[u\otimes\be_i,v\otimes\be_j]=(u\times v)\otimes \sigma_{ij}$.}
\end{table}

\section*{Acknowledgements}

This project started during a visit of JML to the University of Edinburgh
in 2023. JML is grateful to said institution for its support during
the stay.

JMF would like to acknowledge the hospitality he received from the
Departamento de Matemáticas of the Universidade de Santiago de
Compostela and especially from José Carlos Díaz Ramos and the group of
Xeometría e Topoloxía during part of the time it took to write this
paper.

JML has been partially supported by the research grant
PID2022-138988NB-I00 funded by MICIU/AEI/10.13039/501100011033 and by
ERDF, EU, by project ED431C 2023/31 (Xunta de Galicia, Spain), and the
FPI programme (Spain).

\begin{appendices}
  
\section{Bianchi Lie algebras and their automorphisms}
\label{app:Bianchi}

In this appendix we collect some information on three-dimensional Lie
algebras and their automorphisms which we will use repeatedly in the
bulk of the paper.

\subsection{Bianchi Lie algebras}
\label{app:Bianchi-Lie-algebras}

The three-dimensional Lie algebras were classified up to isomorphism
by Bianchi \cite{Bianchi,MR1900159} and we list isomorphism class
representatives in Table~\ref{tab:Bianchi}.

\begin{table}[h!]
  \caption{Isomorphism classes of $3$-dimensional real Lie algebras.}
  \label{tab:Bianchi}
  \centering
  \rowcolors{2}{blue!10}{white}
  \begin{tabular}{*{4}{>{$}l<{$}}l}\toprule
    \multicolumn{1}{c}{Bianchi} &
    \multicolumn{3}{c}{Nonzero brackets} &
    \multicolumn{1}{c}{Structure}\\\midrule
    \hypertarget{BI}{\Romanbar{I}} & & & & abelian \\
    \hypertarget{BII}{\Romanbar{II}} & & & [B_2,B_3]=B_1 & nilpotent, unimodular \\
    \hypertarget{BIII}{\Romanbar{III}} \cong \Romanbar{VI}_1 & & & [B_2,B_3]=2 B_2 & solvable \\
    \hypertarget{BIV}{\Romanbar{IV}} & & [B_1,B_3]=B_1 &  [B_2,B_3]=B_1+B_2 & solvable\\
    \hypertarget{BV}{\Romanbar{V}} & & [B_1,B_3]=B_1 &  [B_2,B_3]=B_2 & solvable\\
    \hypertarget{BVI}{\Romanbar{VI}_{1\neq c\geq 0}} & & [B_1,B_3]=(c-1) B_1 & [B_2,B_3]=(c+1) B_2& solvable\\
    \hypertarget{BVII}{\Romanbar{VII}_{c\geq 0}} & & [B_1,B_3]=c B_1 - B_2 &  [B_2,B_3]=B_1 + c B_2 & solvable\\
    \hypertarget{BVIII}{\Romanbar{VIII}} &  [B_1,B_2]=-B_3 &  [B_1,B_3]=-B_2 &  [B_2,B_3]=B_1 & simple\\
    \hypertarget{BIX}{\Romanbar{IX}}  &   [B_1,B_2]=B_3 &  [B_1,B_3]=-B_2 &  [B_2,B_3]=B_1 & simple\\\bottomrule
  \end{tabular}
  \caption*{The columns are, from left to right, the classical Bianchi
  label for the Lie algebra, the nonzero brackets relative to a basis
  $B_1,B_2,B_3$ and some structural properties which can be
  more or less immediately inferred from the Lie brackets.}
\end{table}

\subsection{Automorphisms of Bianchi Lie algebras}
\label{app:Bianchi-autos}

For every $\fb$ in Table~\ref{tab:Bianchi}, we now list the
automorphisms $\Aut(\fb) \subseteq \GL(\fb)$.  By definition, $\varphi
\in \GL(\fb)$ is automorphism if for all $B,B'\in \fb$,
\begin{equation}
  [\varphi(B), \varphi(B')] = \varphi([B,B']).
\end{equation}
Automorphisms form a subgroup $\Aut(\fb) \subseteq \GL(\fb)$.  Choosing a
basis $\{B_1,B_2,B_3\}$ for $\fb$, every automorphism $\varphi$ gives
rise to a matrix $\Phi$ of $\GL(3,\RR)$ by
\begin{equation}
  \varphi(B_i) = B_j \Phi^j{}_i~.
\end{equation}
Since we have chosen a basis for $\fb$, we will give the
automorphism group $\Aut(\fb)$ as a matrix subgroup of $\GL(3,\RR)$.

\subsubsection{Automorphisms of Bianchi I}
\label{app:Bianchi-I-autos}

This is an abelian Lie algebra so that $\Aut(\fb) = \GL(\fb)$.
Choosing a basis $B_1,B_2,B_3$ for $\fb$, any matrix in $\GL(3,\RR)$
is an automorphism.

\subsubsection{Automorphisms of Bianchi II}
\label{app:Bianchi-II-autos}

The automorphism group $\Aut(\fb_{\scriptsize\Romanbar{II}})$ is the following
six-dimensional subgroup of $\GL(3,\RR)$:
\begin{equation}
  \Aut(\fb_{\scriptsize\Romanbar{II}}) = \left\{
    \begin{pmatrix}
      \det A & \ba^t \\ 0 & A
    \end{pmatrix} ~\middle| ~ A \in \GL(2,\RR)\text{ and }\ba \in \RR^2\right\}~.
\end{equation}

\subsubsection{Automorphisms of Bianchi III}
\label{app:Bianchi-III-autos}

The automorphism group $\Aut(\fb_{\scriptsize\Romanbar{III}})$ is the following
four-dimensional subgroup of $\GL(3,\RR)$:
\begin{equation}
  \Aut(\fb_{\scriptsize\Romanbar{III}}) = \left\{
    \begin{pmatrix}
      A & \ba \\ 0 & 1
    \end{pmatrix} ~\middle| ~ A \in \Diag(2,\RR)\text{ and } \ba \in \RR^2\right\}~,
\end{equation}
where $\Diag(2,\RR) \subseteq \GL(2,\RR)$ are the diagonal invertible
$2\times 2$ matrices.

\subsubsection{Automorphisms of Bianchi IV}
\label{app:Bianchi-IV-autos}

The automorphism group $\Aut(\fb_{\scriptsize\Romanbar{IV}})$ is the following
five-dimensional subgroup of $\GL(3,\RR)$:
\begin{equation}
  \Aut(\fb_{\scriptsize\Romanbar{IV}}) = \left\{
    \begin{pmatrix}
      A & \ba \\ 0 & 1
    \end{pmatrix}~ \middle| ~ A= \begin{pmatrix} u & v \\ 0 &
      u \end{pmatrix} \in \GL(2,\RR)\text{ and } \ba \in \RR^2\right\}~.
\end{equation}

\subsubsection{Automorphisms of Bianchi V}
\label{app:Bianchi-V-autos}

The automorphism group $\Aut(\fb_{\scriptsize\Romanbar{V}})$ is the following
six-dimensional subgroup of $\GL(3,\RR)$:
\begin{equation}
  \Aut(\fb_{\scriptsize\Romanbar{V}}) = \left\{
    \begin{pmatrix}
      A & \ba \\ 0 & 1
    \end{pmatrix} ~\middle|~ A \in \GL(2,\RR) \text{ and } \ba \in \RR^2\right\}~.
\end{equation}

\subsubsection{Automorphisms of Bianchi VI}
\label{app:Bianchi-VI-autos}

Here we do not consider $c=1$, which was already discussed as
Bianchi~$\Romanbar{III}$.  Nevertheless we must distinguish between
$c=0$ and $c>0$. In the generic case ($c> 0$), the automorphism
group $\Aut(\fb_{\scriptsize\Romanbar{VI}_{c> 0}})$ is the following
four-dimensional subgroup of $\GL(3,\RR)$:
\begin{equation}
  \Aut(\fb_{\scriptsize\Romanbar{VI}_{c> 0}}) = \left\{
    \begin{pmatrix}
      A & \ba \\ 0 & 1
    \end{pmatrix}~ \middle| ~ A \in \Diag(2,\RR) \text{ and } \ba \in \RR^2\right\}~.
\end{equation}
If $c=0$ then there are additional automorphisms of the form
\begin{equation}
  \left\{
    \begin{pmatrix}
    A & \ba \\ 0 & -1
  \end{pmatrix} ~\middle|~ A= \begin{pmatrix} 0 & u \\ v &
      0 \end{pmatrix} \in \GL(2,\RR)\text{ and } \ba \in \RR^2\right\}~.
\end{equation}

\subsubsection{Automorphisms of Bianchi VII}
\label{app:Bianchi-VII-autos}

Here again we must distinguish between $c=0$ and $c> 0$.  In the
generic case ($c> 0$), the automorphism group
$\Aut(\fb_{\scriptsize\Romanbar{VII}_{c> 0}})$ is the following four-dimensional
subgroup of $\GL(3,\RR)$:
\begin{equation}
\Aut(\fb_{\scriptsize\Romanbar{VII}_{c> 0}}) = \left\{
    \begin{pmatrix}
      A & \ba \\ 0 & 1
    \end{pmatrix}~ \middle| ~ A= \begin{pmatrix} u & v \\ -v &
      u \end{pmatrix} \in \GL(2,\RR)\text{ and } \ba \in \RR^2\right\}~.
\end{equation}
If $c=0$ then there are additional automorphisms of the form
\begin{equation}
  \left\{
    \begin{pmatrix}
    A & \ba \\ 0 & -1
  \end{pmatrix}~ \middle|~ A= \begin{pmatrix} u & v \\ v &
      -u \end{pmatrix} \in \GL(2,\RR)\text{ and } \ba \in \RR^2\right\}~.
\end{equation}

\subsubsection{Automorphisms of Bianchi VIII}
\label{app:Bianchi-VIII-autos}

The automorphism group of $\fb_{\scriptsize\Romanbar{VIII}} \cong
\fsl(2,\RR)$ is ${\mathbb P}\GL(2,\RR)$, which is the subgroup of
$\GL(3,\RR)$ obtained by conjugating $\fsl(2,\RR)$ (in our chosen
basis) by $\GL(2,\RR)$.

\subsubsection{Automorphisms of Bianchi IX}
\label{app:Bianchi-IX-autos}

The automorphism group of $\fb_{\scriptsize\Romanbar{IX}} \cong
\so(3)$ is the adjoint group $\mathsf{SO}(3) \subseteq \GL(3,\RR)$.

\section{An action of $\GL(2,\RR) \ltimes (\RR^2)^*$ on  $\Hom(\sS^{2}\RR^2,\RR^2)$}
\label{app:binary-cubics}

Consider the group $\mathsf{G}=\GL(2,\RR)\ltimes(\RR^2)^*$ and the
vector space $\Hom(\sS^2\RR^2,\RR^2)$.  The group $\mathsf{G}$ has an
action on $\Hom(\sS^2\RR^2,\RR^2)$ described as follows: if
$(h,\mu)\in\mathsf{G}$ and $\sigma\in\Hom(\sS^2\RR^2,\RR^2)$, we
define
\begin{equation}\label{eq:affine-action-on-sigma}
  ((h, \mu) \cdot \sigma)(\bx,\by) = h \sigma(h^{-1} \bx, h^{-1} \by) - \mu(h^{-1}\bx) \by - \mu(h^{-1}\by) \bx.
\end{equation}
for every $\bx,\,\by\in\RR^2$.  In this appendix, we are interested in
finding the possible normal forms of any
$\sigma\in\Hom(\sS^2\RR^2,\RR^2)$ under this action (in other words,
we want to give a unique representative for each $\mathsf{G}$-orbit).

Let us define a ``trace'' $\tr \colon \Hom(\sS^{2}\RR^2,\RR^2) \to (\RR^2)^*$ by
\begin{equation}
 (\tr \sigma)(\bx) = \Tr \hat\sigma_{\bx},
\end{equation}
for every $\bx \in \RR^2$, where $\hat\sigma_{\bx} \in \End(\RR^2)$ is the
endomorphism sending $\by \mapsto \sigma(\bx,\by)$ and $\Tr$ is the usual
trace of endomorphisms.  Relative to canonical dual bases
$\be_1, \be_2$ for $\RR^2$ and $\be^1,\be^2$ for $(\RR^2)^*$,
$\sigma(\be_i, \be_j) = \sigma_{ij}^k \be_k$ and the trace is given by
$\tr \sigma = \sigma_{ij}^j \be^i$.

\begin{lemma}
  For all $(h,\mu) \in \mathsf{G}$ and symmetric bilinear $\sigma :
  \RR^2 \times \RR^2 \to \RR^2$, we have
  \begin{equation}\label{eq:trace-transformation}
    \tr\left( (h,\mu) \cdot \sigma \right) = (h^{-1})^* ( \tr \sigma - 3 \mu).
  \end{equation}
\end{lemma}

\begin{proof}
  For all $\bx \in \RR^2$, we have an identity between endomorphisms
  \begin{equation}
  \widehat{\left( (h,\mu) \cdot \sigma\right)}_{\bx} = h \circ \hat\sigma_{h^{-1}\bx} \circ h^{-1} - \mu(h^{-1} \bx ) \id - \bx\otimes (\mu \circ h^{-1}),
  \end{equation}
  so that taking the trace
  \begin{align*}
      \tr \left( (h,\mu) \cdot \sigma\right)(\bx) &= \Tr \widehat{\left( (h,\mu) \cdot \sigma \right)}_{\bx}\\
                                              &= \Tr \left(h \circ \hat\sigma_{h^{-1}\bx} \circ h^{-1} - \mu(h^{-1} \bx ) \id - \bx \otimes (\mu \circ h^{-1})\right)\\
                                              &= \Tr \hat\sigma_{h^{-1} \bx} - 2 \mu(h^{-1} \bx) - \mu(h^{-1} \bx) & \tag{using linearity and invariance of $\Tr$}\\
                                              &= (\tr \sigma)(h^{-1}\bx) - 3 \mu(h^{-1}\bx)\\
                                              &= (h^{-1})^* (\tr \sigma - 3 \mu) (\bx),
  \end{align*}
which gives the desired result upon abstracting $\bx$.
\end{proof}

It follows from this lemma that, given any $\sigma$, we may act with
$(\id,\tfrac13 \tr \sigma)$ to arrive at a $\sigma'$ with
$\tr \sigma' =0$.  Furthermore, a direct consequence
of~\eqref{eq:trace-transformation} is that the stabiliser of
$U=\ker \tr$ is the canonical subgroup $\GL(2,\RR)$ (embedded in
$\mathsf{G}$ via the map $h\mapsto (h,0)$).  We conclude that the
$\mathsf{G}$-orbits in $\Hom(\sS^2\RR^2, \RR^2)$ are in a bijective
correspondence with the $\GL(2,\RR)$-orbits in $U$, a real
four-dimensional space with basis
\begin{equation}\label{eq:basis-U}
  \begin{split}
    \bsigma_1 &= \be_1 \otimes \be^1 \otimes \be^1 - \be_2 \otimes \be^1 \otimes \be^2 - \be_2 \otimes \be^2 \otimes \be^1\\
    \bsigma_2 &= \be_2 \otimes \be^2 \otimes \be^2 - \be_1 \otimes  \be^1 \otimes \be^2 - \be_1 \otimes \be^2 \otimes \be^1\\
    \bsigma_3 &= \be_1 \otimes \be^2 \otimes \be^2\\
    \bsigma_4 &= \be_2 \otimes \be^1 \otimes \be^1.
  \end{split}
\end{equation}

The subgroup $\SL(2,\RR) \subseteq \GL(2,\RR)$ preserves a symplectic
structure $\epsilon$ on $\RR^2$, which we normalise so that
$\epsilon(\be_i,\be_j) = \epsilon_{ij}$, the Levi-Civita symbol with
$\epsilon_{12}=1$.  Given $\sigma \in \Hom(\sS^2\RR^2,\RR^2)$ we may
contract with $\epsilon$ to give a cubic $q \in \Hom(\sS^2\RR^2
\otimes \RR^2,\RR)$ by
\begin{equation}
  q(\bx,\by,\bz) = \epsilon(\sigma(\bx,\by),\bz).
\end{equation}
This cubic form is symmetric in $\bx\leftrightarrow \by$, but more is
true.

\begin{lemma}
  If $\tr \sigma = 0$ then $q$ is totally symmetric.
\end{lemma}

\begin{proof}
  The condition $\tr \sigma =0$ says that the endomorphisms $\widehat
  \sigma_{\bx}$ belong to $\fsl(2,\RR)$ for every $\bx \in \RR^2$.
  This is equivalent to $\widehat \sigma_{\bx} (\by \wedge \bz) = 0$
  for all $\bx,\by,\bz \in \RR^2$.  Equivalently,
  \begin{equation}
    \sigma(\bx,\by) \wedge \bz + \by \wedge \sigma(\bx, \bz) = 0.
  \end{equation}
  Applying $\epsilon$, we find that $q(\bx,\by,\bz) =
  \q(\bx,\bz,\by)$.  Since $q(\bx,\by,\bz) = q(\by,\bx,\bz)$, this
  says that $q$ is totally symmetric.
\end{proof}

This means that to every traceless $\sigma$ we can associate a
binary cubic form $Q$ via
\begin{equation}
  Q(\bx) = \epsilon(\sigma(\bx,\bx),\bx),\qquad\text{for $\bx \in \RR^2$,}
\end{equation}
from where $q$ can be reconstructed by polarisation.  As a corollary,
the space $U$ of traceless symmetric bilinear maps $\sigma \colon
\RR^2 \times \RR^2 \to \RR^2$ is isomorphic to $\sS^3 (\RR^2)^*$ as a
representation of $\SL(2,\RR)$.

The action of $\GL(2,\RR)$ on $U$ is not the natural one, but is
instead twisted by the determinant, reflecting the fact that
$\epsilon$ is not invariant under $\GL(2,\RR)$.  Indeed, as a
representation of $\GL(2,\RR)$, $U \cong \sS^3 (\RR^2)^* \otimes
\ext{2} \RR^2$.
In classical notation (see, e.g.,
\cite[Definition~2.27]{MR1694364}), $Q$ is a binary cubic form of weight $-1$.
The determination of the normal forms of such $Q$ is a classical problem
in invariant theory and its solution is described, for example, in
\cite[Section~2]{MR1694364}.  Let us summarise the necessary results.

Let $Q$ be a cubic binary form of weight $-1$.  This means that, using
$\bx = (x,y)$ for our coordinates in $\RR^2$,
\begin{equation}
  Q(x,y) = a_3 x^3 + 3 a_2 x^2 y + 3 a_1 x y^2 + a_0 y^3
\end{equation}
and under a $\GL(2,\RR)$ transformation
\begin{equation}
  \label{eq:gl2r-transformation}
  \begin{pmatrix}
    x \\ y
  \end{pmatrix} \mapsto
  \begin{pmatrix}
    x' \\ y'
  \end{pmatrix} = 
  \begin{pmatrix}
    \alpha & \beta \\ \gamma & \delta
  \end{pmatrix}
  \begin{pmatrix}
    x \\ y
  \end{pmatrix},
\end{equation}
$Q \mapsto Q'$, where
\begin{equation}
  Q(x,y) = (\alpha\delta -\beta\gamma)^{-1} Q'(x',y').
\end{equation}

Recall that a homogeneous polynomial function
$J(\ba,\bx) := J(a_0,a_1,a_2,a_3,x,y)$ is said to be a \emph{covariant of $Q$ of
weight $k$} (written $\wt J = k$) if $J \mapsto J'$ under the
transformation~\eqref{eq:gl2r-transformation}, where
\begin{equation}
  J(\ba,\bx) = (\alpha\delta - \beta\gamma)^k \ J'(\ba', \bx').
\end{equation}

A covariant $J$ of a binary form $Q$ is subject to two natural
gradations in addition to its weight: its \emph{degree} $\deg J$ (the
degree as a polynomial in $\bx$) and its \emph{order} $\ord J$ (the
degree as a polynomial in $\ba$).  These are not independent, but
satisfy the relation:
\begin{equation}
  \label{eq:gradations-relation}
  \deg J + 2 \wt J = (\deg Q + 2 \wt Q) \ord J.
\end{equation}

A binary cubic form $Q$ has four independent covariants: $Q$ itself of
order $1$, degree $3$ and weight $-1$; the \emph{discriminant}
$\Delta$, defined by
\begin{equation}
  \Delta := a_0^2 a_3^2 - 6 a_0 a_1 a_2 a_3 + 4 a_0 a_2^3 - 3 a_1^2
  a_2^2 + 4 a_1^3 a_3,
\end{equation}
of degree $0$, order $4$ and weight $2$; the \emph{hessian} $H$
defined by
\begin{equation}
  H := Q_{xx} Q_{yy} - Q_{xy}^2,
\end{equation}
with subscripts denoting partial differentiation, of degree $2$, order
$2$ and weight $0$; and the \emph{jacobian} $T$ of $Q$ and $H$, defined by
\begin{equation}
  T := Q_x H_y - Q_y H_x,
\end{equation}
which we shall not need in our analysis.  The discriminant $\Delta$ of
a cubic $Q$ agrees up to a nonzero multiple with the discriminant of
the hessian $H$ as a binary quadratic form, so we do not get any
further covariants in this way.

With every binary cubic form $Q(x,y)$ there is an associated polynomial
$\hat Q (z) := Q(z,1)$, from where we may reconstruct $Q$ via
\begin{equation}
  Q(x,y) = y^3 \hat Q\left( \frac{x}{y} \right).
\end{equation}
The function $\hat Q$ is really to be understood as a function on the
projective space $\RR\PP^1$ on which $\GL(2,\RR)$ acts via fractional
linear transformations.

The nature of the roots of $\hat Q$ (the points in $\RR\PP^1$ where it
vanishes) is intimately linked to the possible values of the
covariants $\Delta$ and $H$.  Two basic facts are the following: the
discriminant $\Delta = 0$ if and only if $\hat Q$ has a repeated root
and the hessian $H=0$ if and only if $\hat Q$ has a triple root.

If $\Delta \neq 0$, then $\hat Q$ has three simple roots.  There are
two possibilities: either they are all real ($\Delta < 0$) or there
are two complex roots appearing as a complex conjugate pair ($\Delta >
0$).  If $\Delta = 0$, but $H \not\equiv 0$, then $\hat Q$ has a
double root.  If $H\equiv 0$ then $\hat Q$ has a triple root or $Q
\equiv 0$.

It is then a matter of studying the action of $\GL(2,\RR)$ via
fractional linear transformations on the roots of $\hat Q$, thought of
as points in $\RR\PP^1$ or $\CC\PP^1$.  This results in the following
normal forms:
\begin{itemize}
\item if $\Delta < 0$, then $\hat Q$ has three simple real roots which
  we can place at $-1,1,\infty$, resulting in $\hat Q (z) = z^2 -1$
  and hence
  \begin{equation}
    Q(x,y) = y x^2 - y^3;
  \end{equation}
\item if $\Delta > 0$, then $\hat Q$ has two complex conjugate roots
  and a real root, which we can place at $-i, i,\infty$, resulting in
  $\hat Q(z) = z^2 + 1$ and hence
  \begin{equation}
    Q(x,y) = y x^2 + y^3;
  \end{equation}
\item if $\Delta = 0$, but $H \not\equiv 0$, then $\hat Q$ has a
  repeated root, which we can place at $\infty$ while placing the simple
  root at $0$, resulting in $\hat Q(z) = z$ and hence
  \begin{equation}
    Q(x,y) = y^2x;
  \end{equation}
\item if $H \equiv 0$, but $Q \not\equiv 0$, we have a triple root,
  which we can place at $\infty$, resulting in $\hat Q(z)= 1$ and
  hence
  \begin{equation}
    Q(x,y) = y^3;
  \end{equation}
\item and finally we have the trivial case where $Q\equiv 0$.
\end{itemize}

Let us now apply this to the case at hand.  For the binary cubic form
arising from traceless symmetric bilinears $\sigma \colon \RR^2 \times
\RR^2 \to \RR^2$, we find for $\bx = (x,y)$,
\begin{equation}
  \label{eq:cubic-form-in-KLA}
  Q(\bx) = \epsilon(\sigma(\bx,\bx),\bx) = -\sigma_{11}^2 x^3 + 3
  \sigma_{11}^1 x^2 y + 3 \sigma_{12}^1 x y^2 + \sigma_{22}^1 y^3,
\end{equation}
where we have used the following parametrisation:
$\sigma = \sigma_{11}^1 \bsigma_1 - \sigma_{12}^1 \bsigma_2 +
\sigma_{22}^1 \bsigma_3 + \sigma_{11}^2 \bsigma_4$, where $\bsigma_i$
are given by equation~\eqref{eq:basis-U}.  The above normal forms for
$Q$ translate into the following:
\begin{itemize}
\item $Q(x,y)=x^2 y - y^3$ corresponds to
  \begin{equation}
    \sigma = \tfrac13 \bsigma_1 - \bsigma_3;
  \end{equation}
\item $Q(x,y) = x^2 y + y^3$ corresponds to
  \begin{equation}
    \sigma = \tfrac13 \bsigma_1 + \bsigma_3;
  \end{equation}
\item $Q(x,y) = x y^2$ corresponds to
  \begin{equation}
    \sigma = -\tfrac13 \bsigma_2;
  \end{equation}
\item $Q(x,y) = y^3$ corresponds to
  \begin{equation}
    \sigma = \bsigma_3;
  \end{equation}
\item and of course the trivial case $\sigma = 0$.
\end{itemize}

It is perhaps curious that the hessian $H$ of the binary cubic form $Q$ in
equation~\eqref{eq:cubic-form-in-KLA} is given by
\begin{equation}
    \tfrac{1}{18} H(x,y) = \alpha_{11} x^2 + 2 \alpha_{12} xy +  \alpha_{22} y^2,
\end{equation}
where $\alpha_{ij}$ are given by
equation~\eqref{eq:alpha-in-terms-of-sigma}.

\end{appendices}

\providecommand{\href}[2]{#2}\begingroup\raggedright\endgroup


\begin{thebibliography}{10}

\bibitem{Bacry:1968zf}
H.~Bacry and J.-M. L\'evy-Leblond, ``{Possible kinematics},''
  \href{http://dx.doi.org/10.1063/1.1664490}{{\em J. Math. Phys.} {\bfseries 9}
  (1968) 1605--1614}.

\bibitem{Bacry:1986pm}
H.~Bacry and J.~Nuyts, ``{Classification of Ten-dimensional Kinematical Groups
  With Space Isotropy},'' \href{http://dx.doi.org/10.1063/1.527306}{{\em J.
  Math. Phys.} {\bfseries 27} (1986) 2455}.

\bibitem{Figueroa-OFarrill:2017ycu}
J.~M. Figueroa-O'Farrill, ``{Kinematical Lie algebras via deformation
  theory},'' \href{http://dx.doi.org/10.1063/1.5016288}{{\em J. Math. Phys.}
  {\bfseries 59} no.~6, (2018) 061701},
  \href{http://arxiv.org/abs/1711.06111}{{\ttfamily arXiv:1711.06111
  [hep-th]}}.

\bibitem{Figueroa-OFarrill:2017tcy}
J.~M. Figueroa-O'Farrill, ``{Higher-dimensional kinematical Lie algebras via
  deformation theory},'' \href{http://dx.doi.org/10.1063/1.5016616}{{\em J.
  Math. Phys.} {\bfseries 59} no.~6, (2018) 061702},
  \href{http://arxiv.org/abs/1711.07363}{{\ttfamily arXiv:1711.07363
  [hep-th]}}.

\bibitem{Andrzejewski:2018gmz}
T.~Andrzejewski and J.~M. Figueroa-O'Farrill, ``{Kinematical lie algebras in 2
  + 1 dimensions},'' \href{http://dx.doi.org/10.1063/1.5025785}{{\em J. Math.
  Phys.} {\bfseries 59} no.~6, (2018) 061703},
  \href{http://arxiv.org/abs/1802.04048}{{\ttfamily arXiv:1802.04048
  [hep-th]}}.

\bibitem{Figueroa-OFarrill:2018ilb}
J.~Figueroa-O'Farrill and S.~Prohazka, ``{Spatially isotropic homogeneous
  spacetimes},'' \href{http://dx.doi.org/10.1007/JHEP01(2019)229}{{\em JHEP}
  {\bfseries 01} (2019) 229}, \href{http://arxiv.org/abs/1809.01224}{{\ttfamily
  arXiv:1809.01224 [hep-th]}}.

\bibitem{Figueroa-OFarrill:2019sex}
J.~Figueroa-O'Farrill, R.~Grassie, and S.~Prohazka, ``{Geometry and BMS Lie
  algebras of spatially isotropic homogeneous spacetimes},''
  \href{http://dx.doi.org/10.1007/JHEP08(2019)119}{{\em JHEP} {\bfseries 08}
  (2019) 119}, \href{http://arxiv.org/abs/1905.00034}{{\ttfamily
  arXiv:1905.00034 [hep-th]}}.

\bibitem{Morand:2023emw}
K.~Morand, ``{Possible ambient kinematics},''
  \href{http://dx.doi.org/10.1063/5.0159556}{{\em J. Math. Phys.} {\bfseries
  64} no.~11, (2023) 112502}, \href{http://arxiv.org/abs/2305.14371}{{\ttfamily
  arXiv:2305.14371 [hep-th]}}.

\bibitem{Figueroa-OFarrill:2022kcd}
J.~Figueroa-O'Farrill, R.~Grassie, and S.~Prohazka, ``{Lifshitz symmetry: Lie
  algebras, spacetimes and particles},''
  \href{http://dx.doi.org/10.21468/SciPostPhys.14.3.035}{{\em SciPost Phys.}
  {\bfseries 14} no.~3, (2023) 035},
  \href{http://arxiv.org/abs/2206.11806}{{\ttfamily arXiv:2206.11806
  [hep-th]}}.

\bibitem{BoostedLifshitz}
J.~Figueroa-O'Farrill, J.~Fluxman, R.~Grassie, and S.~Prohazka, ``{Homogeneous
  Lifshitz spacetimes with boosts}.'' In preparation, 2026.

\bibitem{MR155871}
G.~M. Mubarakzjanov, ``Classification of real structures of {L}ie algebras of
  fifth order,'' {\em Izv. Vys\v s. U\v cebn. Zaved. Matematika} {\bfseries
  1963} no.~3(34), (1963) 99--106.

\bibitem{Bianchi}
L.~Bianchi, ``Sugli spazi a tre dimensioni che ammettono un gruppo continuo di
  movimenti,'' {\em Memorie di Matematica e di Fisica della Societa Italiana
  delle Scienze, Serie Terza} {\bfseries Tomo XI} (1898) 267--352.

\bibitem{MR1900159}
L.~Bianchi, ``On the three-dimensional spaces which admit a continuous group of
  motions,'' \href{http://dx.doi.org/10.1023/A:1015357132699}{{\em Gen.
  Relativity Gravitation} {\bfseries 33} no.~12, (2001) 2171--2253}. Translated
  from the Italian by R. Jantzen.

\bibitem{MR442045}
J.~Patera and P.~Winternitz, ``Subalgebras of real three- and four-dimensional
  {L}ie algebras,'' \href{http://dx.doi.org/10.1063/1.523441}{{\em J.
  Mathematical Phys.} {\bfseries 18} no.~7, (1977) 1449--1455}.

\bibitem{MR2004893}
R.~O. Popovych, V.~M. Boyko, M.~O. Nesterenko, and M.~W. Lutfullin,
  ``Realizations of real low-dimensional {L}ie algebras,''
  \href{http://dx.doi.org/10.1088/0305-4470/36/26/309}{{\em J. Phys. A}
  {\bfseries 36} no.~26, (2003) 7337--7360}.

\bibitem{chapovskyi}
Y.~Y. Chapovskyi, S.~D. Koval, and O.~Zhur, ``{Subalgebras of Lie algebras.
  Example of $\mathfrak{sl}(3,\mathbb{R})$ revisited},''
  \href{http://dx.doi.org/10.1142/S0219498828500120}{{\em Journal of Algebra
  and Its Applications} (2024) },
  \href{http://arxiv.org/abs/2403.02554}{{\ttfamily arXiv:2403.02554
  [math-ph]}}. (to appear).

\bibitem{goursat}
E.~Goursat, ``Sur les substitutions orthogonales et les divisions
  r\'eguli\`eres de l'espace,''
  \href{http://dx.doi.org/10.24033/asens.317}{{\em Ann. Sci. \'Ecole Norm. Sup.
  (3)} {\bfseries 6} (1889) 9--102}.

\bibitem{Figueroa-OFarrill:2015bsi}
J.~Figueroa-O'Farrill and M.~Ungureanu, ``{Homogeneous M2 duals},''
  \href{http://dx.doi.org/10.1007/JHEP01(2016)150}{{\em JHEP} {\bfseries 01}
  (2016) 150}, \href{http://arxiv.org/abs/1511.03637}{{\ttfamily
  arXiv:1511.03637 [hep-th]}}.

\bibitem{MR142687}
G.~D. Mostow, ``On maximal subgroups of real {L}ie groups,''
  \href{http://dx.doi.org/10.2307/1970295}{{\em Ann. of Math. (2)} {\bfseries
  74} (1961) 503--517}. \url{https://doi.org/10.2307/1970295}.

\bibitem{MR1694364}
P.~J. Olver, \href{http://dx.doi.org/10.1017/CBO9780511623660}{{\em Classical
  invariant theory}}, vol.~44 of {\em London Mathematical Society Student
  Texts}.
\newblock Cambridge University Press, Cambridge, 1999.

\end{thebibliography}

\end{document}